\documentclass{article}
 \usepackage{latexsym}
\usepackage{amsmath}
 \usepackage{amsfonts}
\usepackage{graphicx}
\input amssym.def
\input amssym.tex

\newtheorem{theoreme}{Theorem }[section]
\newtheorem{proposition}[theoreme]{Proposition}
\newtheorem{lemma}[theoreme]{Lemma}
\newtheorem{definition}[theoreme]{Definition}

\newtheorem{remark}[theoreme]{Remark}
\newtheorem{example}[theoreme]{Example}

\newcommand{\beq}{\begin{equation}}
\newcommand{\eeq}{\end{equation}}
\newcommand{\ben}{\begin{enumerate}}
\newcommand{\een}{\end{enumerate}}

\newcommand{\bex}{\begin{example}}
\newcommand{\eex}{\end{example}}
\newcommand{\ber}{\begin{remark}}
\newcommand{\eer}{\end{remark}}

\def\bel{\begin{lemma}}
\def\eel{\end{lemma}}
\def\bet{\begin{theoreme}}
\def\eet{\end{theoreme}}
\def\bed{\begin{definition}}
\def\eed{\end{definition}}
\def\bep{\begin{proposition}}
\def\eep{\end{proposition}}

\def\rr{{\mathbb R}}
\def\zz{{\mathbb Z}}
\def\cc{{\mathbb C}}

\def\i{{\rm i}}

\def\id{{\rm id}}

\def\e{{\rm e}}
\def\d{{\rm d}}

\def\I{{\rm\scriptscriptstyle I}}
\def\II{{\rm\scriptscriptstyle {I}{I}}}
\def\0{{\rm\scriptscriptstyle 0}}

\def\cL{{\mathcal L}}

\def\C{{\cal C}}

\def\F{{\cal F}}

\def\S{{\cal S}}

\def\12{\frac{1}{2}}

\def\qed{\hfill$\Box$\medskip}
\def\proof{\noindent{\bf Proof.}\ \ }
\def\p{ \partial}

\newcounter{smallarabics}
\newenvironment{arabicenumerate}
{\begin{list}{{\normalfont\textrm{(\arabic{smallarabics})}}}
  {\usecounter{smallarabics}\setlength{\itemindent}{0cm}
   \setlength{\leftmargin}{5ex}\setlength{\labelwidth}{4ex}
   \setlength{\topsep}{0.75\parsep}\setlength{\partopsep}{0ex}
   \setlength{\itemsep}{0ex}}}
{\end{list}}

\newcounter{smallroman}
\newenvironment{romanenumerate}
{\begin{list}{{\normalfont\textrm{(\roman{smallroman})}}}
  {\usecounter{smallroman}\setlength{\itemindent}{0cm}
   \setlength{\leftmargin}{5ex}\setlength{\labelwidth}{4ex}
   \setlength{\topsep}{0.75\parsep}\setlength{\partopsep}{0ex}
   \setlength{\itemsep}{0ex}}}
{\end{list}}

\def\tF{{\tilde F}}

\def\cA{{\mathcal A}}
\def\Re{{\rm Re}}
\def\Im{{\rm Im}}

\begin{document}

\title{Hypergeometric type functions\\
and their symmetries\\
}
 \author{J. Derezi\'{n}ski
\\
Department of Mathematical Methods in Physics, \\
Faculty of Physics
\\ University of Warsaw,\\ Ho\.{z}a 74, 00-682 Warszawa, Poland,\\
email jan.derezinski@fuw.edu.pl}
\maketitle

\abstract{The paper is devoted to a systematic and unified discussion of various classes of hypergeometric type equations: the hypergeometric equation, the confluent equation, the $F_1$ equation (equivalent to the Bessel equation), the Gegenbauer equation and the Hermite equation. In particular, recurrence relations of their solutions, their integral representations and discrete symmetries are dicussed. }

\tableofcontents

\section{Introduction}
Following \cite{NU}, we adopt the following terminology. Equations of the form
\beq\left(\sigma(z)\p_z^2+\tau(z)\p_z+
\eta\right) f(z)=0,\ \ \
\label{req}\eeq
where
 $\sigma$ is a polynomial of degree $\leq2$,

\hspace{2.5ex} $\tau$ is a polynomial of degree $\leq1$,

\hspace{2.5ex} $\eta$ is a number,\\
will be  called {\em hypergeometric type equations}, and their solutions ---{\em hypergeometric type functions}. Differential operators
of the form $\sigma(z)\p_z^2+\tau(z)\p_z+
\eta$ will be called {\em hypergeometric type operators}.

The theory of hypergeometric type  functions is one of the oldest and most useful chapters of mathematics. In usual presentations it appears complicated and messy. The main purpose of this paper is an attempt to present its basics in a way that shows clearly its internal structure and beauty.

\subsection{Classification}

Let us start with a short review of basic classes of hypergeometric type equations. 
We will always  assume that $\sigma(z)\neq0$.
Every class, except for (9),
 will be simplified
 by dividing by a constant and an affine change of the complex variable $z$.

\begin{arabicenumerate}

\item
\noindent
{\bf The ${}_2F_1$
or  hypergeometric equation}

\[\left(z(1-z)\p_z^2+(c-(a+b+1)z)\p_z-ab\right)f(z)=0.\]

\item
\noindent
{\bf The ${}_2F_0$  equation}
\[\left(z^2\p_z^2+(1+(1+a+b)z)\p_z+ab\right)f(z)=0.\]

\item
\noindent
{\bf The ${}_1F_1$ 
  or  confluent equation}
 \[
(z\p_z^2+(c-z)\p_z-a)f(z)=0.\]

\item
\noindent
{\bf The ${}_0F_1$ equation
}
 \[
(z\p_z^2+c\p_z-1)f(z)=0.\]
\item
\noindent
{\bf The 
Gegenbauer equation}
\[\left((1-z^2)\p_z^2-(a+b+1)z\p_z-ab\right)f(z)=0.\]

\item
\noindent
{\bf The Hermite equation}
  \[
(\p_z^2-2 z\p_z-2a)f(z)=0.\]

\item
\noindent
{\bf 2nd order Euler equation}
\[
\left(z^2\p_z^2+bz\p_z+a\right)f(z)=0.\]

\item
\noindent
{\bf 1st order Euler equation for the derivative}
 \[
(z\p_z^2+c\p_z)f(z)=0.\]

\item
\noindent
{\bf  2nd order equation with constant coefficients}
 \[
(\p_z^2+c\p_z+a)f(z)=0.\]
\end{arabicenumerate}
 
One can divide these classes into 3 families:
\ben\item (1), (2), (3), (4);
\item (5), (6);
\item (7), (8), (9). 
\een

Each equation in the first family has a solution equal to the hypergeometric function ${}_pF_q$ with appropriate $p,q$. This function gives a name to the corresponding class of equations.

The second family consists of reflection invariant equations.
 
The third family consists of equations solvable in elementary functions. Therefore, it will not be considered in what follows.

 The ${}_2F_0$ and ${}_1F_1$ equation are equivalent by a simple substitution, therefore they can be discussed together.

Up to an affine
 transformation, (5) is 
a subclass of (1). However, it has additional properties, therefore it is useful to discuss it separately.

The main part of our paper  consists of 5 sections corresponding to the classes (1), (2)-(3), (4), (5) and (6).
 The discussion will  be  divided into two levels:
\ben \item Properties of the operator that defines the equation.
\item Properties of  functions solving the equation.
\een

\subsection{Properties of hypergeometric type operators}
\label{Properties of hypergeometric type operators}

We will discuss the following  types of properties of hypergeometric type operators:
\begin{romanenumerate}
\item equivalence between various classes,
\item integral representations of solutions,
\item discrete symmetries,
\item factorizations,
\item commutation relations.
\end{romanenumerate}

Let us give some examples of these properties. All these examples will be related to the ${}_1F_1$ equation.

We have
\begin{eqnarray}\label{cono1}
&&(-w)^{a+1}\left(w^2\p_w^2+(-1+(1+a+b)w)\p_w+ab\right)w^{-a}\\
 \ \ &=&z\p_z^2+(c-z)\p_z-a,\ \ \ \ \ \ \ w=-z^{-1}.\label{cono0}
\end{eqnarray}
Therefore the  ${}_1F_1$ operator, appearing in (\ref{cono0}), is equivalent to the
${}_2F_0$ operator, which is inside the brackets of (\ref{cono1}). This is an example of (i).

As an example of (ii) we quote the following fact: 
The  integral
\beq \int_\gamma t^{a-c}\e^t(t-z)^{-a}\d t\label{cono3}\eeq
is a solution of the ${}_1F_1$ equation provided that the values of the fuction \beq
t\mapsto 
t^{a-c+1}\e^t(t-z)^{-a-1}\label{cono2}\eeq
 at the endpoints of 
the curve $\gamma$ are equal to one another.

Note that the integrand of (\ref{cono3}) is an elementary function. The condition on the curve $\gamma$ can often be  satisfied in a number of non-equivalent ways, giving rise to distinct natural solutions.

An example of (iii) is the following identity:
\begin{eqnarray}\nonumber&&
 w\p_w^2+(c-w)\p_w-a\\&=&-\e^{-z}\left(z\p_z^2+(c-z)\p_z-c+a\right)\e^z,\ \ 
w=-z.\label{cono4}\end{eqnarray}
Thus the ${}_1F_{1}$ operator is transformed into a ${}_1F_{1}$ operator with different parameters.

Here is a pair of examples of (iv):
\begin{eqnarray}\nonumber
&&z(z\partial_z^2+(c-z)\p_z-a)\\\label{facto1}
&=&
\big(z\p_z+a-1\big)\big(z\p_z+c-a-z\big)+(a-1)(a-c)\\
\label{facto2}&=&
\big(z\p_z+c-a-1-z\big)\big(z\p_z+a\big)+a(a+1-c).
\end{eqnarray}

An example of (v) is
\begin{eqnarray}\nonumber&&
 \left(z\partial_z+a\right)z\left(z\p_z^2+(c-z)\p_z-a\right)\\&&=\ \ \ 
z\left(z\p_z^2+(c-z)\p_z-a-1\right)\left(z\partial_z+a\right).\label{cono5}
\end{eqnarray}
On both sides of the identity we see the ${}_1F_1$ operators whose parameters are contiguous.

The commutation properties can be derived from the factorizations. Let us
show, for example, how (\ref{facto1}) and  (\ref{facto2}) imply (\ref{cono5}).
First we  rewrite (\ref{facto1}) as
\begin{eqnarray}\nonumber
&&z(z\partial_z^2+(c-z)\p_z-a-1)\\\label{facto3}
&=&
\big(z\p_z+a\big)\big(z\p_z+c-a-1-z\big)+a(a+1-c).
\end{eqnarray}
Then we multiply (\ref{facto2}) from the left and 
 (\ref{facto3}) from the right by $\left(z\partial_z+a\right)$, obtaining identical right hand sides. This yields (\ref{cono5}).

\subsection{Hypergeometric type functions}

After the analysis of hypergeometric type operators, we discuss hypergeometric type functions, that is, functions annihilated by hypergeometric type operators. In particular, we will distinguish the so-called {\em standard solutions} which have a simple behavior around a singular point of the equation. In particular, if $z_0$ is a regular singular point, the Frobenius method gives us two solutions behaving as $ (z-z_0)^{\lambda_i}$, where $\lambda_1,\lambda_2$ are the  indices of $z_0$.
One can often find  solutions with a simple behavior also around
 irregular singular points. 

For reflection invariant classes (5) and (6) one can also define another pair of natural solutions: the even  solution $S^+$, which we normalize by $S^+(0)=1$, and the odd solution $S^-$, which we normalize by $(S^-)'(0)=2$.

Discrete symmetries can be used to derive properties of hypergeometric type functions. For instance,
(\ref{cono4}) implies that if $f(z)$ solves the confluent equation for parameters $c-a,c$, then so does $\e^zf(-z)$ for the parameters $a,c$. In particular, both functions $F(a;c;z)$ and $\e^{z}F(c-a;c;-z)$ solve the confluent equation for the parameters $a,c$. Both are analytic around $z=0$  and equal $1$ at $z=0$. By the uniqueness of the solution to the Frobenius method they should coincide. Hence we obtain the identity
\beq
F(a;c;z)=\e^zF(c-a;c;-z).\label{cono6}\eeq

Commutation relations are also useful. For example, it follows immediately from (\ref{cono5}) that $(z\partial_z+a)F(a;c;z)$ is a solution of the confluent equation for the parameters $a+1,c$. At zero it is analytic and its value is $a$. Hence we obtain the recurrence relation
\beq
(z\partial_z+a)F(a;c;z)=aF(a+1;c;z).\label{cono7}\eeq

For each class of equations we describe a whole family of recurrence relations. Every such a recurrence relation involves an operator of the following form:
a 1st order differential operator with no dependence on the parameters +
a multiplication operator depending linearly on the parameters.
We will call them {\em basic recurrence relations}. 

Sometimes there also exist  more complicated recurrence relations. We do not give their complete list, we only mention some of their examples. We call them {\em additional recurrence relations}.

Each of the standard solutions has  simple integral representations of the form analogous to  (\ref{cono3}). Each of these integral representations are associated to a pair of (possibly infinite and possibly coinciding) points where the integrand has a singularity. We will use two basic kinds of contours for standard solutions:
\begin{romanenumerate}\item[(a)] The contour starts at one singularity and ends at the other singularity; we assume that at both singularities the analog of (\ref{cono2})  is zero (hence, trivially, has equal values).
\item[(b)]
The contour starts at the first singularity, goes around the second singularity and returns to the first singularity; we assume that the analog of (\ref{cono2})  is zero at the first singularity.
\end{romanenumerate}
If available, we will always treat the type (a) contour as the basic one.

For instance, under appropriate conditions on the parameters, 
the ${}_1F_1$ function has the following  two  integral representations:
\begin{eqnarray}\hbox{type (a):} \ \ \ \ \ \ \ \ \  \int\limits_{[1,+\infty[}\e^{\frac{z}{t}}t^{-c}(t-1)^{c-a-1}\d t
&=&\frac{\Gamma(a)\Gamma(c-a)}{\Gamma(c)}F(a;c;z),\nonumber \\
\hbox{type (b):} \ \  \frac{1}{2\pi\i}\int\limits_{[1,0^+,1]}\e^{\frac{z}{t}}(-t)^{-c}(-t+1)^{c-a-1}\d t
&=& \frac{\Gamma(c-a)}{\Gamma(1-a)\Gamma(c)}F   (a;c;z).\nonumber
\end{eqnarray}
($0^+$ means that we bypass $0$ in the counterclockwise direction; in this case it is equivalent to bypassing $\infty$ in the clockwise direction).

There are various natural ways to normalize  hypergeometric type functions. 
The most obvious normalization for a  solution analytic at a given regular singular point is to demand that its value there is $1$. (For the ${}_2F_0$ equation, the point 0 is not regular singular, however there is a natural generalization of this normalization condition).
For equations (1)--(4), this function will be denoted by the letter $F$, consistently with the conventional usage. (Note the use of the italic font).
In the case of reflection symmetric equations
(5) and (6), we will use the letter $S$. 

However, it is often preferable to use  different normalizations, which involve   appropriate values of the Gamma funtion or its products. Such normalizations arise naturally when we consider integral representations. They will be denoted by ${\bf F}$ for equations (1) -- (4) (a similar notation  can be  found in \cite{NIST}), and ${\bf S}$ for (5) and (6). 
(Note the use of the boldface roman font).
Sometimes there will be  several varieties of these normalizations denoted by an appropriate superscript, related  to various integral representations. The functions with these normalizations  have often better properties than 
the $F$ and $S$ functions. This  is especially visible in recurrence relations, where the coefficient on the right (such as $a$ in (\ref{cono7})) depends on the normalization. 

For example, for the ${}_1F_1$ function we introduce the following normalizations:
\begin{eqnarray*}
{\bf F}(a;c;z)&:=&\frac{1}{\Gamma(c)}F(a;c;z),\\
{\bf F}^\I(a;c;z)&:=&\frac{\Gamma(a)\Gamma(c-a)}{\Gamma(c)}F(a;c;z),
\end{eqnarray*}
the latter suggested by the type (a) integral representation given above.

\subsection{Degenerate case}

For some values of parameters hypergeometric type functions have special properties. This happens in particular when the difference of the indices at a given regular singular point is an integer. Then the two standard solutions related to this point are proportional to one another. We call them {\em degenerate solutions}. (The best known example of such a situation are the Bessel functions of integer parameters).
In this case we  have a simple generating function and an additional integral representation, which involves integrating over a closed loop. 

\subsection{Canonical forms}
\label{s-canon}
Obviously,
  hypergeometric type operators coincide with differential operators of the form
\begin{eqnarray}\label{e1}
&&\sigma(z)\p_z^2+(\sigma'(z)+\kappa(z))\p_z+\12\kappa'+
\lambda\\&=&\p_z\sigma(z)\p_z+\12(\p_z\kappa(z)+\kappa(z)\p_z)+\lambda,\nonumber
\ \ \hbox{where}
\end{eqnarray}
 $\sigma$ is a polynomial of degree $\leq2$,\\$\kappa$ is a polynomial of degree $\leq1$,\\
$\lambda$ is a number.

One can argue that it is natural to use $\sigma,\kappa,\lambda$ to parametrize the hypergeometric
type operators (more natural than $\sigma,\tau,\eta$). (\ref{e1}) will be denoted
$\C(\sigma,\kappa,\lambda;z,\p_z)$, or, for brevity,
$\C(\sigma,\kappa,\lambda)$.
Let $\rho(z)$ be a solution of the equation
\beq(\sigma(z)\p_z-\kappa(z))\rho(z)=0.\label{e1b}\eeq
(Note that equation (\ref{e1b}) is solvable in elementary functions).
We have the identity
\beq
\C(\sigma,\kappa,\lambda)=
\rho^{-1}(z)\p_z\sigma(z)\rho(z)\p_z+\12\kappa'+\lambda,
\label{e1a}\eeq
We will call $\rho$ the {\em natural weight}. To justify this name note that if $\lambda$ is real, $\sigma,\kappa$ are real   and $\rho$ is positive and nonsingular on $]a,b[\subset \rr$, then $\C(\sigma,\kappa,\lambda)$ is Hermitian on the weighted space $L^2(]a,b[,\rho)$, when as the domain we take $C_{\rm c}^\infty(]a,b[)$.

It is sometimes useful to replace the   operator $\C(\sigma,\kappa,\lambda)$ with
\begin{eqnarray}\label{adva}
\rho(z)^{\frac12}\C(\sigma,\kappa,\lambda)\rho(z)^{-\frac12}
&=&
\partial_z\sigma(z)\partial_z-\frac{\kappa(z)^2}{4\sigma(z)}+\lambda.
\end{eqnarray}
We will call (\ref{adva}) the {\em balanced form of 
$\C(\sigma,\kappa,\lambda)$}.

Sometimes one replaces (\ref{req}) by the 1-dimensional Schr\"odinger equation
\beq \big(\partial_z^2-V(z)\big)f=0,\label{schro}\eeq
where
\begin{eqnarray*}
V(z):&=&
\frac12\big(\sigma(z)^{-1}\sigma'(z)\big)'
+\frac14\big(\sigma(z)^{-1}\sigma'(z)\big)^2+\frac{\kappa(z)^2}{4\sigma(z)^2}
-\frac{\lambda}{\sigma(z)}.
\end{eqnarray*}
(\ref{schro}) is equivalent to (\ref{req}), because
\begin{eqnarray}\label{adva1}
&&\sigma(z)^{-\frac12}\rho(z)^{\frac12}\C(\sigma,\kappa,\lambda)\rho(z)^{-\frac12}
\sigma(z)^{-\frac12}
=\partial_z^2-V(z),\end{eqnarray}
It   will be called the {\em Schr\"odinger-type form of the equation}
$\C(\sigma,\kappa,\lambda)f=0$.

  Some of the
  symmetries of hypergeometric type equations are obvious in
the balanced and Schr\"odinger-type forms
 these forms. This is
  partly due to the fact that they do not change when we switch the
  sign in front of $\kappa$.  This is a serious advantage of these forms.

In the literature various forms of hypergeometric type equations are
used. Instead of
the Gegenbauer equation one usually finds its balanced form, called
the {\em associated Legendre equation}. The {\em modified Bessel equation} and the {\em Bessel equation}, equivalent to
the rarely used ${}_0F_1$ equation, is the balanced form of a
special case of the ${}_1F_1$ equation. Instead of the ${}_1F_1$
equation one often finds its Schr\"odinger-type form, 
the {\em Whittaker equation}. This usage, due
 mostly to historical traditions, makes the subject more complicated
 than  necessary.

We will always use (\ref{req}) as the basic form. Its main advantage is
that
in almost all cases the equation in the form  (\ref{req}) has at least
one  solution analytic around a given finite singular point. Even in the case of the ${}_2F_0$
equation, whose all  solutions have a branch point at $0$, there exists a distinguished
solution particularly well behaved at zero.

\subsection{Hypergeometric type polynomials}
\label{Hypergeometric type polynomials}
{\em Hypergeometric type polynomials}, that is,
polynomial solutions of hypergeometric type equations deserve a separate analysis. They have traditional names involving various 19th century mathematicians. Note in particular that the (rarely used) polynomial cases of the ${}_2F_0$ function are called {\em Bessel polynomials}, however they do not have a direct relation to the better-known Bessel functions.

 There exists 
a well-known elegant approach to their theory that allows us to derive most of their basic properties in a unified way, see e.g. \cite{NU,R}. Let us sketch this approach.

Fix $\sigma,\kappa,\rho$, as in Subsect. \ref{s-canon}.
For any $n=0,1,2,\dots$ we
define
\beq
P_n(\sigma,\rho;z)
:=\frac{1}{n!}\rho^{-1}(z)\p_z^n\rho(z)\sigma^n(z).\label{rodrig}\eeq
We will call (\ref{rodrig}) a {\em Rodriguez-type formula}, since it is a generalization of the Rodriguez formula for Legendre polynomials.

One can show that  $P_n$ solves the equation
\beq
\Big(\sigma(z)\p_z^2+(\sigma'(z)+\kappa(z))\p_z
-n(n+1)\frac{\sigma''}{2}-n\kappa'\Big)P_n(\sigma,\rho;z)=0.
\label{pop1}\eeq
$P_n$ is a polynomial, typically of degree $n$, more precisely its
degree is given as follows:
\ben \item If $\sigma''=\kappa'=0$, then $\deg P_n=0$.
\item If $\sigma''\neq0$ and $-\frac{2\kappa'}{\sigma''}-1=m$ is a positive
  integer, then
 \[\deg P_n=\left\{\begin{array}{ll}
n,& n=0,1,\dots,m;\\n-m-1,&n=m+1,m+2,\dots.\end{array}\right.\] 
\item Otherwise, $\deg P_n=n$.
\een
We  have a generating function
\[\frac{\rho(z+t\sigma(z))}{\rho(z)}=\sum_{n=0}^\infty t^nP_n(\sigma,\rho\sigma^{-n};z),\]
 an integral representation
\beq
P_n(\sigma,\rho;z)=\frac{1}{2\pi \i}\rho^{-1}(z)\int\limits_{[z^+]}
\sigma^n(z+t)\rho(z+t)t^{-n-1}\d t
\label{rod}\eeq
and  recurrence relations
\begin{eqnarray*}
\bigl(\sigma(z)\partial_z+(\kappa(z)-n\sigma'(z)\big)
P_n(\sigma,\rho\sigma^{-n};z) &=&P_{n+1}(\sigma,\rho\sigma^{-n-1};z),\\
\partial_zP_{n+1}(\sigma,\rho\sigma^{-n-1};z)&=&\Big(-n\frac{\sigma''}{2}+\kappa'\Big)
P_n(\sigma,\rho\sigma^{-n};z).
\end{eqnarray*}

In almost all sections we devote a separate subsection to the
corresponding class of polynomials. Beside the properties that follow
immediately from the unified theory presented above we describe
additional properties valid in a given class. 

The   ${}_0F_1$ equation does not have polynomial solutions, hence the corresponding section is the only one without a subsection about polynomials.

Another special situation arises in the case of the Gegenbauer equation. The standard Gegenbauer polynomials found in the literature
do not have the normalization given by the Rodriguez-type formula. The Rodriguez-type formula yields the Jacobi polynomials, which for $\alpha=\beta$ coincide with the Gegenbaquer polynomials up to a nontrivial coefficient. Thus for the Gegenbauer equation it is natural to consider two classes of polynomials differing by normalization. This is related to an interesting symmetry called the Whipple transformation, which is responsible for two kinds of integral representations.

\subsection{Parametrization}

Each class  (1)--(6) depends on a number of complex parameters, denoted by Latin letters belonging to the set $\{a,b,c\}$. They will be called the {\em classical parameters}. They are  convenient when we discuss  power series expansions of standard solutions.

Unfortunately, the classical parameters are not convenient to describe discrete symmetries. Therefore, for each class (1)--(6) we introduce an alternative set of parameters, which we will call the  {\em Lie-algebraic parameters}. They will be denoted by Greek letters such as $\alpha,\beta,\mu,\theta,\lambda$, and will be given by certain linear (possibly, inhomogeneous) combinations of the classical parameters. Discrete symmetries of hypergeometric type  equations will simply involve  signed permutations of the Lie algebraic parameters -- in the classical parameters they look much more complicated. Recurrence relations also become simpler in the Lie-algebraic parameters.

For polynomials of hypergeometric type a third kind of parametrization is traditionally used. They are characterized by their degree $n$, which coincides  with $-a$, where $a$ is one of the classical parameters. The Lie-algebraic parameters appearing inside  the 1st order part of the equation are used as the remaining parameters.

Let us stress that all these parametrizations  are natural and useful. Therefore, we sometimes face the dilemma which parametrization to use for a given set of identities. We usually try to choose the one that  gives the simplest formulas.

We sum up the information about various parametrizations in the following table:
\[\begin{array}{cccccc}
\hbox{Equation}& \begin{array}{c}\hbox{classical}\\\hbox{parameters}\end{array}
&\begin{array}{c}\hbox{Lie-algebraic}\\\hbox{parameters}\end{array}&\hbox{Polynomial}&
\begin{array}{c}\hbox{parameters}\\\hbox{for polynomials}\end{array}\\[2ex]\hline
\\[2ex]
{}_2F_1&
a,b,c&\begin{array}{c}\alpha=c-1\\
\beta=a+b-c\\
\gamma=b-a\end{array}
&\hbox{Jacobi}
&\begin{array}{c}\alpha=c-1\\
\beta=a+b-c\\
n=-a\end{array}
\\[5ex]
{}_2F_0&
a,b&\begin{array}{c}\theta=-1+a+b\\\alpha=a-b
\end{array}
&\hbox{Bessel}
&\begin{array}{c}
\theta=-1+a+b\\
n=-a\end{array}
\\[5ex]
{}_1F_1&
a,c&\begin{array}{c}\theta=-c+2a\\
\alpha=c-1\end{array}
&\hbox{Laguerre}
&\begin{array}{c}
\alpha=c-1\\n=-a\end{array}\\[5ex]
{}_0F_1&
c&
\alpha=c-1
&-----
&-----\\[5ex]
\hbox{Gegenbauer}&
a,b&\begin{array}{c}
\alpha=\frac{a+b-1}{2}\\
\lambda=\frac{b-a}{2}
\end{array}
&\begin{array}{c}
\hbox{$\alpha=\beta$ Jacobi}\\\hbox{or Gegenbauer}\end{array}
&\begin{array}{c}
\alpha=\frac{a+b-1}{2}\\
n=-a\end{array}\\[5ex]
\hbox{Hermite}&
a&\lambda =a-\frac12
&\hbox{Hermite}
&n=-a
\end{array}\]

\subsection{Group-theoretical background}

Identities for hypergeometric type operators and functions 
 have  a high degree of symmetry. 
Therefore, it is natural to expect that a certain group-theoretical structure is responsible for these identities.

There exists a large literature about the relations between special functions and the group theory \cite{V,Wa,M1, VK}. Nevertheless, as far as we know, the arguments found in the literature give a rather incomplete explanation of the properties   that we describe. In a seperate publication \cite{DM} we would like to present a group-theoretical approach to hypergeometric type functions with, we believe, a  more satisfactory  justification of their high symmetry.
Below we would like to briefly sketch the main ideas of \cite{DM}. 

Each  hypergeometric type equation can be obtained by separating  the variables of a certain 2nd order PDE of the complex variable with constant coefficients. One can introduce the  Lie algebra of generalized symmetries of this PDE. In this Lie algebra we fix
 a certain maximal commutative  algebra, which we will call the ``Cartan algebra''. Operators whose adjoint action is diagonal in the ``Cartan algebra'' will be called ``root operators''. 
Automorphisms of the Lie algebra leaving invariant the ``Cartan algebra'' will be called ``Weyl symmetries''.

(Note that in some cases the Lie algebra of symmetries is simple, and then the names {\em Cartan algebra}, {\em root operators} amd {\em Weyl symmetries} correspond to the standard names. In other cases the Lie algebra is non-semisimple, and then the names are less standard -- this is the reason for the quotation marks that we use).

Now the parameters of hypergeometric type equation  can be interpreted as the eigenvalues of elements of the ``Cartan algebra''. In particular, the Lie agebraic parameters correspond to a certain natural choice of the ``Cartan algebra''.
 Each  recurrence relation  is related to a ``root operator''.
Finally, each symmetry of a hypergeometric type operator corresponds to a Weyl symmetry of the Lie algebra.

We can distinguish 3 kinds of PDE's with constant coefficients:
\ben \item The  {\em Helmholtz equation} on $\cc^n$ given by $\Delta_n+1$, whose Lie algebra of symmetries is $\cc^n\rtimes so(n,\cc)$;
\item The {\em Laplace equation} on $\cc^n$ given by $\Delta_n$, whose Lie algebra of  generalized symmetries is $so(n+2,\cc)$
\item The  {\em heat equation} on $\cc^n\oplus\cc$ given by $\Delta_n+\partial_s$, whose Lie algebra of generalized symmetries is $sch(n,\cc)$ (the so-called {\em (complex) Schr\"odinger Lie algebra}.\een
 Separating the variables in these equations usually leads to differential equations with many variables. Only in a few cases it leads to ordinary differential equations, which turn out to be of hypergeometric type. Here is a table of these cases:
\[\begin{array}{ccccc}
\hbox{PDE}&\begin{array}{c}\hbox{Lie}\\ \hbox{algebra}\end{array}
&\begin{array}{c}\hbox{dimension of}\\ \hbox{Cartan algebra}\end{array}
&\begin{array}{c}\hbox{discrete}\\ \hbox{symmetries}\end{array}&
\hbox{equation}\\[1ex]
\hline\\[1ex]
\Delta_2+1&\cc^2\rtimes so(2,\cc)&1&\zz_2&{}_0F_1;
\\[1.5ex]
\Delta_4&so(6,\cc)&3&\hbox{cube}& {}_2F_1; \\[1.5ex]
\Delta_3&so(5,\cc)&2&\hbox{square}&\hbox{Gegenbauer};\\[1.5ex]
\Delta_2+\p_s&sch(2,\cc)&2&\zz_2\times
\zz_2&
{}_1F_1\hbox{ or }{}_2F_0;\\[1.5ex]
\Delta_1+\p_s&sch(1,\cc)&1&\zz_4&
\hbox{Hermite}.
\end{array}\]


\subsection{Comparison with the literature}

There exist many works  that discuss hypergeometric type functions, e.g.
\cite{NIST,Ho,MOS,AAR,R,WW,Ol,Tr}. Some of them are meant to be encyclopedic collections of formulas, other try to show   mathematical structure that underlies their properties. 

In our opinion, this work differs substantially from the existing literature.
In our presentation we try to follow the intrinsic logic of the subject, without too much regard for  the  traditions.
If possible, we apply the same pattern  to each class of hypergeometric type equations.
This  sometimes
 forces us to introduce unconventional notation.

We believe that the intricacy of usual presentations of hypergeometric type functions can be partly explained by historical reasons. In the literature various classes of these functions are often described with help of different conventions. 
Sometimes we  will give short remarks devoted to the conventions found
in the literature. These remarks will  always be clearly separated from
the main text.

Of course, our presentation does not contain all useful identities and properties of hypergeometric functions. Some of them are on purpose left out, e.g. the so-called addition formulas. We restrict ourselves to what we view as the  most basic theory. On the other hand, we try to be complete for each type of properties that we consider.

Our work is  strongly inspired by the book by Nikiforov and Uvarov
\cite{NU}, who tried to develop a unified approach to hypergeometric
type functions. They stressed in particular the role of integral
representations and of recurrence relations. 

Another important influence are
 the works of Miller \cite{M1,M2} who stressed
the Lie-algebraic structure behind the recurrence relations.

The method of factorization can be traced back  at least to \cite{IH}.

\medskip

\noindent{\small{\bf Acknowledgement.} I acknowledge the help of Laurent
 Bruneau, Micha\l{} Godli\'{n}ski, and especially Micha\l{} Wrochna and Przemys\l{}aw Majewski who
  proofread parts of previous versions of this work.

The research of the author was supported in part by the National
Science
 Center (NCN) grant No.  2011/01/B/ST1/04929.}

\section{Preliminaries}

In this section we fix basic terminology, notation and collect a number of well known useful facts, mostly from complex analysis. It is supposed to serve as a reference and can be skipped at the first reading.

\subsection{Differential equations}

The main object of our paper are ordinary homogeneous 2nd order
linear differential
equations in the complex domain, that is equations of the form
\beq \left(a(z)\partial^2_z+b(z)\partial_z+c(z)\right)\phi(z)=0.\label{equo}\eeq
 It will be convenient to treat (\ref{equo}) as the problem of
finding the kernel of the operator
\beq\cA(z,\partial_z):=
a(z)\partial^2_z+b(z)\partial_z+c(z).\label{equo2}\eeq
We will then say that {\em the equation (\ref{equo}) is given by the operator (\ref{equo2})}. When we do not consider the change of the variable, we will often write $\cA$ for $\cA(z,\partial_z)$.

\subsection{The principal branch of the logarithm and the power function}
\label{a.1}
The function 
\[\{z\in\cc\ :\ -\pi< \Im z<\pi\}\ni z\mapsto \e^z\in\cc\backslash]-\infty,0]\]
is bijective. Its inverse will be called the {\em principal branch of the logarithm} 
and will be denoted simply $\log z$.

If $\mu\in\cc$ then the {\em principal branch of the power function} is defined
as
\[\cc\backslash]-\infty,0]\ni z\mapsto z^\mu:=\e^{\mu\log z}.\]

Consequently, if $\alpha\in\cc\backslash\{0\}$, then
the functions $\log(\alpha(z-z_0))$ and $(\alpha(z-z_0))^\mu$
have the domain
$\cc\backslash(z_0+\alpha^{-1}]-\infty,0])$.

Of course, if needed we will use the analytic continuation to extend the
definition of the logarithm and  the power function beyond 
$\cc\backslash]-\infty,0]$ onto the appropriate covering of
    $\cc\backslash\{0\}$. 

\subsection{Contours}
\label{a.3}
We will write 
\[f(z)\Big|_{z_0}^{z_1}:=f(z_1)-f(z_0).\]
In particular, if $]0,1[\ni t\mapsto \gamma(t)\in\cc$ is a curve, then
\beq f(z)\Big|_{\gamma(0)}^{\gamma(1)}
=\int_\gamma f'(z)\d z.\label{mainth}\eeq

In order to avoid making pictures,
we will use special notation for contours of integration.

Broken lines will be denoted as in the following example:
\[[w_0,u,w_1]:=[w_0,u]\cup[u,w_1].\]

\begin{center}
\includegraphics[width=8cm,totalheight=2cm]{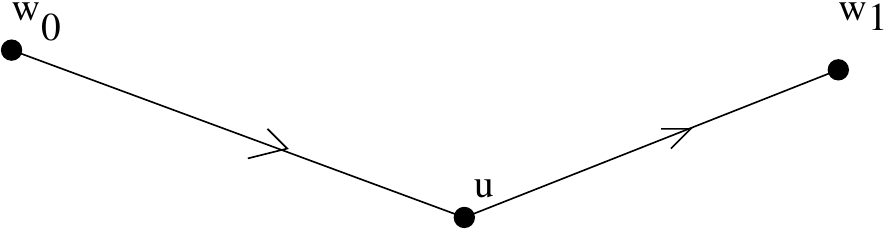}
\end{center}

This contour may be inappropriate if the function has a nonintegrable
singularity at $u$. Then we might want
to bypass $u$ with a small arc counterclockwise or clockwise. In such a case we can use the
curves
\begin{eqnarray}&&\label{zn3}
[w_0,u^+,w_1].
\end{eqnarray}
\begin{center}
\includegraphics[width=12cm,totalheight=1.6cm]{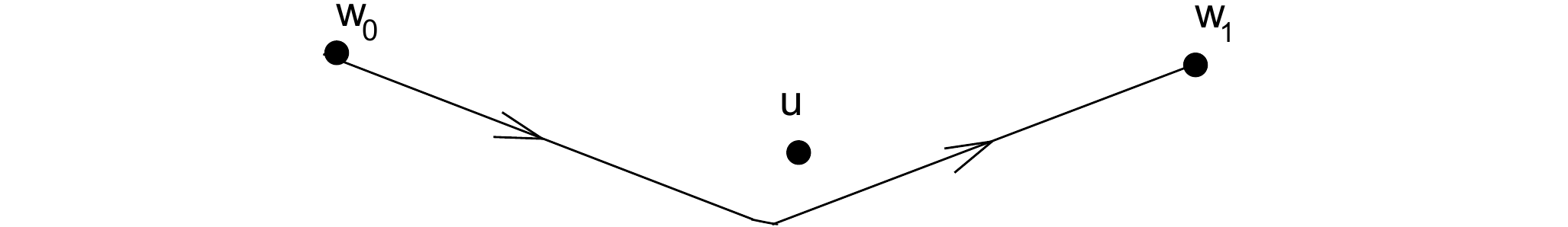}
\end{center}

\begin{eqnarray}\label{z4}
&&[w_0,u^-,w_1].
\end{eqnarray}
\begin{center}
\includegraphics[width=12cm,totalheight=1.4cm]{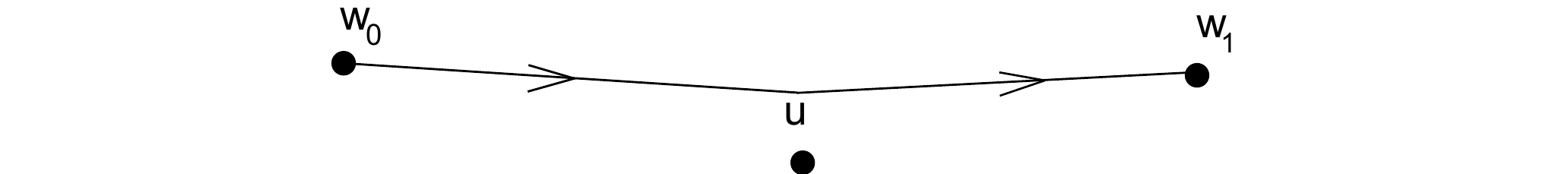}
\end{center}

We may want to bypass a group of points, say $u_1,u_2$. Such contours are denoted by
\[[w_0,(u_0,u_1)^+,w_1],\]
\begin{center}
\includegraphics[width=12cm,totalheight=2cm]{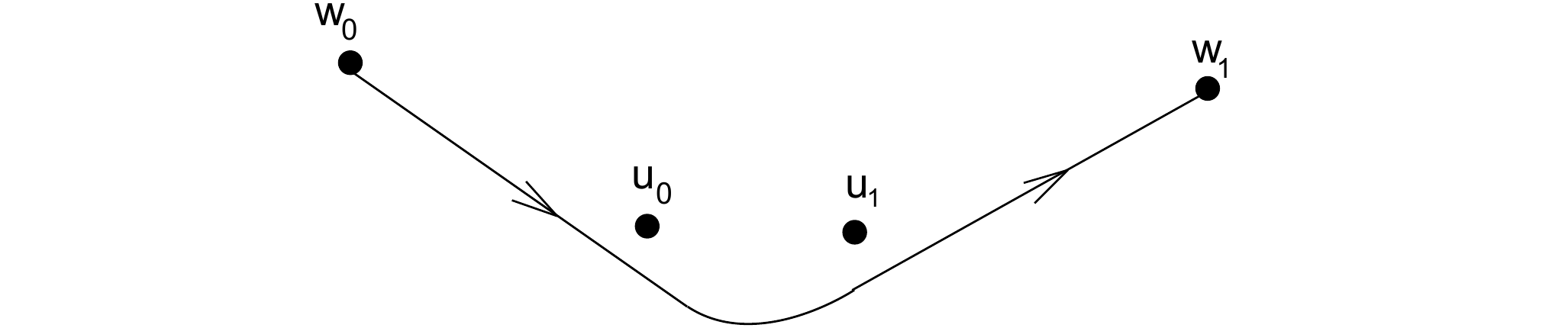}
\end{center}

\[[w_0,(u_0,u_1)^-,w_1].\]

\begin{center}
\includegraphics[width=12cm,totalheight=2cm]{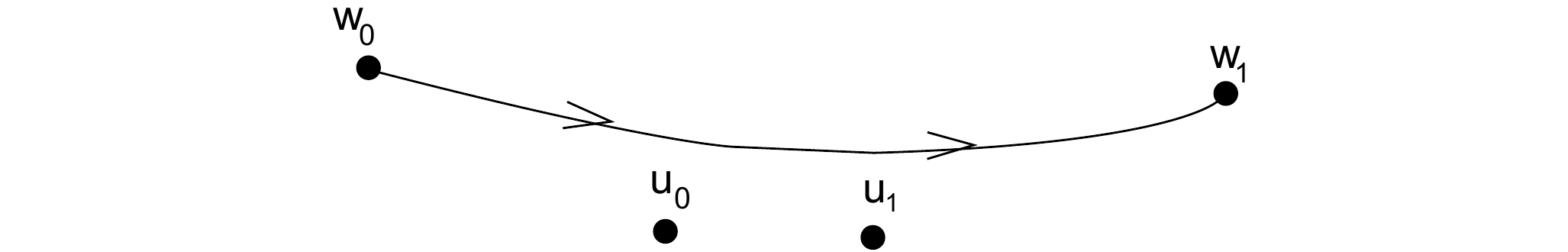}
\end{center}

A small counterclockwise/clockwise loop around $u$ is denoted 
\[[u^+],\hskip 26ex [u^-]\]
\begin{center}
\includegraphics[width=7cm,totalheight=1cm]{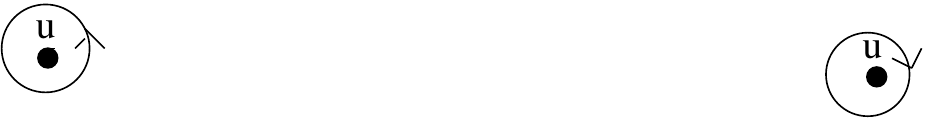}
\end{center}

A counterclockwise/clockwise loop around a group of points, say, $u_1,u_2$ is denoted
\[[(u_1,u_2)^+],\hskip 16ex [(u_1,u_2)^-].\]
\begin{center}
\includegraphics[width=12cm,totalheight=2cm]{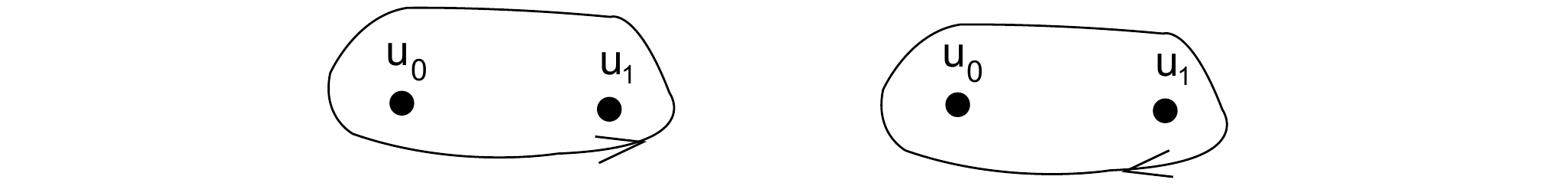}
\end{center}

A half-line starting at $u$ and inclined at the angle $\phi$ is denoted
\beq[u,\e^{\i\phi}\infty[:=\{u+\e^{\i\phi}t\ :\ t>0\}:\eeq
\vskip 0.2cm
\begin{center}
\includegraphics[width=12cm,totalheight=2cm]{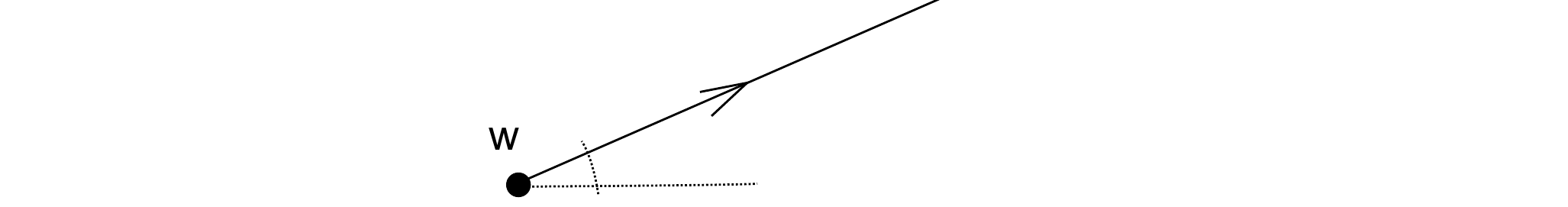}
\end{center}

We will also need  slightly more complicated contours:
\begin{eqnarray}
&&[(u+\e^{\i\phi}\cdot 0)^+,w]
\nonumber
\end{eqnarray}
\begin{center}
\includegraphics[width=12cm,totalheight=1cm]{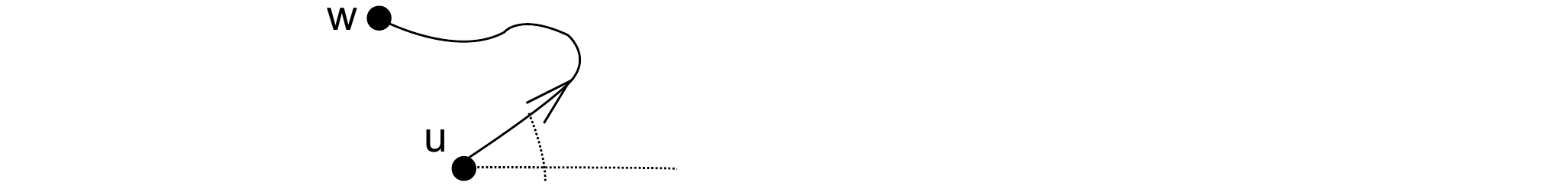}
\end{center}
Here,
the contour departs from  $u$ at the angle $\phi$, then it bypasses $u$
with a small arc counterclockwise
and then it goes in the direction of $w$.

The following countour has the shape of a kidney:
\begin{eqnarray}&&
[(u+\e^{\i\phi}\cdot 0)^+]
\nonumber\end{eqnarray}
\begin{center}
\includegraphics[width=12cm,totalheight=2cm]{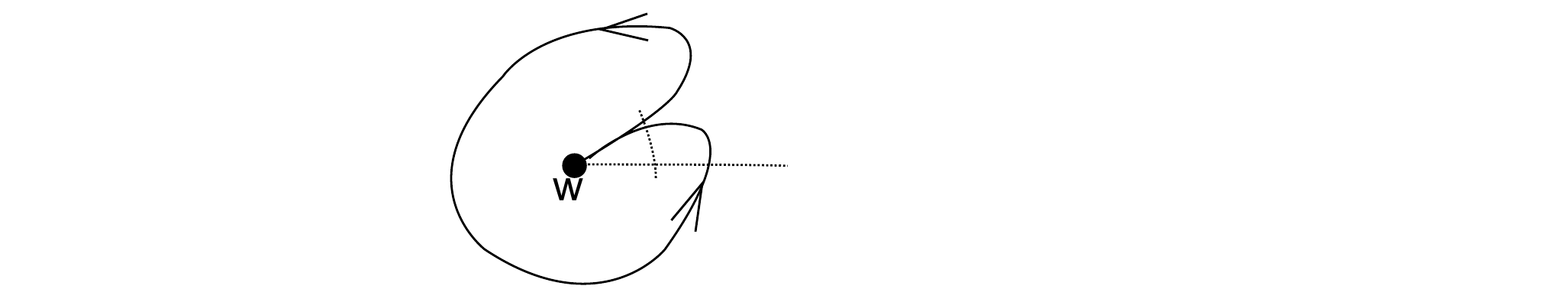}
\end{center}
This contour  departs from  $u$ at the angle $\phi$, then it
goes around  $u$ and  returns to $u$ again at the angle $\phi$.

Instead of $u+\e^{\i0}\cdot 0$ we will write $u+0$. Likewise, instead of
$u+\e^{\i\pi}\cdot 0$ we will write $u-0$.

\subsection{Reflection invariant  differential equations}
\label{a.5}

Consider a 2nd order differential operator
\beq\p_z^2+b(z)\p_z+c(z).\label{e5}\eeq
Assume that (\ref{e5}) 
is invariant w.r.t. the reflection $z\mapsto -z$. This means that
for some functions $\pi$, $\rho$ we have
\[b(z)=z\pi(z^2),\ \ \ \ c(z)=\rho(z^2).\]
Then it is natural to make a quadratic change of coordinates:
\begin{eqnarray}\nonumber
&&\p_z^2+b(z)\p_z+c(z)\\&=&
4u\left(\p_u^2+\Big(\frac{1}{2u}+\frac{\pi(u)}{2}\Big)\p_u+\frac{\rho(u)}{4u}\right),\label{refl1}\\[2ex]\nonumber
&&z^{-1}(\p_z^2+b(z)\p_z+c(z))z\\&=&
4u\left(\p_u^2+\Big(\frac{3}{2u}+\frac{\pi(u)}{2}\Big)\p_u++\frac{\pi(u)+\rho(u)}{4u}\right),\label{refl2}\end{eqnarray}
where 
\[u=z^2,\ \ \ z=\sqrt u.\]
Thus if $g_+(u)$, resp. $g_-(u)$ satisfy
\begin{eqnarray*}
\left(\p_u^2+\Big(\frac{1}{2u}+\frac{\pi(u)}{2}\Big)\p_u+\frac{\rho(u)}{4u}\right)
g_+(u)&=&0,\\
\left(\p_u^2+\Big(\frac{3}{2u}+\frac{\pi(u)}{2}\Big)\p_u+\frac{\pi(u)+\rho(u)}{4u}\right)g_-(u)&=&0,\end{eqnarray*}
then $g_+(z^2)$ is an even solution, resp. $zg_-(z^2)$ is an odd solution of the equation given by   (\ref{e5}).

Note that if $\pi,\rho$ are holomorphic, then $0$ is a regular singular point of  (\ref{refl1}) with indices $0,\12$ and of 
 (\ref{refl2}) with indices $0,-\12$.

\subsection{Regular singular points}
\label{a4}
 
 In this subsection we recall well known facts about regular singular points
of differential equations

We will write
\[f(z)\sim (z-z_0)^\lambda\ \ \ \ \hbox{at}\ \ \ z_0\]
if $f(z)(z-z_0)^{-\lambda}$ is analytic at $z_0$ and
$\lim\limits_{z\to z_0}f(z)(z-z_0)^{-\lambda}=1$.
In particular, we write
\[f(z)\sim 1\ \ \ \ \hbox{at}\ \ \ z_0\]
if $f$ is analytic in a neighborhood of $z_0$ and $f(z_0)=1$.

An equation given by the operator
\beq
\p_z^2+b(z)\p_z+c(z)
\label{e9}\eeq with meromorphic coefficients $a(z)$,
$c(z)$  has a {\em regular singular point at $z_0$} if
\[b_0:=\lim_{z\to z_0}b(z)(z-z_0),\ \ \ \ 
c_0:=\lim_{z\to z_0}c(z)(z-z_0)^2\]
exist. The {\em indices $\lambda_1$, $\lambda_2$ of $z_0$} are the solutions of the
indicial equation
\[\lambda(\lambda-1)+b_0\lambda+c_0=0.\]
 
\bet[The Frobenius method] If
$\lambda_1-\lambda_2\neq-1,-2,\cdots$, then
 there exists a unique solution $f(z)$ of the equation given by (\ref{e9}) such that 
$f(z)\sim(z-z_0)^{\lambda_1}$ at $z_0$. 
\label{fro}\eet

The case
$\lambda_1-\lambda_2\in\zz$ is called the {\em degenerate case.} In this case the Frobenius method gives one solution corresponding to the point $z_0$.

Likewise, (\ref{e9}) has a
 {\em regular singular point at $\infty$} if
\[\tilde b_0:=\lim_{z\to \infty}b(z)z,\ \ \ \ 
\tilde c_0:=\lim_{z\to \infty}c(z)z^2\]
exist. The {\em indices $\tilde\lambda_1$, $\tilde\lambda_2$ of $\infty$}  are the solutions of the
indicial equation
\[\tilde\lambda(\tilde\lambda+1)-\tilde b_0\tilde\lambda+\tilde c_0=0.\]
\bet[The Frobenius method at infinity] If 
$-\tilde\lambda_1+\tilde\lambda_2\neq-1,-2,\cdots$, then
 there exists a unique solution $\tilde f_1(z)$ of (\ref{e9}) such that 
$\tilde f_1(z)\sim z^{-\tilde\lambda_1}$ at $\infty$.
\eet

Note the identity
\begin{eqnarray}\label{shifto}&&
(z-z_0)^{-\theta}\left(\p_z^2+b(z)\p_z+c(z)\right)(z-z_0)^\theta\\
&=&\partial_z^2+\big(2\theta(z-z_0)^{-1}+b(z)\big)\partial_z+
(\theta^2-\theta)(z-z_0)^{-2}+\theta b(z)(z-z_0)^{-1}+c(z).
\nonumber\end{eqnarray}
If $z_0$ is a regular singular point, then  the corresponding  indices of (\ref{shifto})
equal those of (\ref{e9}) $+\theta$. Likewise, if $\infty$ is a regular sigular point, then  the corresponding indices are shifted by $-\theta$. The indices corresponding to other points are left unchanged.

\subsection{The Gamma function}

In this section we collect basic identities related to {\em Euler's Gamma function} that we will use.


\begin{eqnarray}{\bf Relationship\ to\ factorial}&&\Gamma(n+1)=n!,\ \ n=0,1,2,\dots,\label{g5a}\\
{\bf Recurrence\ \  relation}&&\Gamma(z+1)=z\Gamma(z),\label{g5}\\
{\bf Reflection\ \  formula}&&
 \Gamma(z)\Gamma(1-z)=\frac{\pi}{\sin\pi z},\\
{\bf II \   \ Euler's\ \ integral.}&&
\Gamma(z):=\int_0^\infty\e^{-t}t^{z-1}\d t,\ \ \ \Re z>0
\label{gam}
\\
{\bf Hankel's\ \  formula.}
&&
\frac{1}{\Gamma(-z+1)}=\frac{1}{2\pi \i}\!\!\!\!\!\!\!\!
\int\limits_{[-\infty,0^+,-\infty[}\!\!\!\!\!\!\!\e^{t}t^{z-1}\d t,\\
{\bf Legendre's\  \  formula}&&
2^{2z-1}\Gamma(z)\Gamma\left(z+1\slash2\right)=
\sqrt{\pi}\Gamma(2z).
\label{g17aa}\end{eqnarray}
{\bf I  Euler's integral and its consequences.}\begin{eqnarray}
\frac{\Gamma(u)\Gamma(v)}{\Gamma(u+v)}
&=&\int_0^1 t^{u-1}(1-t)^{v-1}\d t
\ \ \Re u>0,\ \Re v>0,\label{b0}
\\[3ex]
\frac{\Gamma(u)\Gamma(v)}{\Gamma(u+v)}\frac{\sin\pi
u}{\sin\pi(u+v)}&&\label{b1}\\
=\
\frac{\Gamma(1-u-v)\Gamma(v)}{\Gamma(1-u)}&
=&\int_1^\infty t^{u-1}(t-1)^{v-1}\d t
\ \ \Re v>0,\ \Re (1-u-v)>0,
\nonumber\\[3ex]
\frac{\Gamma(-u)\Gamma(-v)}{\Gamma(-u-v)}\frac{\sin\pi
u}{\sin\pi(u+v)}&&\label{b2}\\
=\ \frac{\Gamma(-u-v+1)}
{\Gamma(-u+1)\Gamma(-v+1)}&=&
\frac{1}{2\pi \i}
\int\limits_{]-\infty,0^+,-\infty[}t^{u-1}(1-t)^{v-1}\d t\nonumber\\&=&
\frac{1}{2\pi \i}
\int\limits_{]\infty,1^-,\infty[}t^{u-1}(1-t)^{v-1}\d t,\ \ \ \Re(-u-v+1)>0.\nonumber
\\[3ex]
\frac{\Gamma(u)\Gamma(v)}{\Gamma(u+v)}\sin\pi u&&\label{put}\\
=\ \frac{\Gamma(v)}{\Gamma(1-u)\Gamma(u+v)}&=&
\frac{1}{2\pi \i}
\int\limits_{]1,0^+,1]}t^{u-1}(1-t)^{v-1}\d t,\ \ \ \Re v>0.\nonumber
\\[3ex]
\frac{\Gamma(u)\sqrt\pi}{\Gamma(u+\frac12)}&=&
\int_{-1}^1(1-s^2)^{u-1}\d s,\\[3ex]
\frac{\Gamma(u)\sqrt\pi}{2\cos\pi
  u\Gamma(u+\frac12)}
&=&\int_1^\infty(s^2-1)^{u-1}\d s.\end{eqnarray}



\subsection{The Pochhammer symbol}
\label{a.2}
If $a\in\cc$ and $n\in\zz$, then the so-called {\em Pochhammer symbol}
is defined as follows:
\[\begin{array}{ll}
(a)_0=1,\\[3MM]
(a)_n:=a(a+1)\dots(a+n-1),&\ n=1,2,\dots\\[3mm]
(a)_n:=\frac{1}{(a-n)\dots(a-1)},&\ n=\dots,-2,-1.
\end{array}\]
Note the identities
\begin{eqnarray}
&&(a)_n=\frac{\Gamma(a+n)}{\Gamma(a)}=
(-1)^n\frac{\Gamma(1-a)}{\Gamma(1-a-n)}=(-1)^n(1-n-a)_n,\nonumber\\
&&(1-z)^{-a}=\sum\limits_{n=0}^\infty\frac{(a)_n}{n!}z^n,\ \ \ |z|<1,\label{double1}\\
&&(1/2)_nn!=\frac{(2n)!}{2^{2n}},\ \ \ (3/2)_nn!=\frac{(2n+1)!}{2^{2n}}.\label{double}
\end{eqnarray}

\section{The ${}_2F_1$ or the hypergeometric equation}
\label{s3}
\subsection{Introduction}
Let $a,b,c\in\cc$.
Traditionally, the {\em hypergeometric equation} is given by the operator
\beq
\F(a,b;c;z,\p_z):=z(1-z)\p_z^2+\big(c-(a+b+1)z\big)\p_z-ab.\label{hy1-tra}\eeq

The {
\em classical parameters $a,b,c$}  will be often replaced by another set of parameters  $\alpha,\beta,\mu\in\cc$, called 
{\em Lie-algebraic}.  They are related to one another by
\[\begin{array}{rll}
\alpha:=c-1,&\ \ \beta: =a+b-c,&\ \ \mu:=b-a;\\[2ex]
\label{newnot}
a=\frac{1+\alpha+\beta -\mu}{2},&\ \ b=\frac{1+\alpha+\beta +\mu}{2},&\ \ c=1+\alpha.
\end{array}\]
In the Lie-algebraic parameters the hypergeometric operator (\ref{hy1-tra}) becomes
\begin{eqnarray}&&
\F_{\alpha,\beta ,\mu}(z,\p_z)\label{hy1}\\
&=&z(1-z)\p_z^2+
\big((1+\alpha)(1-z)-(1+\beta )z\big)\p_z+\frac14\mu^2
-\frac14(\alpha+\beta +1)^2.\nonumber\end{eqnarray}
The Lie-algebraic parameters  have an interesting interpretation in terms of the natural basis of the Cartan algebra of the Lie algebra $so(6)$ \cite{DM}.

The singular points of the hypergeometric operator are located at $0,1,\infty$. All of them are regular singular. The indices of these points are
\begin{center}
\begin{tabular}{ccc}
$z=0$&$z=1$&$z=\infty$\\\hline\\
$1-c=-\alpha$ & $c-a-b=-\beta$ &$a=\frac{1+\alpha+\beta-\mu}{2}$\\
$0$&$0$&$ b=\frac{1+\alpha+\beta+\mu}{2}$
\end{tabular} 
\end{center}
Thus the Lie-algebraic parameters are the differences of the indices.

The hypergeometric operator remains the same if we interchange $a$  and $b$ (replace $\mu$ with $-\mu$).

\subsection{Integral representations}

\bet
Let $[0,1]\ni  t\mapsto\gamma(t)$ satisfy
\[t^{b-c+1}(1-t)^{c-a}(t-z)^{-b-1}\Big|_{\gamma(0)}^{\gamma(1)}=0.\]
Then
\beq 
\F(a,b;c;z,\p_z)
\int_\gamma t^{b-c}(1-t)^{c-a-1}(t-z)^{-b}\d t=0.
\label{f4}\eeq
\label{intr}\eet

\proof 
We check that for any contour $\gamma$ (\ref{f4}) equals
\[-b \int_\gamma
\Big(\p_t t^{b-c+1}(1-t)^{c-a}(t-z)^{-b-1}\Big)\d t
.\]
\qed

Analogous (and nonequivalent) 
integral representations can be obtained by interchanging $a$ and $b$
in Theorem \ref{intr}.

\subsection{Symmetries}
\label{s-sym}

To every permutation of the set
of singularities   $\{0,1,\infty\}$ we can associate
exactly one homography
 $z\mapsto w(z)$.
Using the method described at the end of Subsect. \ref{a4}, with every such a homography we can associate 4 substitutions 
 that preserve the 
form of the hypergeometric equation. 
Altogether there are $6\times4=24$ substitutions. They form a group isomorphic
to the group of proper symmetries of the cube. If we  take into account
the fact that replacing $\mu$ with $-\mu$ is also an obvious symmetry of the
hypergeometric equation, then we obtain a group of $2\times 24=48$
elements, isomorphic to the group of all (proper and improper) symmetries of
a cube, which is the Weyl group of $so(6)$.
 
Below we describe  the table of
 symmetries of the hypergeometric operator
 except for those obtained by switching the sign of the last
 parameter.
We fix the sign of the last parameter by demanding that the number of minus signs is even.
 
Note that the table looks much simpler in the Lie-algebraic parameters than in the classical parameters.

All the operators below equal $\F_{\alpha,\beta ,\mu}(w,\p_w)$ for the
corresponding $w$:
\[\begin{array}{rrcl}
w=z:&&&\\[1ex]
&&\F_{\alpha,\beta ,\mu}(z,\p_z),&
\\[1ex]
&(-z)^{-\alpha}(z-1)^{-\beta }&\F_{-\alpha,-\beta ,\mu}(z,\p_z)&(-z)^{\alpha}(z-1)^{\beta }
\\[1ex]
&(z-1)^{-\beta }&\F_{\alpha,-\beta ,-\mu}(z,\p_z)&(z-1)^{\beta },\\[1ex]
&(-z)^{-\alpha}&\F_{-\alpha,\beta ,-\mu}(z,\p_z)&
(-z)^{\alpha};\\
w=1-z:&&&\\
&&\F_{\beta ,\alpha,\mu}(z,\p_z),&
\\[1ex]
&(z-1)^{-\alpha}(-z)^{-\beta }&\F_{-\beta ,-\alpha,\mu}(z,\p_z)&(z-1)^{\alpha}(-z)^{\beta },
\\[1ex]
&(z-1)^{-\alpha}&\F_{\beta ,-\alpha,-\mu}(z,\p_z)&
(z-1)^{\alpha},\\[1ex]
&(-z)^{-\beta }&\F_{-\beta ,\alpha,-\mu}(z,\p_z)
&(-z)^{\beta };\\
w=\frac{1}{z}:&&&\\
&(-z)^{\12(\alpha+\beta +\mu +1)}&(-z)\F_{\mu ,\beta ,\alpha}(z,\p_z)&
(-z)^{\12(-\alpha-\beta -\mu -1)},\\[1ex]
&(-z)^{\12(\alpha+\beta -\mu +1)}(z-1)^{-\beta }&
(-z)\F_{-\mu ,-\beta ,\alpha}(z,\p_z) &(-z)^{\12(-\alpha-\beta +\mu -1)}(z-1)^{\beta },
\\[1ex]
&(-z)^{\12(\alpha+\beta +\mu +1)}(z-1)^{-\beta }&
(-z)
\F_{\mu ,-\beta ,-\alpha}(z,\p_z)& (-z)^{\12(-\alpha-\beta -\mu -1)}(z-1)^{\beta }
,\\[1ex]
&(-z)^{\12(\alpha+\beta -\mu +1)}&(-z)\F_{-\mu ,\beta ,-\alpha}(z,\p_z)
& (-z)^{\12(-\alpha-\beta +\mu -1)}
;\\
w=1-\frac{1}{z}:&&&\\
&(-z)^{\12(\alpha+\beta +\mu +1)}&
(-z)\F_{\mu ,\alpha,\beta }(z,\p_z)&(-z)^{\12(-\alpha-\beta -\mu -1)},\\[1ex]
&(-z)^{\12(\alpha+\beta -\mu +1)}(z-1)^{-\alpha}&
(-z)\F_{-\mu ,-\alpha,\beta }(z,\p_z) &(-z)^{\12(-\alpha-\beta +\mu -1)}(z-1)^{\alpha},\\[1ex]
&(-z)^{\12(\alpha+\beta +\mu +1)}(z-1)^{-\alpha}&
(-z) \F_{\mu ,-\alpha,-\beta }(z,\p_z)& (-z)^{\12(-\alpha-\beta -\mu -1)}(z-1)^{\alpha},
\\[1ex]
 &(-z)^{\12(\alpha+\beta -\mu +1)}&(-z)\F_{-\mu ,\alpha,-\beta }(z,\p_z)&
  (-z)^{\12(-\alpha-\beta +\mu -1)};\\
w=\frac{1}{1-z}:&&&\\
& (z-1)^{\12(\alpha+\beta +\mu +1)}&
(z-1)\F_{\beta ,\mu ,\alpha}(z,\p_z)& (z-1)^{\12(-\alpha-\beta -\mu -1)},\\[1ex]
&  (-z)^{-\beta }(z-1)^{\12(\alpha+\beta -\mu +1)}
&(z-1)\F_{-\beta ,-\mu ,\alpha}(z,\p_z) & (-z)^{\beta }(z-1)^{\12(-\alpha-\beta +\mu -1)},\\[1ex]
& (z-1)^{\12(\alpha+\beta -\mu +1)}
&(z-1)\F_{\beta ,-\mu ,-\alpha}(z,\p_z) & (z-1)^{\12(-\alpha-\beta +\mu -1)},\\[1ex]
&(-z)^{-\beta }(z-1)^{\12(\alpha+\beta +\mu +1)}&(z-1)
\F_{-\beta ,\mu ,-\alpha}(z,\p_z)&
 (-z)^{\beta }(z-1)^{\12(-\alpha-\beta -\mu -1)};
\\
w=\frac{z}{z-1}:&&\\
&(z-1)^{\12(\alpha+\beta +\mu+1)}&(z-1)\F_{\alpha,\mu,\beta }(z,\p_z)&
(z-1)^{\12(-\alpha-\beta -\mu-1)}
,\\[1ex]
&  (-z)^{-\alpha}(z-1)^{\12(\alpha+\beta -\mu+1)}&(z-1)\F_{-\alpha,-\mu,\beta }(z,\p_z) &
 (-z)^{\alpha}(z-1)^{\12(-\alpha-\beta +\mu-1)}
,\\[1ex]
& (z-1)^{\12(\alpha+\beta -\mu+1)}&
(z-1)\F_{\alpha,-\mu,-\beta }(z,\p_z) & (z-1)^{\12(-\alpha-\beta +\mu-1)}
,\\[1ex]
& (-z)^{-\alpha}(z-1)^{\12(\alpha+\beta +\mu +1)}
&(z-1)\F_{-\alpha,\mu ,-\beta }(z,\p_z)&
  (-z)^{\alpha}(z-1)^{\12(-\alpha-\beta -\mu -1)}.
\end{array}\]

\subsection{Factorization and commutation relations}
\label{commu}

The   hypergeometric operator can be factorized in several ways:
\begin{eqnarray*}
\F_{\alpha,\beta,\mu}&=&
\Big(z(1-z)\partial_z+\big((1+\alpha)(1-z)-(1+\beta)  z\big)\Big)\partial_z\\&&-\frac14(\alpha+\beta+\mu+1)(\alpha+\beta-\mu+1),\\
&=&
\partial_z\Big(z(1-z)\partial_z+\big(\alpha(1-z)-\beta  z\big)\Big)\\&&-\frac14(\alpha+\beta+\mu-1)(\alpha+\beta-\mu-1),\\
&=&
\Big((1-z)\partial_z-\beta-1\Big)\Big(z\partial_z+\alpha\Big)\\
&&-\frac14(\alpha+\beta+\mu+1)(\alpha+\beta-\mu+1),\\
&=&
\Big(z\partial_z+\alpha+1\Big)\Big((1-z)\partial_z-\beta\Big)\\
&&-\frac14(\alpha+\beta+\mu+1)(\alpha+\beta-\mu+1);
\end{eqnarray*}\begin{eqnarray*}
z\F_{\alpha,\beta,\mu}&=&
\Big(z\partial_z+\frac12(\alpha+\beta+\mu-1)\Big)
\Big(z(1-z)\partial_z+\frac12(1-z)(\alpha+\beta-\mu+1)-\beta\Big)\\
&&-\frac14(\alpha+\beta+\mu-1)(\alpha-\beta-\mu+1),\\
&=&
\Big(z(1-z)\partial_z+\frac12(1-z)(\alpha+\beta-\mu+1)-\beta-1\Big)
\Big(z\partial_z+\frac12(\alpha+\beta+\mu+1)\Big)
\\
&&-\frac14(\alpha+\beta+\mu+1)(\alpha-\beta-\mu-1),\\
&=&
\Big(z\partial_z+\frac12(\alpha+\beta-\mu-1)\Big)
\Big(z(1-z)\partial_z+\frac12(1-z)(\alpha+\beta+\mu+1)-\beta\Big)\\
&&-\frac14(\alpha+\beta-\mu-1)(\alpha-\beta+\mu+1),\\
&=&
\Big(z(1-z)\partial_z+\frac12(1-z)(\alpha+\beta+\mu+1)-\beta-1\Big)
\Big(z\partial_z+\frac12(\alpha+\beta-\mu+1)\Big)
\\
&&-\frac14(\alpha+\beta-\mu+1)(\alpha-\beta+\mu-1);
\end{eqnarray*}\begin{eqnarray*}
(z-1)\F_{\alpha,\beta,\mu}&=&
\Big((z-1)\partial_z+\frac12(\alpha+\beta+\mu-1)\Big)
\Big(z(1-z)\partial_z+\frac12z(-\alpha-\beta+\mu-1)+\alpha\Big)\\
&&-\frac14(\alpha+\beta+\mu-1)(\alpha-\beta+\mu-1),\\
&=&
\Big(z(1-z)\partial_z+\frac12z(-\alpha-\beta+\mu-1)+\alpha+1\Big)
\Big((z-1)\partial_z+\frac12(\alpha+\beta+\mu+1)\Big)
\\
&&-\frac14(\alpha+\beta+\mu+1)(\alpha-\beta+\mu+1),\\
&=&
\Big((z-1)\partial_z+\frac12(\alpha+\beta-\mu-1)\Big)
\Big(z(1-z)\partial_z+\frac12z(-\alpha-\beta-\mu-1)+\alpha\Big)\\
&&-\frac14(\alpha+\beta-\mu-1)(\alpha-\beta-\mu-1),\\
&=&
\Big(z(1-z)\partial_z+\frac12z(-\alpha-\beta-\mu-1)+\alpha+1\Big)
\Big((z-1)\partial_z+\frac12(\alpha+\beta-\mu+1)\Big)
\\
&&-\frac14(\alpha+\beta-\mu+1)(\alpha-\beta-\mu+1).
\end{eqnarray*}
One  way of showing the above factorizations is as follows: We
start with deriving  the first one,
 and then we apply the symmetries of Subsect. \ref{s-sym}.

 The factorizations can be used to derive the following commutation relations:
\[\begin{array}{rrl}
&
 \p_z&\F_{\alpha,\beta ,\mu }\\[1ex]
&=\ \ \ \ \ \ \F_{\alpha+1,\beta +1,\mu } &\p_z,\\[3ex]
&(z(1-z)\p_z+(1-z)\alpha-z\beta )&\F_{\alpha,\beta ,\mu }\\[1ex]
&=\ \ \ \ \ \ \F_{\alpha-1,\beta -1,\mu }&(z(1-z)\p_z+(1-z)\alpha-z\beta ),\\[3ex]
& ((1-z)\p_z -\beta )&\F_{\alpha,\beta ,\mu }\\[1ex]
&=\ \ \ \ \ \ \F_{\alpha+1,\beta -1,\mu }& ((1-z)\p_z -\beta ),
\\[3ex]&(z\p_z+\alpha)&\F_{\alpha,\beta ,\mu }
\\[1ex]
&=\ \ \ \ \ \ \F_{\alpha-1,\beta +1,\mu }&(z\p_z+\alpha);
\end{array}\]
\[\begin{array}{rrl}
&(z\p_z+\12 (\alpha+ \beta +\mu +1))&z\F_{\alpha,\beta ,\mu }
\\[1ex]
&=\ \ \ \ \ \ z\F_{\alpha,\beta +1,\mu +1}&(z\p_z+\12 (\alpha+ \beta +\mu +1)),
\\[3ex]
&(z(1{-}z)\p_z{+}\frac12(1{-}z)(\alpha{+}\beta {-}\mu {+}1){-}\beta )&z\F_{\alpha,\beta ,\mu }\\[1ex]
&=\ \ \ \ \ \ z\F_{\alpha,\beta {-}1,\mu {-}1}&
(z(1{-}z)\p_z{+}\frac12(1{-}z)(\alpha{+}\beta {-}\mu {+}1){-}\beta ),\\[3ex]
&( z\p_z{+}\12(\alpha{+}\beta {-}\mu {+}1))&z\F_{\alpha,\beta ,\mu }\\[1ex]
&=\ \ \ \ \ \ z\F_{\alpha,\beta {+}1,\mu {-}1}&( z\p_z{+}\12(\alpha{+}\beta {-}\mu {+}1),
\\[3ex]
&
 (z(z{-}1)\p_z{-}\12(1{-}z)(\alpha{+}\beta {+}\mu {+}1)
{+}\beta )&z\F_{\alpha,\beta ,\mu }\\[1ex]
&=\ \ \ \ \ \ z\F_{\alpha,\beta -1,\mu +1}&(z(z{-}1)\p_z{-}\12(1{-}z)(\alpha{+}\beta {+}\mu {+}1)
{+}\beta );
\end{array}\]
\[\begin{array}{rrl}
&((z-1)\p_z+\12(\alpha+\beta +\mu +1))&(1-z)\F_{\alpha,\beta ,\mu }\\[1ex]
&=\ \ \ \ \ \ (1-z)\F_{\alpha+1,\beta ,\mu +1}& ((z-1)\p_z+\12(\alpha+\beta +\mu +1)
,
\\[3ex]
&(z(1{-}z)\p_z{-}\12z(\alpha{+}\beta {-}\mu {+}1)
{+}\alpha)&(1-z)\F_{\alpha,\beta ,\mu }
\\[1ex]
&=\ \ \ \ \ \ (1-z)\F_{\alpha-1,\beta ,\mu -1}&(z(1{-}z)\p_z{-}\12z(\alpha{+}\beta {-}\mu {+}1)
{+}\alpha),
\\[3ex]
&((z-1)\p_z+\12(\alpha+\beta -\mu +1))&
(1-z)\F_{\alpha,\beta ,\mu }
\\[1ex]&=\ \ \ \ \ \ 
(1-z)\F_{\alpha+1,\beta ,\mu -1}&((z-1)\p_z+\12(\alpha+\beta -\mu +1)),
\\[3ex]
&(z(z{-}1)\p_z{+}\12z(\alpha{+}\beta {+}\mu {+}1)
-\alpha)&(1-z)\F_{\alpha,\beta ,\mu }
\\[1ex]
&=\ \ \ \ \ \ (1-z)\F_{\alpha-1,\beta ,\mu +1}&(z(z{-}1)\p_z{+}\12z(\alpha{+}\beta {+}\mu {+}1)
{-}\alpha).
\end{array}\]
Each of these commutation relations corresponds to a
root of the Lie
algebra $so(6)$.

\subsection{Canonical forms}

The natural weight of the 
 hypergeometric operator is $z^\alpha(1-z)^\beta$, so that
\[\F_{\alpha,\beta,\mu}=
z^{-\alpha}(1-z)^{-\beta}\partial_zz^{\alpha+1}(1-z)^{\beta+1}\partial_z
+\frac{\mu^2}{4}-\frac{(\alpha+\beta+1)^2}{4}.\]
The balanced form of the 
 hypergeometric operator is
\begin{eqnarray*}
&&z^{\frac{\alpha}{2}}(1-z)^{\frac{\beta}{2}}\F_{\alpha,\beta,\mu}
z^{-\frac{\alpha}{2}}(1-z)^{-\frac{\beta}{2}}\\
&=&
\partial_zz(1-z)\partial_z-\frac{\alpha^2}{4z}-\frac{\beta^2}{4(1-z)}
+\frac{\mu^2-1}{4}.
\end{eqnarray*}
Note that the symmetries $\alpha\to-\alpha$, $\beta\to-\beta$ and $\mu\to-\mu$ are obvious in the balanced form.

\ber In the literature, the balanced form of the hypergeometric equation is sometimes called the {\em generalized associated Legendre equation}. Its standard form according to \cite{NIST} is
\beq
(1-w^2)\partial_w^2-2w\partial_w+\nu(\nu+1)-\frac{\mu_1^2}{2(1-w)}-\frac{\mu_1^2}{2(1+w)}.
\eeq Thus $z=\frac{w+1}{2}$, moreover, $\mu_1$, $\mu_2$ and $\nu$ correspond to $\beta,\alpha$ and  $\frac{\mu}{2}-\frac12$.
\eer

\subsection{The hypergeometric function}
\label{ss-hyp}

$0$ is a regular singular point of 
the hypergeometric equation.  Its indices are 
 $0$ and $1-c$. 
The Frobenius method implies that, for $c\neq 0,-1,-2,\dots$,  the unique solution of
the hypergeometric equation equal to $1$
 at $0$ is given by the series 
\[F(a,b;c;z)=\sum_{j=0}^\infty
\frac{(a)_j(b)_j}{
(c)_j}\frac{z^j}{j!},\]
convergent for $|z|<1$. The function extends to the whole complex plane cut at $[1,\infty[$ and is 
called the {\em hypegeometric function}.
Sometimes it is more convenient to consider the function
\[ {\bf F}  (a,b;c;z):=\frac{F(a,b,c,z)}{\Gamma(c)}
=\sum_{j=0}^\infty
\frac{(a)_j(b)_j}{
\Gamma(c+j)}\frac{z^j}{j!}\]
 defined for all $a,b,c\in\cc$.
Another useful function proportional to ${}_2F_1$ is
\[ {\bf F}^\I  (a,b;c;z):=\frac{\Gamma(a)\Gamma(c-a)}{\Gamma(c)}
F(a,b;c;z)
=\sum_{j=0}^\infty
\frac{\Gamma(a+j)\Gamma(c-a)(b)_j}{
\Gamma(c+j)}\frac{z^j}{j!}.
\]It has
the integral representation
\begin{eqnarray}\label{eqa1}
&&\int_1^\infty t^{b-c}(t-1)^{c-a-1}(t-z)^{-b}\d t\\
&=&
 {\bf F}^\I  (a,b;c;z),\ \ \ \ \Re(c-a)>0,\ \Re a>0,\ \ \ z\not\in[1,\infty[.
\nonumber\end{eqnarray}

Indeed, by Theorem \ref{intr} the left hand side of  (\ref{eqa1})
 is annihilated by
the hypergeometric operator
(\ref{hy1-tra}). Besides, by  (\ref{b1}) it equals
 $\frac{\Gamma(a)\Gamma(c-a)}{\Gamma(c)}$ at $0$. So does
the right hand side. Therefore, Equation (\ref{eqa1}) follows by the uniqueness of the solution by the Frobenius method.

Another, closely related 
integral representation is
\beq
\frac{\sin\pi a}{\pi}{\bf F}^{\I}(a,b;c;z)=\frac{1}{2\pi\i}\int\limits_{[1,(z,0)^+,1]} (-t)^{b-c}(1-t)^{c-a-1}(z-t)^{-b}\d t.
\label{bequ}\eeq
It is proven essentially in the same way as (\ref{eqa1}), except that instead of
(\ref{b1}) we use (\ref{put}). 


We have the identities
\begin{eqnarray}\nonumber
&&F(a,b;c;z)\\\nonumber
&=&(1-z)^{c-a-b}F\left(c-a,c-b;c;z\right)\\\nonumber
&=&(1-z)^{-a}F\left(a,c-b;c;\frac{z}{z-1}\right)
\\\label{stan}
&=&(1-z)^{-b}F\left(c-a,b;c;\frac{z}{z-1}\right).
\end{eqnarray}
In fact, by the 3rd, 9th and 11th symmetry of Subsect. \ref{s-sym} all
these functions are annihilated by the hypergeometric operator. All of
them are $\sim1$ at $1$. Hence, by the uniqueness of the Frobenius method they coincide, at least for $c\neq0,-1,\dots$. By continuity, the identities hold for all $c\in\cc$.

Let us introduce new notation for various varieties of the
hypergeometric function involving the Lie-algebraic parameters instead of the
classical parameters. 
\begin{eqnarray*}
 F_{\alpha,\beta ,\mu }(z)&=&F\Bigl(
\frac{1+\alpha+\beta -\mu}{2},\frac{1+\alpha+\beta +\mu}{2};1+\alpha;z\Bigr),\\
 {\bf F}_{\alpha,\beta ,\mu }(z)&=&{\bf F} \Bigl(
\frac{1+\alpha+\beta -\mu}{2},\frac{1+\alpha+\beta +\mu}{2};1+\alpha;z\Bigr)\\
&=&
\frac{1}{\Gamma(\alpha+1)} F_{\alpha,\beta ,\mu }(z),\\
 {\bf F}_{\alpha,\beta ,\mu }^{\I}(z)&=&{\bf F}^{\I}\Bigl(
\frac{1+\alpha+\beta -\mu}{2},\frac{1+\alpha+\beta +\mu}{2};1+\alpha;z\Bigr)\\
&=&
\frac{\Gamma\big(\frac{1+\alpha+\beta-\mu}{2}\big)\Gamma\big(\frac{1+\alpha-\beta+\mu}{2}\big)}{\Gamma(\alpha+1)}
F_{\alpha,\beta ,\mu }(z).
\end{eqnarray*}

\subsection{Standard solutions -- Kummer's table}

To each of the singular points $0,1,\infty$ we can associate two solutions corresponding to its indices. Thus we obtain $3\times2=6$ solutions, which we will call {\em standard solutions}. Using the identites (\ref{stan}), each solution can be written in 4 distinct ways (not counting the trivial change of the sign in front of the last parameter). Thus we obtain a list of $6\times 4=24$ expressions for solutions of the hypergeometric equation, called sometimes {\em Kummer's table}. 

We describe the standard solutions to the hypergeometric equation in this section.
We will use consistently the Lie-algebraic parameters, which give much simpler expressions.

It follows from Thm \ref{intr} that for  appropriate contours $\gamma$  integrals of the form
\beq
\int_\gamma t^{\frac{-1-\alpha+\beta +\mu }{2}}
(t-1)^{\frac{-1+\alpha-\beta +\mu }{2}}(t-z)^{\frac{-1-\alpha-\beta -\mu }{2}}\d t\label{inte1}\eeq
are solutions of the hypergeometric equation. The integrand
has four singularities: $\{0,1,\infty,z\}$. It is natural to chose
$\gamma$  as the interval joining a pair of  singularities. This choice leads to $6$ standard solutions with the $\rm I$-type normalization.

\subsubsection{Solution $\sim1$ at $0$}
\label{ss-1}

If $\alpha\neq-1,-2,\dots$, then the following function is the unique solution  $\sim1$ at $0$: 
\begin{eqnarray*}
&& F_{\alpha,\beta ,\mu }(z)\\
&=&(1-z)^{-\beta } F_{\alpha,-\beta ,-\mu }(z)\\
&=&(1-z)^{\frac{-1-\alpha-\beta +\mu }{2}}
 F_{\alpha,-\mu ,-\beta }(\frac{z}{z-1})\\
&=&(1-z)^{\frac{-1-\alpha-\beta -\mu }{2}} F_{\alpha,\mu ,\beta }(\frac{z}{z-1}).
\end{eqnarray*}

An integral representation for $\Re(1+\alpha)> |\Re(\beta -\mu )|$:
\begin{eqnarray*}
\int_1^\infty t^{\frac{-1-\alpha+\beta +\mu }{2}}
(t-1)^{\frac{-1+\alpha-\beta +\mu }{2}}(t-z)^{\frac{-1-\alpha-\beta -\mu }{2}}\d t
&=& {\bf F}^\I  _{\alpha,\beta ,\mu }(z),\\&&
\ \ \ z\not\in[1,\infty[.\nonumber
\end{eqnarray*}

Note that all the identities of this subsubsection are the transcriptions of identities of Subsect.
\ref{ss-hyp} to the Lie-algebraic parameters.

\subsubsection{Solution $\sim z^{-\alpha}$ at $0$}
\label{ss-2}

If $\alpha\neq 1,2,\dots$, then the following function is the unique solution behaving
as $z^{-\alpha}$ at $0$:
\begin{eqnarray*}
&&z^{-\alpha} F_{-\alpha,\beta ,-\mu }(z)\\
&=&z^{-\alpha}(1-z)^{-\beta } F_{-\alpha,-\beta ,\mu }(z)\\
&=&z^{-\alpha}(1-z)^{\frac{-1+\alpha-\beta +\mu }{2}}
 F_{-\alpha,-\mu ,\beta }(\frac{z}{z-1})\\
&=&z^{-\alpha}(1-z)^{\frac{-1+\alpha-\beta -\mu }{2}} F_{-\alpha,\mu ,-\beta }(\frac{z}{z-1}).
\end{eqnarray*}

Integral representations for $\Re(1-\alpha)> |\Re(\beta -\mu )|$:
\begin{eqnarray*}
\int_0^z t^{\frac{-1-\alpha+\beta +\mu }{2}}
(1-t)^{\frac{-1+\alpha-\beta +\mu }{2}}(z-t)^{\frac{-1-\alpha-\beta -\mu }{2}}\d t
&=&
z^{-\alpha} {\bf F}^\I  _{-\alpha,\beta ,-\mu }(z),\\&&z\not\in]{-}\infty,0]{\cup}[1,\infty[
;\nonumber\\
\int_z^0 (-t)^{\frac{-1-\alpha+\beta +\mu }{2}}
(1-t)^{\frac{-1+\alpha-\beta +\mu }{2}}(t-z)^{\frac{-1-\alpha-\beta -\mu }{2}}\d t
&=&
(-z)^{-\alpha} {\bf F}^\I  _{-\alpha,\beta ,-\mu }(z),\\&&z\not\in[0,\infty[
.\nonumber
\end{eqnarray*}
To check these identities we note first that the integrals are solutions of the hypergeometric equation. By substituting $t=zs$ we easily  check that they have the correct behavior at zero. 

Of course, it is elementary to pass from the first identity, which is adapted to the region on the right of the singularity $z=0$  to the second, adapted to the region on the left of the singularity. For convenience we give both identities.

\subsubsection{Solution $\sim1$  at $1$}
\label{ss-3}

If $\beta \neq-1,-2,\dots$, then the following function is the unique solution 
 $\sim1$ at $1$:
\begin{eqnarray*}
&& F_{\beta ,\alpha,\mu }(1-z)\\
&=&z^{-\alpha} F_{\beta ,-\alpha,-\mu }(1-z)\\
&=&z^{\frac{-1-\alpha-\beta +\mu }{2}}
 F_{\beta ,-\mu ,-\alpha}(1-z^{-1})\\
&=&z^{\frac{-1-\alpha-\beta -\mu }{2}}
 F_{\beta ,\mu ,\alpha}(1-z^{-1})
\end{eqnarray*}

Integral representation for $\Re(1+\beta )> |\Re(\alpha-\mu )|$:
\begin{eqnarray*}
\int_{-\infty}^0 (-t)^{\frac{-1-\alpha+\beta +\mu }{2}}
(1-t)^{\frac{-1+\alpha-\beta +\mu }{2}}(z-t)^{\frac{-1-\alpha-\beta -\mu }{2}}\d t
&=&
 {\bf F}^\I  _{\beta ,\alpha ,\mu }(1-z),\\
&&z\not\in]-\infty,0].
\end{eqnarray*}

\subsubsection{Solution $\sim(1-z)^{-\beta }$ at $1$}
\label{ss-4}
If $\beta \neq1,2,\dots$, then the following function is the unique solution of the
hypergeometric equation $\sim(1-z)^{-\beta }$ at $1$:

\begin{eqnarray*}
&&(1-z)^{-\beta } F_{-\beta ,\alpha,-\mu }(1-z)\\
&=&z^{-\alpha}(1-z)^{-\beta } F_{-\beta ,-\alpha,\mu }(1-z)\\
&=&z^{\frac{-1-\alpha+\beta -\mu }{2}}(1-z)^{-\beta }
 F_{-\beta ,\mu ,-\alpha}(1-z^{-1})\\
&=&z^{\frac{-1-\alpha+\beta +\mu }{2}}(1-z)^{-\beta }
 F_{-\beta ,-\mu ,\alpha}(1-z^{-1}).
\end{eqnarray*}

Integral representations for $\Re(1-\beta )> |\Re(\alpha+\mu )|$:
\begin{eqnarray*}
\int_z^1 t^{\frac{-1-\alpha+\beta +\mu }{2}}
(1-t)^{\frac{-1+\alpha-\beta +\mu }{2}}(t-z)^{\frac{-1-\alpha-\beta -\mu }{2}}\d t
&=&(1-z)^{-\beta } {\bf F}^\I  _{-\beta ,\alpha,-\mu }(1-z),\\
&&z\not\in]-\infty,0]\cup[1,\infty[;\\
\int_1^z t^{\frac{-1-\alpha+\beta +\mu }{2}}
(t-1)^{\frac{-1+\alpha-\beta +\mu }{2}}(z-t)^{\frac{-1-\alpha-\beta -\mu }{2}}\d t
&=&(z-1)^{-\beta } {\bf F}^\I  _{-\beta ,\alpha,-\mu }(1-z),\\
&&z\not\in]-\infty,1].
\end{eqnarray*}

\subsubsection{Solution $\sim z^{-a}$ at $\infty$}
\label{ss-5}

If $\mu \neq 1,2\dots$, then the following function is the unique solution of the
hypergeometric equation $\sim (-z)^{-a}=(-z)^{\frac{-1-\alpha-\beta +\mu }{2}}$ at $\infty$:

\begin{eqnarray*}
&&(-z)^{\frac{-1-\alpha-\beta +\mu }{2}} F_{-\mu ,\beta ,-\alpha}(z^{-1})\\
&=&(-z)^{\frac{-1-\alpha+\beta +\mu }{2}}(1-z)^{-\beta } F_{-\mu ,-\beta ,\alpha}(z^{-1})\\
&=&(1-z)^{\frac{-1-\alpha-\beta +\mu }{2}} F_{-\mu ,\alpha,-\beta }((1-z)^{-1})\\
&=&(-z)^{-\alpha}(1-z)^{\frac{-1+\alpha-\beta +\mu }{2}} F_{-\mu ,-\alpha,\beta }((1-z)^{-1}).
\end{eqnarray*}

Integral representations for $\Re(1-\mu )> |\Re(\alpha+\beta )|$:
\begin{eqnarray*}
\int_z^\infty t^{\frac{-1-\alpha+\beta +\mu }{2}}
(t-1)^{\frac{-1+\alpha-\beta +\mu }{2}}(t-z)^{\frac{-1-\alpha-\beta -\mu }{2}}\d t
&=&
z^{\frac{-1-\alpha-\beta -\mu }{2}} {\bf F}^\I   _{-\mu ,\beta ,-\alpha}(z^{-1})
,\\
&&z\not\in]-\infty,1];\\
\int_{-\infty}^z (-t)^{\frac{-1-\alpha+\beta +\mu }{2}}
(1-t)^{\frac{-1+\alpha-\beta +\mu }{2}}(z-t)^{\frac{-1-\alpha-\beta -\mu }{2}}\d t
&=&
(-z)^{\frac{-1-\alpha-\beta -\mu }{2}} {\bf F}^\I   _{-\mu ,\beta ,-\alpha}(z^{-1})
,\\
&&z\not\in]0,\infty].
\end{eqnarray*}

\subsubsection{Solution $\sim z^{-b}$ at $\infty$}
\label{ss-6}

If $\mu \neq-1,-2,\dots$, then the following function is the unique solution of the
hypergeometric equation
$\sim(-z)^{-b}=(-z)^{\frac{-1-\alpha-\beta -\mu }{2}}$
at $\infty$:

\begin{eqnarray*}
&&(-z)^{\frac{-1-\alpha-\beta -\mu }{2}} F_{\mu ,\beta ,\alpha}(z^{-1})\\
&=&(-z)^{\frac{-1-\alpha+\beta -\mu }{2}}(1-z)^{-\beta } F_{\mu ,-\beta ,-\alpha}(z^{-1})\\
&=&(1-z)^{\frac{-1-\alpha-\beta -\mu }{2}} F_{\mu ,\alpha,\beta }((1-z)^{-1})\\
&=&(-z)^{-\alpha}(1-z)^{\frac{-1+\alpha-\beta -\mu }{2}} F_{\mu ,-\alpha,-\beta }((1-z)^{-1})
\end{eqnarray*}

Integral representations for $\Re(1+\mu )> |\Re(\alpha-\beta )|$:
\begin{eqnarray*}
\int_0^1 t^{\frac{-1-\alpha+\beta -\mu }{2}}
(1-t)^{\frac{-1+\alpha-\beta +\mu }{2}}(t-z)^{\frac{-1-\alpha-\beta -\mu }{2}}\d t
&=&
(-z)^{\frac{-1-\alpha-\beta +\mu }{2}} {\bf F}^\I  _{\mu ,\beta ,\alpha}(z^{-1}),\\
&&z\not\in[0,\infty[;\\
\int_0^1 t^{\frac{-1-\alpha+\beta -\mu }{2}}
(1-t)^{\frac{-1+\alpha-\beta +\mu }{2}}(z-t)^{\frac{-1-\alpha-\beta -\mu }{2}}\d t
&=&
z^{\frac{-1-\alpha-\beta +\mu }{2}} {\bf F}^\I  _{\mu ,\beta ,\alpha}(z^{-1}),\\
&&z\not\in[-\infty,1[.
\end{eqnarray*}
\subsection{Connection formulas}
We use the solutions $\sim 1$ and $\sim z^{-\alpha}$ at $0$ 
as the basis. We show how the other solutions decompose in this basis.

For the first pair of relations we assume that $
 z\not\in]-\infty,0]{\cup}[1,\infty[$:
\begin{eqnarray*}
  {\bf F}  _{\beta ,\alpha,\mu }(1-z)
&=&
\frac{\pi}{\sin\pi(-\alpha)
\Gamma\left(\frac{1-\alpha+\beta -\mu }{2}\right)
\Gamma\left(\frac{1-\alpha+\beta +\mu }{2}\right)} {\bf F}  _{\alpha,\beta ,\mu }(z)\\
&&\!\!\!\!\!\!\!\!\!\!\!\!\!\!\!+\frac{\pi}{\sin\pi \alpha
\Gamma\left(\frac{1+\alpha+\beta -\mu }{2}\right)
\Gamma\left(\frac{1+\alpha+\beta +\mu }{2}\right)}z^{-\alpha} {\bf F}  _{-\alpha,\beta ,-\mu }(z),\\[3ex]
(1-z)^{-\beta }  {\bf F}  _{-\beta ,\alpha,-\mu }(1-z)
&=&
\frac{\pi}{\sin\pi(-\alpha)
\Gamma\left(\frac{1-\alpha-\beta +\mu }{2}\right)
\Gamma\left(\frac{1-\alpha-\beta -\mu }{2}\right)} {\bf F}  _{\alpha,\beta ,\mu }(z)\\
&&\!\!\!\!\!\!\!\!\!\!\!\!\!\!\!+\frac{\pi}{\sin\pi \alpha
\Gamma\left(\frac{1+\alpha-\beta +\mu }{2}\right)
\Gamma\left(\frac{1+\alpha-\beta -\mu }{2}\right)}z^{-\alpha} {\bf F}  _{-\alpha,\beta ,-\mu }(z)
.\end{eqnarray*}
For the second pair we assume that $ z\not\in[0,\infty[$
\begin{eqnarray*}
(-z)^{\frac{-1-\alpha-\beta +\mu }{2}}  {\bf F}  _{-\mu ,\beta ,-\alpha}(z^{-1})
&=&
\frac{\pi}{\sin\pi(-\alpha)
\Gamma\left(\frac{1-\alpha-\beta -\mu }{2}\right)
\Gamma\left(\frac{1-\alpha+\beta -\mu }{2}\right)} {\bf F}  _{\alpha,\beta ,\mu }(z)\\
&&\!\!\!\!\!\!\!\!\!\!\!\!\!\!\!\!\!\!\!\!\!\!\!\!\!+\frac{\pi}{\sin\pi \alpha
\Gamma\left(\frac{1+\alpha+\beta -\mu }{2}\right)
\Gamma\left(\frac{1+\alpha-\beta -\mu }{2}\right)}(-z)^{-\alpha}
     {\bf F}  _{-\alpha,\beta ,-\mu }(z),
\\(-z)^{\frac{-1-\alpha-\beta -\mu }{2}}  {\bf F}  _{\mu ,\beta ,\alpha}(z^{-1})
&=&
\frac{\pi}{\sin\pi(-\alpha)
\Gamma\left(\frac{1-\alpha-\beta +\mu }{2}\right)
\Gamma\left(\frac{1-\alpha+\beta +\mu }{2}\right)} {\bf F}  _{\alpha,\beta ,\mu }(z)\\
&&\!\!\!\!\!\!\!\!\!\!\!\!\!\!\!\!\!\!\!\!\!\!\!\!\!+\frac{\pi}{\sin\pi \alpha
\Gamma\left(\frac{1+\alpha+\beta +\mu }{2}\right)
\Gamma\left(\frac{1+\alpha-\beta +\mu }{2}\right)}(-z)^{-\alpha} {\bf F}  _{-\alpha,\beta ,-\mu }(z)
.\end{eqnarray*}
The connection formulas are easily derived from the integral representations by looking at the behavior around $0$.

\subsection{Recurrence relations}
\label{s3.15}
The following recurrence relations follow easily from the commutation
relations of Subsect. \ref{commu}:
\begin{eqnarray*}
\p_z {\bf F}^\I  _{\alpha,\beta ,\mu }(z)&=&
\frac{1{+}\alpha{+}\beta {+}\mu }{2} {\bf F}^\I  _{\alpha+1,\beta +1,\mu }(z),\\[2ex]
(z(1{-}z)\p_z{+}\alpha(1{-}z){-}\beta z) {\bf F}^\I  _{\alpha,\beta ,\mu }(z)&=&
 \frac{{-}1{+}\alpha{+}\beta{+}\mu}{2}{\bf F}^\I  _{\alpha-1,\beta -1,\mu }(z),\\[3ex]
((1-z)\p_z-\beta ) {\bf F}^\I  _{\alpha,\beta ,\mu }(z)&=&\frac{1{+}\alpha{-}\beta {-}\mu }{2}{\bf F}^\I  _{\alpha{+}1,\beta {-}1,\mu }(z),\\[2ex]
(z\p_z+\alpha) {\bf F}^\I  _{\alpha,\beta ,\mu }(z)&=& \frac{1{+}\alpha{-}\beta{+}\mu}{2}{\bf F}^\I  _{\alpha-1,\beta +1,\mu }(z),\\[4ex]
 \left(z\p_z+\frac{1+\alpha+\beta +\mu }{2}\right) {\bf F}^\I  _{\alpha,\beta ,\mu }(z)&=&
\frac{1+\alpha+\beta +\mu }{2} {\bf F}^\I  _{\alpha,\beta +1,\mu +1}(z)
, \\[2ex]
 \left(z(1{-}z)\p_z{-}\beta {+}\frac{1{+}\alpha{+}\beta {-}\mu }{2}(1{-}z)
\right) {\bf F}^\I  _{\alpha,\beta ,\mu }(z)&=&\frac{1+\alpha-\beta -\mu }{2} {\bf F}^\I  _{\alpha,\beta -1,\mu -1}(z),\\[3ex]
 \left(z\p_z+\frac{1+\alpha+\beta -\mu }{2}\right) {\bf F}^\I  _{\alpha,\beta ,\mu }(z)&=&
\frac{1+\alpha+\beta -\mu }{2} {\bf F}^\I  _{\alpha,\beta +1,\mu -1}(z)
, \\[2ex]
 \left(z(1{-}z)\p_z{-}\beta {+}\frac{1{+}\alpha{+}\beta {+}\mu }{2}(1{-}z)\right) {\bf F}^\I  _{\alpha,\beta ,\mu }(z)&=&
\frac{1+\alpha-\beta +\mu }{2}
 {\bf F}^\I  _{\alpha,\beta -1,\mu +1}(z)
 ,\\[4ex]
\left((z-1)\p_z+\frac{1+\alpha+\beta +\mu }{2}\right) {\bf F}^\I  _{\alpha,\beta ,\mu }(z)&=&
\frac{1{+}\alpha{+}\beta {+}\mu }{2}
{\bf F}^\I  _{\alpha{+}1,\beta ,\mu {+}1}(z)
, \\[2ex]
 \left(z(1{-}z)\p_z{+}\alpha{-}\frac{1{+}\alpha{+}\beta {-}\mu }{2}z\right)
 {\bf F}^\I  _{\alpha,\beta ,\mu }(z)&=&\frac{{-}1{+}\alpha{-}\beta{+}\mu}{2} {\bf F}^\I  _{\alpha-1,\beta ,\mu -1}(z),
\\[3ex]
 \left((z-1)\p_z+\frac{1+\alpha+\beta -\mu }{2}\right)
 {\bf F}^\I  _{\alpha,\beta ,\mu }(z)&=&\frac{1{+}\alpha{-}\beta {-}\mu }{2}
 {\bf F}^\I  _{\alpha{+}1,\beta ,\mu {-}1}(z)
, \\[2ex]
 \left(z(1{-}z)\p_z{+}\alpha{-}\frac{1{+}\alpha{+}\beta {+}\mu }{2}z\right) {\bf F}^\I  _{\alpha,\beta ,\mu }(z)&=&\frac{{-}1{+}\alpha{+}\beta{-}\mu}{2} {\bf F}^\I  _{\alpha-1,\beta ,\mu +1}(z)
 .
 \end{eqnarray*}

\subsection{Additional recurrence relations}
\label{addi}

There exist other, more complicated recurrence relations for hypergeometric functions, for example
\begin{eqnarray}
&&\Big(\frac{(1{+}\alpha{+}\beta {+}\mu )({-}1{-}\alpha{+}\beta {-}\mu )}{4}\nonumber\\
&&+\frac{(1{+}\alpha{+}\beta {+}\mu )(\mu {+}1)}{2}z
-(1{+}\mu )z(1{-}z)\p_z\Big){\bf F}_{\alpha,\beta ,\mu }\nonumber
\\
&=&\frac{(1{+}\alpha{+}\beta {+}\mu )({-}1{-}\alpha{+}\beta {-}\mu )}{4}{\bf F}_{\alpha,\beta ,\mu +2}(z),
\label{r1}\\[5ex]
&&\Big(\frac{(1{+}\alpha{+}\beta {-}\mu )({-}1{-}\alpha{+}\beta {+}\mu )}{4}\nonumber\\
&&+\frac{(1{+}\alpha{+}\beta {-}\mu )(-\mu {+}1)}{2}z
-(1{-}\mu )z(1{-}z)\p_z\Big){\bf F}_{\alpha,\beta ,\mu }\nonumber
\\
&=&\frac{(1{+}\alpha{+}\beta {-}\mu )({-}1{-}\alpha{+}\beta {+}\mu )}{4}{\bf F}_{\alpha,\beta ,\mu -2}(z)
\label{r2}.
\end{eqnarray}
Note that (\ref{r1}) follows from the 6th and 7th recurrence relation,
 and (\ref{r2}) follows from the 5th and 8th of Subsect.
 \ref{s3.15}.

\subsection{Degenerate case}

$\alpha=m\in\zz$ is the degenerate case of the hypergeometric equation at $0$.
We have then
\[{\bf F}
(a,b;1+m;z)=\sum_{n=\max(0,-m)}\frac{(a)_n(b)_n}{n!(m+n)!}z^n.\]
This easily implies the identity
\beq(a-m)_m(b-m)_m{\bf F}(a,b;1+m;z)=z^{-m}
{\bf F}(a-m,b-m;1-m;z).
\label{pwr}\eeq
Thus the two standard solutions determined by the behavior at zero are proportional to one another. 

One can also see the degenerate case in the integral representation (\ref{f4}).
If we go around $0,z$, the phase of the integrand changes by $\e^{\i2\pi c}=\e^{\i2\pi\alpha}$. Therefore, if 
 $\alpha=m\in\zz$, then the loop around $0,z$ is closed on the Riemann surface of the integrand. 

 We have an additional integral representation and a generating function:
\begin{eqnarray*}
\frac{1}{2\pi \i}\int\limits_{[(0,z)^+]}(1-t)^{-a}(1-z\slash t)^{-b}
t^{-m-1}\d t&=&
(a)_m {\bf F}_{m,a+b-1,-a+b-m}(z)
\\&=&
z^{-m}(b)_{-m} {\bf F}  _{-m,a+b-1,a-b+m}(z),\\[4ex]
(1-t)^{-a}(1-z\slash t)^{-b}
&=&\sum_{m\in\zz}t^m(a)_m {\bf F}  _{m,a+b-1,b-a-m}(z).\nonumber
\end{eqnarray*}
To see the integral representation  we note that the integral on the l.h.s. is annihilated by the hypergeometric operator. Then we check that its value at zero equals
\[\frac{1}{2\pi\i}\int\limits_{[0^+]}(1-t)^{-a}t^{-m-1}\d t=\frac{(a)_m}{m!},\]
see (\ref{double1}). 

The second identity follows from (\ref{pwr}). Another way to see it is
to make the substitution $t=\frac{z}{s}$. Note that $[(0,z)^+]$ becomes $[(\infty,1)^+]$, which coincides with $[(0,z)^-]$. Then we change the sign in front of the integral and the orientation of the contour of integration, obtaining
\[\frac{z^{-m}}{2\pi\i}\int_{[(0,z)^+]}(1-s)^{-b}(1-z/s)^{-a}s^{-m-1}\d s.\]
Finally, we apply the first integral representation again.

The generating function follows from the integral representation.
\subsection{Jacobi polynomials}

If $-a=n=0,1,\dots$, then  hypergeometic functions are polynomials. 
We will call them the {\em Jacobi polynomials}.

Following
Subsect. 
\ref{Hypergeometric type polynomials}, the Jacobi polynomials are
  defined by the
 Rodriguez-type formula
\[
\begin{array}{rl}R^{\alpha,\beta}_n(z):=
\frac{(-1)^n}{n!}z^{-\alpha}(z-1)^{-\beta}\p_z^nz^{\alpha+n}(z-1)^{\beta+n}
.\end{array}\]

\ber
In most of the literature, the  Jacobi polynomials are slightly different:
\[P_n^{\alpha,\beta}(z):=R_n^{\alpha,\beta}\Big(\frac{1-z}{2}\Big)
=(-1)^nR_n^{\beta,\alpha}\Big(\frac{1+z}{2}\Big).\]
\label{jaco}\eer
The equation:
\begin{eqnarray*}
&&0=\F(-n,1+\alpha+\beta+n;\beta+1;z,\p_z)P_n^{\alpha,\beta}(z)\\
&&=\Big(z(1-z)\p_z^2+
\big((1+\alpha)(1-z)-(1+\beta
)z\big)\p_z+n(n+\alpha+\beta+1)\Big)P_n^{\alpha,\beta}(z).
\end{eqnarray*}
Generating functions:
\begin{eqnarray*}
(1+t(1-z))^{\alpha}(1-tz)^{\beta}
&=&\sum\limits_{n=0}^\infty t^nR_n^{\alpha-n,\beta-n}(z),\\
(1+zt)^{-1-\alpha-\beta}(1+t)^{\alpha}
&=&\sum\limits_{n=0}^\infty t^nR_n^{\alpha-n,\beta}(z),\\
(1+(z-1)t)^{-1-\alpha-\beta}(1-t)^{\beta}
&=&\sum\limits_{n=0}^\infty t^nR_n^{\alpha,\beta-n}(z).
\end{eqnarray*}
Integral representations:
\begin{eqnarray*}
R_n^{\alpha,\beta}(z)
&=&\frac{1}{2\pi \i}\int\limits_{[0^+]}
(1+(1-z)t)^{\alpha+n}(1-zt)^{\beta+n}t^{-n-1}\d t\\
&=&\frac{1}{2\pi \i}\int\limits_{[0^+]}
(1+zt)^{-\alpha-\beta-n-1}(1+t)^{\alpha+n}t^{-n-1}\d t
\\
&=&\frac{1}{2\pi \i}\int\limits_{[0^+]}
(1+(z-1)t)^{-\alpha-\beta-n-1}(1-t)^{\beta+n}t^{-n-1}\d t.\end{eqnarray*}
Discrete symmetries:
\begin{eqnarray*}
R_n^{\alpha,\beta}(z)
&=&(1-z)^n
R_n^{\alpha,-1-\alpha-\beta-2n}\Big(\frac{z}{z-1}\Big)\\
=\ (-1)^nR_n^{\beta,\alpha}(1-z)
&=& (-z)^nR_n^{\beta,-1-\alpha-\beta-2n}\Big(\frac{z-1}{z}\Big)
\\=\ 
z^nR_n^{-1-\alpha-\beta-2n,\beta}\Big(\frac{1}{z}\Big)&=&
(z-1)^nR_n^{-1-\alpha-\beta-2n,\alpha}\Big(\frac{1}{1-z}\Big).
\end{eqnarray*}
Recurrence relations:
\begin{eqnarray*}
\p_zR_n^{\alpha,\beta}(z)&=&-(\alpha+\beta+n+1)
R_{n-1}^{\alpha+1,\beta+1}(z),\\
(z(1-z)\p_z-\alpha(z-1)-\beta z)R_n^{\alpha,\beta}(z)&=&(n+1)
R_{n+1}^{\alpha-1,\beta-1}(z),\\[4ex]
((1-z)\p_z-\beta)R_n^{\alpha,\beta}(z)&=&-(\beta+n)
R_n^{\alpha-1,\beta+1}(z),\\
(z\p_z+\alpha)R_n^{\alpha,\beta}(z)&=&(\beta+n)R_n^{\alpha-1,\beta+1}(z),\\[6ex]
 \left(z\p_z-n\right)R_n^{\alpha,\beta}(z)&=&-(\alpha+n)
R_{n-1}^{\alpha,\beta+1}(z)
, \\
 \left(z(1-z)\p_z+1+\alpha+n-(1+\alpha+\beta+n)z\right)
R_n^{\alpha,\beta}(z)&=&(n+1)
R_{n+1}^{\alpha,\beta-1}(z)
 ,\\[4ex]
 \left(z\p_z+1+\alpha+\beta+n\right)R_n^{\alpha,\beta}(z)&=&(1+\alpha+\beta+n)
R_n^{\alpha,\beta+1}(z)
, \\
 \left(z(1-z)\p_z-n-\beta+nz\right)R_n^{\alpha,\beta}(z)&=&
-(\beta+n)R_n^{\alpha,\beta-1}(z),
\\[6ex]
 \left((z-1)\p_z-n\right)R_n^{\alpha,\beta}(z)&=&
(\beta+n)R_{n-1}^{\alpha+1,\beta}(z) 
, \\
 \left(z(1-z)\p_z+\alpha-(1+\alpha+\beta+n)z\right)R_n^{\alpha,\beta}(z)
&=&(n+1)R_{n+1}^{\alpha-1,\beta}(z)
 ,\\[4ex]
 \left((z-1)\p_z+1+n+\alpha+\beta\right)R_n^{\alpha,\beta}(z)
&=&(1+n+\alpha+\beta)R_n^{\alpha+1,\beta}(z)
, \\  
 \left(z(1-z)\p_z+\alpha+nz\right)R_n^{\alpha,\beta}(z)&=&
(n+\alpha)R_n^{\alpha-1,\beta}(z).
\end{eqnarray*}

 The first, second, resp. third integral representation is
easily seen to be equivalent to the first, second, resp. third
generating function. The first follows immediately from the
Rodriguez-type formula.

The symmetries can be interpreted as a subset of Kummer's table. The first line corresponds to the symmetries of the solution regular at $0$, see 
(\ref{stan}) (or Subsubsect. \ref{ss-1}). Note that from 4 expressions
in (\ref{stan}) only the first and the third survive, since $n=-a$ should not change. The second line corresponds to the solution regular at 
$1$ (Subsubsect. \ref{ss-3}{}), finally the third line to the solution $\sim z^{-a}=z^n$ (Subsubsect. \ref{ss-5}{}).

The differential equation, the Rodriguez-type formula, the
first generating function, the  first integral representation and the
first pair of recurrence relations are special cases of the corresponding formulas of Subsect. 
\ref{Hypergeometric type polynomials}.

Note that Jacobi polynomials are regular at $0$, $1$, and behave as $z^n$ in infinity. Thus (up to  coefficients) they coincide with the 3 standard solutions.
They have the following
values at $0$, $1$ and the behavior at $\infty$:
\begin{eqnarray*}R_n^{\alpha,\beta}(0)&=&\frac{(\alpha+1)_n}{n!},\ \ \ 
R_n^{\alpha,\beta}(1)\ =\ (-1)^n\frac{(\beta+1)_n}{n!},\\
\lim\limits_{z\to\infty}\frac{R_n^{\alpha,\beta}(z)}{z^n}&=&
(-1)^n\frac{(\alpha+\beta+n+1)_n}{n!}.\end{eqnarray*}

We have several alternative expressions for Jacobi polynomials:
\begin{eqnarray*}
R_n^{\alpha,\beta}(z)&:=&\lim\limits_{\nu\to n}
(-1)^n(\nu-n)
{\bf F}_{\alpha,\beta,2\nu+\alpha+\beta+1}^{\I}(z)
=\frac{(\alpha+1)_n}{n!}F_{\alpha,\beta,2n+\alpha+\beta+1}(z)\\
&=&\frac{\Gamma(\alpha+1+n)}{\Gamma(\alpha+1)\Gamma(n+1)}F(-n,n+\alpha+\beta+1;\alpha+1;z)\\
&=&\sum\limits_{j=0}^n
\frac{(1+\alpha+j)_{n-j}(1+\alpha+\beta+n)_j}{j!(n-j)!}(-z)^j.
\end{eqnarray*}

One way to derive the first of the above identities is to use integral representation
(\ref{bequ}). Using that $a$ is an integer we can replace the open curve $[1,(0,z)^+,1]$
with a closed loop $[\infty^-]$:
\begin{eqnarray*}&&
\lim\limits_{\nu\to n}
(-1)^n(\nu-n)
{\bf F}_{\alpha,\beta,2\nu+\alpha+\beta+1}^{\I}(z)\\
&=&
\lim\limits_{\nu\to n}\frac{\sin \nu\pi}{\pi}{\bf F}_{\alpha,\beta,2\nu+\alpha+\beta+1}^{\I}(z)\\
&=&\frac{1}{2\pi\i}
\int_{[\infty^-]} (-s)^{\beta+n}(1-s)^{\alpha+n}
(z-s)^{-1-\alpha-\beta-n}\d s.
\end{eqnarray*}
Then, making the substitions
$s=z-\frac1t$, $s=zt$, resp. $s=(z-1)t$ we obtain the 1st, 2nd, resp. 3rd integral representation.

Additional identities valid in the degenerate case:
\begin{eqnarray*}
R_n^{\alpha,\beta}(z)&=&\frac{(n+1)_\alpha}{(\beta+n+1)_\alpha}(-z)^{-\alpha}
R_{n+\alpha}^{-\alpha,\beta}(z),\ \ \alpha\in\zz;\\
R_n^{\alpha,\beta}(z)&=&\frac{(n+1)_\beta}{(\alpha+n+1)_\beta}(1-z)^{-\beta}
R_{n+\beta}^{\alpha,-\beta}(z),\ \ \beta\in\zz;\\
R_n^{\alpha,\beta}(z)&=&
(-z)^{-\alpha}
(1-z)^{-\beta} 
R_{n+\alpha+\beta}^{-\alpha,-\beta}(z),\ \ \alpha,\beta\in\zz.
\end{eqnarray*}

There is a region where Jacobi polynomials are zero. This happens iff
$\alpha,\beta\in\zz$ and $\alpha,\beta$ are in the triangle
\begin{eqnarray}\nonumber
0&\leq&\alpha+n,\\\nonumber
0&\leq&\beta+n,\\
0&\leq&-\alpha-\beta-n-1.\label{mumu}
\end{eqnarray}

In the analysis of symmetries of Jacobi polynomials it is useful to go back to the Lie-algebraic parameters, more precisely, to set $\mu:=-\alpha-\beta-2n-1$.
Then (\ref{mumu}) acquires a more symmetric form, since we can replace its last condition by
 \[0 \leq\mu+n.\]

One can distinguish 3  strips where Jacobi polynomials have special properties. Note that the intersection of the strips below is precisely the triangle described in (\ref{mumu}). 
\ben
\item $\mu\in\zz$ and $-n\leq\mu\leq-1$ or, equivalently,  $\alpha+\beta\in\zz$ and $-2n\leq\alpha+\beta\leq-n-1$. Then $R_n^{\alpha,\beta}=0$ or
\[\deg R_n^{\alpha,\beta}=\mu+n=-\alpha-\beta-n-1.\]
\item $\alpha\in\zz$ and $-n\leq\alpha\leq-1$. Then $R_n^{\alpha,\beta}=0$ or
\[R_n^{\alpha,\beta}=z^{-\alpha} W,\ \ \ \ W\ \ \hbox{not divisible by} \ z.\]
\item $\beta\in\zz$ and $-n\leq\beta\leq-1$. Then $R_n^{\alpha,\beta}=0$ or
\[R_n^{\alpha,\beta}=(z-1)^{-\beta} V,\ \ \ \ V\ \ \hbox{not divisible by} \ z-1.\]
\een
These regions are presented in the following picture:
\begin{center}\includegraphics[width=7cm,totalheight=7cm]{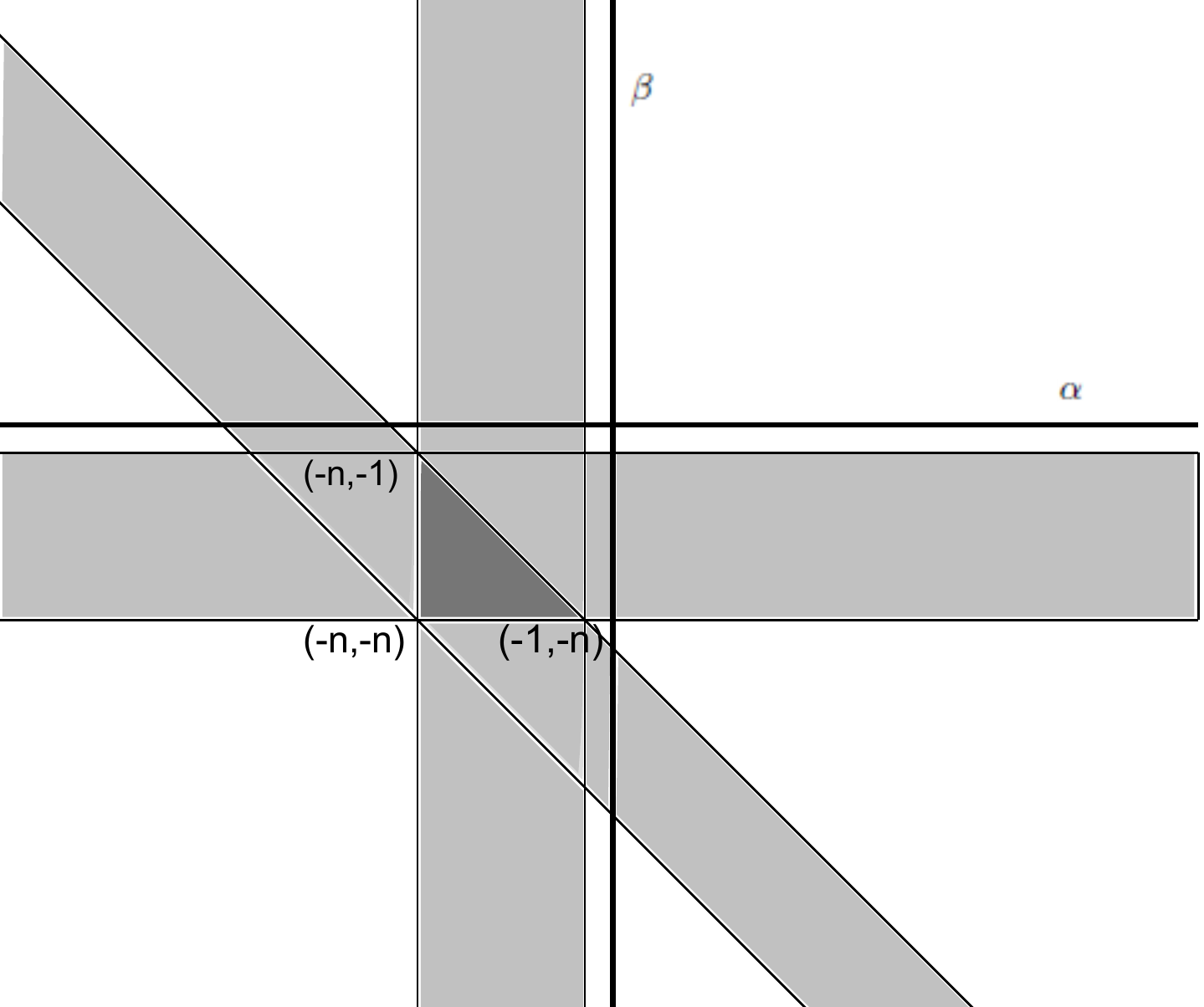}
\end{center}

Finally Jacobi polynomials satisfy some identities related to Subsect. \ref{addi}.
An additional generating function:
\begin{eqnarray}\nonumber
2^{\alpha+\beta}r^{-1}(1-t+r)^{-\alpha}(1+t+r)^{-\beta}
&=&\sum\limits_{n=0}^\infty t^nR_n^{\alpha,\beta}(z),\\
\hbox{where}\ \ \ r=\sqrt{(1-t)^2+4zt}.&&
\label{gg1}\end{eqnarray}
Additional recurrence relations:
\begin{eqnarray*}
\Big((n+\alpha+\beta+1)\big((n+\beta+1)-
(2n+\alpha+\beta+2)z\big) \hspace{-10ex} &&\\
\!\!\!\!\!+(2n+\alpha+\beta+2)z(1-z)\p_z\Big)R_n^{\alpha,\beta}(z)
&=&(n+\alpha+\beta+1)(n+1)R_{n+1}^{\alpha,\beta}(z),\\[3ex]
\Big(n\big((n+\alpha)-(2n+\alpha+\beta)z\big)\ \ \ \ \ \ \ \ \ &&\\-(2n+\alpha+\beta)z(1-z)\p_z\Big)
R_n^{\alpha,\beta}(z)&=&(n+\alpha)(n+\beta)R_{n-1}^{\alpha,\beta}(z).
\end{eqnarray*}

\subsection{Special cases}

Beside the polynomial and degenerate cases, the hypergeometric equation has a number of other special cases. In their description  most of the time we will use  the Lie-algebraic parameters,
which are here more convenient than the classical parameters.

\subsubsection{Gegenbauer equation through an affine transformation}

Consider a
hypergeometric equation whose two parameters coincide up to a sign.
After applying an appropriate symmetry we can assume that they are at the first and second place, and that they are equal to one another. In other words,
 $\alpha=\beta$. A simple affine transformation (\ref{ha2}) can be then applied to obtain a reflection invariant equation called the Gegenbauer equation.  We study it separately in Sect. \ref{s7}. 

\subsubsection{Gegenbauer equation through a quadratic transformation}
Hypergeometric equations  with one of the  parameters equal to
$\frac12$ or $-\frac12$ also enjoy special properties. After applying,
if needed, one of the symmetries, we can assume that $\mu=\pm\frac12$.
Then identity (\ref{1/2}) or (\ref{-1/2}) leads to the Gegenbauer equation.
\subsubsection{Chebyshev equation}
Even more special properties have equations with a pair of parameters  $\pm\frac12$. After applying  one of the symmetries we can assume that $\alpha=\beta=\frac12$. Thus we are reduced to the Chebyshev equation of the first kind; see  (\ref{cheb1}). Another option is to reduce it to the Chebyshev equation of the second kind, which corresponds to $\alpha=\beta=-\frac12$;
see (\ref{cheb2}).

\subsubsection{Legendre equation}
Let $\cL$ be the sublattice of $\zz^3$ consisting of points whose
sum of
coordinates is even. It is a sublattice of $\zz^3$ of degree $2$.
By using recurrence relations of
Subsect. \ref{s3.15} we can pass from hypergeometric functions with
given Lie-algebraic parameters $(\alpha,\beta,\mu)$ to parameters from
$(\alpha,\beta,\mu)+\cL$. 

This is especially useful in the degenerate
case, when some of the parameters are integers. In particular, if two
of the parameters are integers, by applying recurrence relations we
can make both of them zero. By applying an appropriate symmetry
we can assume that $\alpha=\beta=0$. Thus we obtain  the
Legendre equation, see (\ref{legen}).

\subsubsection{Elementary solutions}
One can easily check that
\[F(a,b;b;;z)=F_{b-1,a,b-a}(z)=(1-z)^{-a}.\]
Therefore, using  Kummer's table and recurrence relations we see that if
\beq\epsilon_1\alpha+\epsilon_2\beta+\epsilon_3\mu\ \ \hbox{ is an odd
  integer for some}\ \ \
\epsilon_1,\epsilon_2,\epsilon_3\in\{-1,1\}\label{condi}\eeq
then $F_{\alpha,\beta,\mu}$ is an elementary function involving power
functions, but not logarithms.

\subsubsection{Fully degenerate case}
An interesting situation arises if $\alpha,\beta,\mu\in\zz$,
that is, we have the degenerate case at all singular points. We can distinguish two situations:
\ben
\item If $\alpha+\beta+\mu$ is even, by walking on the lattice $\cL$ we can reduce ourselves to the
equation for the {\em complete elliptic integral}, which corresponds to $\alpha=\beta=\mu=0$.
\item If $\alpha+\beta+\mu$ is odd, by walking on the lattice $\cL$ 
we  can reduce ourselves to the equation for the Legendre polynomial of degree $0$, which corresponds to
 $\alpha=\beta=0$,
$\mu=1$. This equation is solved by
\begin{eqnarray*}F_{0,0,1}(z)&=F(0,1;1;z)&=1,\\
z^{-1}F_{1,0,0}(1-z^{-1})&=z^{-1}F(1,1;2;1-z^{-1})&=\log(z-1)-\log z,
\end{eqnarray*}
where we used Kummer's table and
\[F(1,1;2;w)=-w^{-1}\log(1-w).\]
\een

\section{The    ${}_1F_1$ and  ${}_2F_0$ equation}
\label{s4}

\subsection{The ${}_1F_1$  equation }
Let $a,c\in\cc$.
The {\em confluent} or the {\em ${}_1F_1$  equation} is given by the operator
\beq
\F(a;c;z,\p_z):=z\p_z^2+(c-z)\p_z-a.
\label{f1c}\eeq
This  equation is a limiting case of the hypergeometric
equation: 
\[\lim_{b\to\infty}
\frac1b\F(a,b;c;z/b,\p_{z/b})=
\F(a;c;z,\p_z).\]

\subsection{The ${}_2F_0$ equation}

Parallel to the ${}_1F_1$ equation we will  consider
the {\em ${}_2F_0$ equation}, given by the operator
\beq
\F(a,b;-;z,\p_z):=z^2\p_z^2+(-1+(1+a+b)z)\p_z+ab,
\label{g8}\eeq
where $a,b\in \cc$.
This
 equation is another limiting case of the hypergeometric
equation: 
\[\lim_{c\to\infty}
\F(a,b;c;cz,\p_{(cz)})=
-\F(a,b;-;z,\p_z).\]


\subsection{Equivalence of the ${}_1F_1$ and  ${}_2F_0$ equation}

Note that
\[\F(a,b;-;z,\p_z)=w^2\p_w^2+(-w^2+(1-a-b)w)\p_w+ab\,\]
where $w=-z^{-1}$, $z=-w^{-1}$.
Moreover,
\beq
(-z)^{a+1}\F(a,b;-;z,\p_z)(-z)^{-a}
=\F(a;1+a-b;w,\p_w).
\label{g8a}\eeq
Hence the ${}_2F_0$ equation is equivalent to the  ${}_1F_1$ 
equation. We will treat the ${}_1F_1$ equation as the principal one.

The relationship between the parameters is
\[c=1+a-b,\ \ \ \ b=1+a-c.\]

\subsection{Lie-algebraic parameters}

Instead of the classical
 parameters we usually prefer  the Lie-algebraic parameters $\alpha,\theta$:
\[\begin{array}{lll}
\alpha:=c-1=a-b,&& \theta: =-c+2a=-1+a+b;\\[2ex]
a=\frac{1+\alpha+\theta }{2},& b=\frac{1 -\alpha+\theta}{2},\ \ \ &
 c=1+\alpha.
\end{array}
\]
In these parameters the  ${}_1F_1$   operator (\ref{f1c}) becomes
\begin{eqnarray*}
\F_{\theta ,\alpha}(z,\p_z)
&=&z\p_z^2+(1+\alpha-z)\p_z-\frac{1}{2}(1+\theta +\alpha)
,\end{eqnarray*}
and the  ${}_2F_0$   operator (\ref{g8}) becomes
\begin{eqnarray*}
\tilde\F_{\theta,\alpha}(z,\p_z)&=&
z^2\p_z^2+(-1+(2+\theta)z)\p_z+\frac14(1+\theta)^2-\frac14\alpha^2.
\end{eqnarray*}

The Lie-algebraic parameters  have an interesting interpretation in terms of a natural basis of  a ``Cartan algebra'' of the Lie algebra $sch(2)$ \cite{DM}.

\subsection{Integral representations}
Two kinds of integral representations of solutions to the 
 ${}_1F_1$ equation
are described below:
\bet\ben \item
Let $[0,1]\ni t\mapsto\gamma(t)$ satisfy
\[t^{a-c+1}\e^t(t-z)^{-a-1}\Big|_{\gamma(0)}^{\gamma(1)}=0.\]
Then
\beq
\F(a;c;z,\p_z)
\int_\gamma t^{a-c}\e^t(t-z)^{-a}\d t=0.\label{dad1}\eeq
\item
Let  $[0,1]\ni t\mapsto\gamma(t)$ satisfy
\[\e^{\frac{z}{t}}t^{-c}(1-t)^{c-a}\Big|_{\gamma(0)}^{\gamma(1)}=0.\]
Then
\beq
\F(a;c;z,\p_z)\int_\gamma\e^{\frac{z}{t}}t^{-c}(1-t)^{c-a-1}\d t
=0.
\label{dad}
\eeq
\een\label{dad4}\eet

\proof
We check that for any contour $\gamma$ the l.h.s of (\ref{dad1}) and (\ref{dad}) equal 
\begin{eqnarray*}
&&
-a\int_\gamma\Big(\p_t
t^{a-c+1}\e^t(t-z)^{-a-1}\Big)\d t,\\
&&-\int_\gamma\Big(\p_t
\e^{\frac{z}{t}}t^{-c}(1-t)^{c-a}\Big)\d t
\end{eqnarray*}
respectively. \qed

For solutions of the  ${}_2F_0$ equation we  also have two kinds of
integral representations:
\bet
Let  $[0,1]\ni t\mapsto\gamma(t)$ satisfy
\[\e^{-\frac{1}{t}}t^{b-a-1}(t-z)^{-b-1}
\Big|_{\gamma(0)}^{\gamma(1)}=0.\]
Then
\beq\F(a,b;-;z,\p_z)\int_\gamma\e^{-\frac{1}{t}}t^{b-a-1}(t-z)^{-b}\d t
\label{dad2}\eeq
\eet

\proof
We check that for any contour $\gamma$ (\ref{dad2}) equals
\[-b\int_\gamma\Big(\p_t\e^{-\frac{1}{t}}t^{b-a-1}(t-z)^{-b-1}\Big)\d t.\]
\qed

The second integral representation is obtained if we interchange $a$ and $b$.

\subsection{Symmetries}
\label{symcom1}

The following operators equal $\F_{\theta ,\alpha}(w,\p_w)$ for the appropriate $w$:
\[\begin{array}{rrcl}
w=z:&&&\\
&&\F_{\theta ,\alpha}(z,\p_z),&\\[1ex]
&z^{-\alpha}&\F_{\theta ,-\alpha}(z,\p_z)&z^{\alpha},\\
w=-z:&&&\\
&-\e^{-z}&\F_{-\theta ,\alpha}(z,\p_z)&\e^z,
\\[1ex]
&-\e^{-z}z^{-\alpha}&\F_{-\theta ,-\alpha}(z,\p_z)&\e^zz^{\alpha}
.\end{array}\label{newnot1}\]
The third symmetry is sometimes called the {\em 1st Kummer transformation}.

Symmetries of the ${}_1F_1$ operators can be interpreted as the ``Weyl group'' of the Lie algebra $sch(2)$.

\subsection{Factorizations  and commutation relations}
\label{symcom1a}
There are several ways of factorizing the ${}_1F_1$ operator.
\begin{eqnarray*}
\F_{\theta,\alpha}&=&\Big(z\partial_z+1+\alpha-z\Big)\partial_z-\frac12(\theta+\alpha+1),\\
&=&\partial_z\Big(z\partial_z+\alpha-z\Big)-\frac12(\theta+\alpha-1),\\
&=&\Big(z\partial_z+1+\alpha\Big)\Big(\partial_z-1\Big)+\frac12(-\theta+\alpha+1),\\
&=&\Big(\partial_z-1\Big)\Big(z\partial_z+\alpha\Big)+\frac12(-\theta+\alpha-1);
\end{eqnarray*}
\begin{eqnarray*}
z\F_{\theta,\alpha}&=&\Big(z\partial_z+\frac12(\theta+\alpha-1)\Big)\Big(
z\partial_z+\frac12(-\theta+\alpha+1)-z\Big)
\\&&-\frac14(-\theta+\alpha+1)(\theta+\alpha-1),\\
&=&\Big(
z\partial_z+\frac12(-\theta+\alpha-1)-z\Big)
\Big(z\partial_z+\frac12(\theta+\alpha+1)\Big)
\\&&-\frac14(-\theta+\alpha-1)(\theta+\alpha+1).
\end{eqnarray*}

One can use the factorizations to derive the following commutation relations: 
\[\begin{array}{rrl}
&\p_z&\F_{\theta ,\alpha}\\[1ex]
&=\ \ \ \F_{\theta +1,\alpha+1}&\p_z,\\[3ex]
&(z\p_z+\alpha-z)&\F_{\theta ,\alpha}\\[1ex]
&=\ \ \ \F_{\theta -1,\alpha-1}&(z\p_z+\alpha-z),\\[3ex]
&(z\p_z+\alpha)&
\F_{\theta ,\alpha}\\[1ex]
&=\ \ \ \F_{\theta +1,\alpha-1}&(z\p_z+\alpha),\\[3ex]
& (\p_z-1)&\F_{\theta ,\alpha},\\[1ex]
&=\ \ \ \F_{\theta -1,\alpha+1}&(\p_z-1);\\[3ex]
&\big( z\p_z+\12(\theta + \alpha+1)\big)&z\F_{\theta ,\alpha}\\[1ex]
&=\ \ \  z\F_{\theta +2,\alpha}&\big( z\p_z+\12(\theta + \alpha+1)\big),\\[3ex]
&\big(z\p_z+\12(-\theta +\alpha+1)-z)
&z\F_{\theta ,\alpha}\\[1ex]
&=\ \ \ z\F_{\theta -2,\alpha}&\big(z\p_z+\12(-\theta +\alpha+1)-z\big).
\end{array}\]
Each of these commutation relations can be associated with
a ``root'' of the Lie algebra 
 $sch(2)$.

\subsection{Canonical forms}

The natural weight of the ${}_1F_{1}$  operator equals $z^{\alpha}\e^{-z}$, so that
\[\F_{\theta,\alpha}=
z^{-\alpha}\e^z\partial_zz^{\alpha+1}\e^{-z}\partial_z
-\frac12(1+\alpha+\theta).\]
The balanced form  of the ${}_1F_{1}$  operator is
\begin{eqnarray*}
z^{\frac{\alpha}{2}}\e^{-\frac{z}{2}}\F_{\theta,\alpha}
z^{-\frac{\alpha}{2}}\e^{\frac{z}{2}}&=&
\partial_zz\partial_z-\frac{z}{4}-\frac{\theta}{2}-\frac{\alpha^2}{4z}.
\end{eqnarray*}

\ber
We have
\begin{eqnarray*}
2z^{\frac{\alpha}{2}-1}\e^{-\frac{z}{2}}\F_{0,\alpha}(z,\partial_z)
z^{-\frac{\alpha}{2}}\e^{\frac{z}{2}}&=&
\partial_w^2+\frac{1}{w}\partial_w
-1-\frac{\alpha^2}{w^2},\ \ \ \ z=2w;\\
2\i z^{\frac{\alpha}{2}-1}\e^{-\frac{z}{2}}\F_{0,\alpha}(z,\partial_z)
z^{-\frac{\alpha}{2}}\e^{\frac{z}{2}}&=&
\partial_u^2+\frac{1}{u}\partial_u
+1-\frac{\alpha^2}{u^2},\ \ \ \ z=2\i u.
\end{eqnarray*}
which are the operators for the {\em modified Bessel} and {\em Bessel equations}.
Thus both these equations  essentially coincide with the balanced form
of the ${}_1F_1$ equation with $\theta=0$. We will discuss them further
in Rem. \ref{bessel0}.
\label{bess}
\eer

The Schr\"odinger form of the ${}_1F_1$ equation is
\begin{eqnarray}
z^{\frac{\alpha}{2}-\frac12}\e^{-\frac{z}{2}}\F_{\theta,\alpha}
z^{-\frac{\alpha}{2}-\frac12}\e^{\frac{z}{2}}&=&
\p_z^2
-\frac{1}{4}-\frac{\theta}{2z}+\Big(\frac14-\frac{\alpha^2}{4}\Big)\frac{1}{z^2}.\label{whitta}
\end{eqnarray}

\ber
In the literature the equation given by (\ref{whitta}) is often called the {\em Whittaker equation}. Its standard form is
\[\p_z^2
-\frac{1}{4}+\frac{\kappa}{z}+\Big(\frac14-\mu^2\Big)\frac{1}{z^2}.\label{witta}
\]
Thus,  $\kappa$, $\mu$ correspond to $-\frac{\theta}{2}$, $\frac{\alpha}{2}$. 
\eer

The natural weight of  the ${}_2F_0$  operator
equals $z^{\theta}\e^{\frac{1}{z}}$, so that
\[\tilde\F_{\theta,\alpha}=
z^{-\theta}\e^{-\frac{1}{z}}\partial_zz^{\theta+2}\e^{\frac{1}{z}}\partial_z
+\frac{(1+\theta)^2}{4}-\frac{\alpha^2}{4}.\]
The balanced form of the ${}_2F_0$  operator is
\begin{eqnarray}
z^{\frac{\theta}{2}}\e^{\frac{1}{2z}}\tilde\F_{\theta,\alpha}
z^{-\frac{\theta}{2}}\e^{-\frac{1}{2z}}
&=&
\partial_zz^2\partial_z-\frac{1}{4z^2}+\frac{\theta}{2z}+\frac{1-\alpha^2}{4}.
\label{besslo}\end{eqnarray}

The symmetries $\alpha\mapsto-\alpha$, as well as $(z,\theta)\mapsto(-z,-\theta)$ are obvious in both balanced forms and in the Whittaker equation.

\subsection{The ${}_1F_1$ function}

Equation  (\ref{f1c}) has a regular singular point at $0$.
Its indices at $0$ are equal $0$, $1-c$.
For $c\neq 0,-1,-2,\dots$, the unique solution of the confluent equation analytic
at $0$  and equal to 1 at 0 is called 
the ${}_1F_1$ hypergeometric function or
the  confluent function. It is equal to
\[F(a;c;z):=\sum_{n=0}^\infty
\frac{(a)_n}{
(c)_n}\frac{z^n}{n!}.\]
It is defined for
$c\neq0,-1,-2,\dots$.
Sometimes it is more convenient to consider
the function
\[ {\bf F}  (a;c;z):=\frac{F(a;c;z)}{\Gamma(c)}=
\sum_{n=0}^\infty
\frac{(a)_n}{
\Gamma(c+n)}\frac{z^n}{n!}.\]
Another useful function proportional to ${}_1F_1$ is
\begin{eqnarray*}
 {\bf F}^\I   (a;c;z)&:=&\frac{\Gamma(a)\Gamma(c-a)}{\Gamma(c)}F(a;c;z).
\end{eqnarray*}

The confluent function can be obtained as the limit of the hypergeometric
function: 
\[F(a;c;z)=\lim_{b\to\infty}F(a,b;c;z/b).\]

It satisfies the so-called {\em Kummer's identity}:
\beq F(a;c;z)=\e^z F\left(c-a;c;-z\right).
\eeq

Integral representations for all parameters
\begin{eqnarray*}
\frac{1}{2\pi \i}\int\limits
_{]-\infty,(0,z)^+,-\infty[} t^{a-c}\e^t(t-z)^{-a}\d t
&=& {\bf F}  (a;c;z),\end{eqnarray*}
for $\Re a>0,\ \Re (c-a)>0$
\begin{eqnarray*}\int\limits_{[1,+\infty[}\e^{\frac{z}{t}}t^{-c}(t-1)^{c-a-1}\d t
&=& {\bf F}^\I   (a;c;z),\end{eqnarray*}
and for $\Re (c-a)>0$
\begin{eqnarray}\frac{1}{2\pi\i}\int\limits_{[1,0^+,1]}\e^{\frac{z}{t}}(-t)^{-c}(-t+1)^{c-a-1}\d t
&=& \frac{\sin\pi a}{\pi}{\bf F}^{\I}   (a;c;z).
\label{kumme}\end{eqnarray}

In the Lie-algebraic parameters:
\begin{eqnarray*} 
F_{\theta ,\alpha}(z)&:=&F\Bigl(\frac{1+\alpha+\theta }{2};1+\alpha;z\Bigr)
,\\
 {\bf F}  _{\theta ,\alpha}(z)&:=&
 {\bf F}  \Bigl(\frac{1+\alpha+\theta }{2};1+\alpha;z\Bigr)\\
&=&
 \frac{1}{\Gamma(\alpha+1)}F_{\theta ,\alpha}(z),\\
 {\bf F}^\I   _{\theta ,\alpha}(z)&:=&
 {\bf F}^\I   \Bigl(\frac{1+\alpha+\theta }{2};1+\alpha;z\Bigr)\\&=&
\frac{\Gamma(\frac{1+\alpha+\theta}{2})
\Gamma(\frac{1+\alpha-\theta}{2})}{\Gamma(\alpha+1)}
F_{\theta ,\alpha}(z).\end{eqnarray*}

\ber In the literature the ${}_1F_1$ function is often called {\em Kummer's function} and denoted
\[M(a,c,z):=F(a;c;z).\]
One also  uses the {\em Whittaker function of the 1st kind}
\[M_{\kappa,\mu}(z):=\exp(-z/2)z^{\mu+1/2}M\Big(\mu-\kappa+\frac12,1+2\mu,z\Big),\]
which solves the Whittaker equation. \eer

\subsection{The ${}_2F_0$
 function}

We define, for $z\in\cc\backslash[0,+\infty[$,
\[F(a,b;-;z):=\lim_{c\to\infty}F(a,b;c;cz),\]
where $|\arg c-\pi|<\pi-\epsilon$, $\epsilon>0$.
It extends to an analytic function on the universal cover of
$\cc\backslash\{0\}$ 
with a branch point of an infinite order at 0.
It has the following asymptotic expansion:
\[
F(a,b;-;z)\sim\sum_{n=0}^\infty\frac{(a)_n(b)_n}{n!}z^n,
\ |\arg z-\pi|<\pi-\epsilon.
\]
Sometimes instead of ${}_2F_0$ it is useful to consider the function
\begin{eqnarray*}
 {\bf F}^\I   (a,b;-;z)&:=&\Gamma(a)F(a,b;-;z).
\end{eqnarray*}

We have an integral representation for $\Re a>0$
\[\int_0^\infty
\e^{-\frac{1}{t}}t^{b-a-1}(t-z)^{-b}\d t
= {\bf F}^\I   (a,b;-;z), \ \ z\not\in[0,\infty[,
\]and without a restriction on parameters
\[\frac{1}{2\pi\i}\int\limits_{[0,z^+,0]}
\e^{-\frac{1}{t}}t^{b-a-1}(t-z)^{-b}\d t
= \frac{\sin\pi a}{\pi}{\bf F}^\I   (a,b;-;z), \ \ z\not\in[0,\infty[.
\]

When we use the Lie-algebraic parameters, we denote the ${}_2F_0$ function by
$\tilde F$ and $\tilde {\bf F}$. The tilde is needed to  avoid the confusion with the ${}_1F_1$ function:
\begin{eqnarray*}
\tilde F_{\theta ,\alpha}(z)&:=&F\Bigl(\frac{1+\alpha+\theta }{2},\frac{1-\alpha+\theta }{2};-;z\Bigr),\\
\tilde {\bf F}^\I  _{\theta ,\alpha}(z)&:=& {\bf F}^\I  \Bigl(\frac{1+\alpha+\theta }{2},
\frac{1 -\alpha+\theta}{2};-;z\Bigr)\\
&=& \Gamma\Big(\frac{1-\alpha+\theta}{2}\Big)\tilde F_{\theta ,\alpha}(z)  .
\end{eqnarray*}

\ber In the literature the ${}_2F_0$ function is seldom used.
Instead one uses {\em Tricomi's function} 
\[U(a,c,z):=z^{-a}F(a;a-b-1;-;z^{-1}).\]
It is one of  solutions of the ${}_1F_1$ equation, which we will discuss  in
Subsubsect \ref{s4.6-}.
One also  uses the {\em Whittaker function of the 2nd kind}
\[W_{\kappa,\mu}(z):=\exp(-z/2)z^{\mu+1/2}U\Big(\mu-\kappa+\frac12;1+2\mu;z\Big),\]
which solves the Whittaker equation. \eer

\subsection{Standard solutions}

The ${}_1F_1$ equation has two singular points. $0$ is a regular
singular point and with each of its two indices we can associate the
corresponding solution. $\infty$ is not a regular singular point. However  we can define two solutions with a simple behavior around $\infty$. Altogether we obtain 4 {\em standard solutions}, which we will describe in this subsection.

It follows by Thm \ref{dad4} that, for appropriate contours $\gamma_1$, $\gamma_2$, the integrals
\begin{eqnarray*}
&&\int\limits
_{\gamma_1}
t^{\frac{-1+\theta -\alpha}{2}}\e^t(t-z)^{\frac{-1-\theta -\alpha}{2}}\d t,
\\
&&
\int\limits_{\gamma_2}\e^{\frac{z}{t}}t^{-1-\alpha}(t-1)^{\frac{-1-\theta +\alpha}{2}}
\d t
\end{eqnarray*}
solve the
 ${}_1F_1$ equation. 

In the first integral  the natural candidates for the endpoints of the intervals of integration are $\{-\infty,0,z\}$. We will see that all 4 standard solutions can be obtained as such integrals.

In the second integral the natural candidates for endpoints are
$\{1,0-0,\infty\}$. (Recall from Subsect \ref{a.3} that $0-0$ denotes $0$ approached from the left). The 4 standard solutions can  be obtained  also from the integrals with these endpoints.

\subsubsection{Solution $\sim1$ at $0$}

For $\alpha\neq-1,-2,\dots$, the only solution $\sim 1$ around $0$
is
\begin{eqnarray*}
F_{\theta ,\alpha}(z)&=&\e^zF_{-\theta ,\alpha}(-z).\end{eqnarray*}

The first integral representation is valid for all parameters:
\begin{eqnarray*}
\frac{1}{2\pi \i}\int\limits
_{]-\infty,(0,z)^+-\infty[}
t^{\frac{-1+\theta -\alpha}{2}}\e^t(t-z)^{\frac{-1-\theta -\alpha}{2}}\d t
&=& {\bf F}  _{\theta ,\alpha}(z).
\end{eqnarray*}
The second is valid for $\Re(1+\alpha)>|\Re \theta |$:
\begin{eqnarray*}
\int\limits_{[1,+\infty[}\e^{\frac{z}{t}}t^{-1-\alpha}(t-1)^{\frac{-1-\theta +\alpha}{2}}
\d t
&=& {\bf F}^\I   _{\theta ,\alpha}(z).\end{eqnarray*}

\subsubsection{Solution $\sim z^{-\alpha}$ at $0$}
\label{s4.6}

If $\alpha\neq1,2,\dots$, then the
 only solution of the confluent equation behaving as
 $z^{-\alpha}$  at
 $0$ is equal to
\begin{eqnarray*}
z^{-\alpha} F _{\theta ,-\alpha}(z)&=&z^{-\alpha}
\e^z F _{-\theta ,-\alpha}(-z).\end{eqnarray*}

Integral representation 
for $\Re(1-\alpha)>|\Re \theta |$:
\begin{eqnarray*}
\int_0^z
t^{\frac{-1+\theta -\alpha}{2}}\e^t(z-t)^{\frac{-1-\theta -\alpha}{2}}\d t
&=&
z^{-\alpha} {\bf F}^\I   _{\theta ,-\alpha}(z),\ \ z\not\in]-\infty,0];\\
\int_z^0
(-t)^{\frac{-1+\theta -\alpha}{2}}\e^t(t-z)^{\frac{-1-\theta -\alpha}{2}}\d t
&=&
(-z)^{-\alpha} {\bf F}^\I   _{\theta ,-\alpha}(z),\ \ z\not\in[0,\infty[;
\end{eqnarray*}
and without a restriction on parameters:
\begin{eqnarray*}
\frac{1}{2\pi\i}
\int\limits_{(0-0)^+}
\e^{\frac{z}{t}}t^{-1-\alpha}(1-t)^{\frac{-1-\theta +\alpha}{2}}
\d t
&=&z^{-\alpha} {\bf F}  _{\theta ,-\alpha}(z),\ \ \Re z>0.
\end{eqnarray*}

\subsubsection{Solution $\sim z^{-a}$  at $+\infty$ }
\label{s4.6-}

The following solution to the confluent equation
behaves as   $\sim z^{-a}=z^{-\frac{1+\theta +\alpha}{2}}$ 
 at $+\infty$ for $|\arg z|<\pi-\epsilon$:
\begin{eqnarray*}
z^{\frac{-1-\theta -\alpha}{2}}\tilde F_{\theta ,\pm \alpha}(-z^{-1}).
&&
\end{eqnarray*}

Integral representations for $\Re(1+\theta -\alpha)>0$:
\begin{eqnarray*}
\int_{-\infty}^0
(-t)^{\frac{-1+\theta -\alpha}{2}}\e^t(z-t)^{\frac{-1-\theta -\alpha}{2}}\d t
&=&
z^{\frac{-1-\theta -\alpha}{2}} \tilde {\bf F}^\I  _{\theta ,\alpha}(-z^{-1}),\ \ \ \ 
\ \ z\not\in]-\infty,0];\end{eqnarray*}
and, for $ \Re(1+\theta +\alpha)>0$:
\begin{eqnarray*}\int_{-\infty}^0
\e^{\frac{z}{t}}(-t)^{-1-\alpha}(1-t)^{\frac{-1+\theta +\alpha}{2}}
\d t&=&
z^{\frac{-1-\theta -\alpha}{2}}
\tilde {\bf F}^\I  _{\theta ,-\alpha}(-z^{-1}),\ \ \ \ \Re z>0.
\end{eqnarray*}

\subsubsection{Solution $\sim  (-z)^{-b}\e^z$ at $-\infty$ }

The following solution to the confluent equation
behaves as   $\sim (-z)^{-b}\e^z= (-z)^{\frac{1+\theta -\alpha}{2}}\e^z$ 
 at $\infty$ for $|\arg z-\pi|<\pi-\epsilon$:
\begin{eqnarray*}
\e^z(-z)^{\frac{-1-\theta -\alpha}{2}}\tilde F_{-\theta ,\pm\alpha}(z^{-1})
.&&
\end{eqnarray*}

Integral representation for $\Re(1+\theta +\alpha)>0$:
\begin{eqnarray*}
\int_{-\infty}^z
(-t)^{\frac{-1+\theta -\alpha}{2}}\e^t(z-t)^{\frac{1-\theta -\alpha}{2}}\d t
&=&
\e^z(-z)^{\frac{-1-\theta -\alpha}{2}} \tilde {\bf F}^\I  _{-\theta ,-\alpha}(z^{-1}),\ \ \ \ \ z\not\in[0,\infty[;\end{eqnarray*}
and for $\Re(1+\theta -\alpha)>0$:
\begin{eqnarray*}\int_0^1
\e^{\frac{z}{t}}t^{-1-\alpha}(1-t)^{\frac{-1+\theta +\alpha}{2}}
\d t&=&
\e^z(-z)^{\frac{-1-\theta -\alpha}{2}}
\tilde {\bf F}^\I  _{-\theta ,\alpha}(z^{-1}),\ \ \ \Re z<0.
\end{eqnarray*}

\subsection{Connection formulas}

We decompose  standard solutions in pair of solutions with a simple behavior around zero.

\begin{eqnarray*}
z^{\frac{-1-\theta -\alpha}{2}}\tilde F_{\theta ,\pm \alpha}(-z^{-1})
&=&\frac{\pi}{\sin{\pi(-\alpha)}\Gamma\left(\frac{1+\theta -\alpha}{2}\right)}
 {\bf F}  _{\theta ,\alpha}(z)\\
&&+\frac{\pi}{\sin\pi \alpha\Gamma\left(\frac{1+\theta +\alpha}{2}\right)}
z^{-\alpha} {\bf F}  _{\theta ,-\alpha}(z),\\
\e^z(-z)^{\frac{-1-\theta -\alpha}{2}}\tilde F_{-\theta ,\pm\alpha}(z^{-1})
&=&\frac{\pi}{\sin{\pi(-\alpha)}
\Gamma\left(\frac{1-\theta -\alpha}{2}\right)}
 {\bf F}  _{\theta ,\alpha}(z)\\
&&+\frac{\pi}
{\sin\pi \alpha\Gamma\left(\frac{1-\theta +\alpha}{2}\right)}
(-z)^{-\alpha} {\bf F}_{\theta ,-\alpha}(z).
\end{eqnarray*}

\subsection{Recurrence relations}
The following recurrence relations follow easily
from the commutation relations of Subsect. \ref{symcom1a}:

\begin{eqnarray*}
,\\[3ex]
 \p_z {\bf F}  _{\theta ,\alpha}(z)&=&\frac{1+\theta +\alpha}{2} {\bf F}  _{\theta +1,\alpha+1}(z),
 \\
 \left(z\p_z+\alpha-z\right) {\bf F}  _{\theta ,\alpha}(z)&=& {\bf F}  _{\theta -1,\alpha-1}(z),\\[3ex]
 \left(z\p_z+\alpha\right) {\bf F}  _{\theta ,\alpha}(z)&=& {\bf\rm
   F}  _{\theta +1,\alpha-1}(z),\\
 \left(\p_z-1\right) {\bf F}  _{\theta ,\alpha}(z)&=&\frac{-1+\theta -\alpha}{2} {\bf F}  _{\theta -1,\alpha+1}(z),
\\[3ex]
 \left(z\p_z+\frac{1+\theta +\alpha}{2}\right) {\bf F}  _{\theta ,\alpha}(z)&=&\frac{1+\theta +\alpha}{2} {\bf F}  _{\theta +2,\alpha}(z),
\\
 \left(z\p_z+\frac{1-\theta +\alpha}{2}-z\right) {\bf F}  _{\theta ,\alpha}(z)&
=&\frac{1-\theta +\alpha}{2} {\bf F}  _{\theta
  -2,\alpha}(z).
\end{eqnarray*}

The recurrence relations for the ${}_2F_0$ functions
are similar:
\begin{eqnarray*}
\left(z\p_z+\frac{1+\theta +\alpha}{2}\right)\tilde {\bf F}_{\theta ,\alpha}^\I(z)&=&\frac{1+\theta +\alpha}{2}
\tilde {\bf F}_{\theta +1,\alpha+1}^\I(z),\\
\left(z^2\p_z-1+\frac{1+\theta -\alpha}{2}z\right)\tilde {\bf F}_{\theta ,\alpha}^\I(z)&=&-\tilde {\bf F}_{\theta -1,\alpha-1}^\I(z),\\[2ex]
\left(z\p_z+\frac{1+\theta -\alpha}{2}\right)\tilde {\bf F}_{\theta ,\alpha}^\I(z)&=&
\tilde {\bf F}_{\theta +1,\alpha-1}^\I(z),\\
\left(z^2\p_z-1+\frac{1+\theta +\alpha}{2}z\right)\tilde {\bf F}_{\theta ,\alpha}^\I(z)&=&
\frac{1-\theta+\alpha}{2}\tilde {\bf F}_{\theta -1,\alpha+1}^\I(z),\\[2ex]
\p_z\tilde{\bf F}_{\theta ,\alpha}^\I(z)&=&\frac{1{+}\theta {+}\alpha}{2}
\tilde {\bf F}_{\theta +2,\alpha}^\I(z),\\
(z^2\p_z-1-\theta z)\tilde {\bf F}_{\theta ,\alpha}^\I(z)&=&\frac{1-\theta+\alpha}{2}\tilde {\bf F}_{\theta -2,\alpha}^\I(z).
\end{eqnarray*}
\subsection{Additional recurrence relations}

There exists an additional pair of recurrence relations:
\begin{eqnarray*}
\left((1{-}\alpha)z^2\p_z{+}\frac{(1{-}\alpha)(1-\alpha+\theta)}{2}z{+}\frac{{-}1{-}\theta {+}\alpha}{2}\right)
\tilde F_{\theta ,\alpha}(z)&=&\frac{-1{-}\theta {+}\alpha}{2}\tilde F_{\theta ,\alpha-2}(z),\\[2ex]
\left(
(1{+}\alpha)z^2\p_z{+}\frac{(1{+}\alpha)(1+\alpha+\theta)}{2}z{+}\frac{{-}1{-}\theta {-}\alpha}{2}\right)
\tilde F_{\theta ,\alpha}(z)&=&
\frac{{-}1{-}\theta {-}\alpha}{2}\tilde F_{\theta ,\alpha+2}(z).
\end{eqnarray*}

\subsection{Degenerate case}

$\alpha=m\in\zz$ is the degenerate case of the confluent equation at $0$.
We have then
\[{\bf F}
(a;1+m;z)=\sum_{n=\max(0,-m)}\frac{(a)_n}{n!(m+n)!}z^n.\]
This easily implies the identity
\[(a-m)_m{\bf F}(a;1+m;z)=z^{-m}
{\bf F}(a-m;1-m;z).
\]Thus the two standard solutions determined by the behavior at zero are proportional to one another. 

One can also see the degenerate case in the integral representations: 
\begin{eqnarray*}\frac{1}{2\pi\i}\int\limits
_{[(z,0)^+]}\e^{t}(1-z\slash t)^{-a}
t^{-m-1}\d t&=&
 {\bf F}  _{-1+2a-m ,m}(z)\\
&=&
\left(a\right)_{-m} z^{-m} {\bf F}  _{-1+2a-m ,-m}(z),\\[4ex]
\frac{1}{2\pi\i}\int\limits_{[(0,1)^+]}
\e^{z\slash t}(1-t)^{-a}t^{-m-1}\d t
&=&\left(a\right)_{m}
 {\bf F}  _{-1+2a+m ,m}(z)\\
&=&
 z^{-m} {\bf F}  _{-1+2a+m ,-m}(z)
.\end{eqnarray*}
The corresponding  generating functions are
\begin{eqnarray*}
\e^{t}(1-z\slash t)^{-a}&=&
\sum_{m\in\zz}t^m {\bf F}  _{-1+2a-m,m}(z),\\
\e^{z\slash t}(1-t)^{-a}&=&
\sum_{m\in\zz}t^m(a)_{m} {\bf F}  _{-1+2a+m,m}(z).\end{eqnarray*}



\subsection{Laguerre polynomials}
\label{s9.4}

${}_1F_1$ functions for $-a=n=0,1,2,\dots$ are polynomials. 
They are
known as  {\em Laguerre polynomials}.

Following
Subsect. 
\ref{Hypergeometric type polynomials}, they can be defined by the following version of the Rodriguez-type formula:
\[
L_n^\alpha(z):=\frac{1}{n!}\e^zz^{-\alpha}\p_z^n\e^{-z}z^{n+\alpha}.
\]
The differential equation:
\begin{eqnarray*}
\F(-n;\alpha+1;z,\p_z)L_n^\alpha(z)&&\\
=\ \Big(z\partial_z^2+(1+\alpha-z)\partial_z+n\Big)L_n^\alpha(z)&=&0.
\end{eqnarray*}
Generating functions:
\[\begin{array}{l}
\e^{-tz}(1+t)^{\alpha}
=\sum\limits_{n=0}^\infty t^n L_n^{\alpha-n}(z),\\[5mm]
(1-t)^{-\alpha-1}
\exp{\frac{tz}{t-1}}=\sum\limits_{n=0}^\infty t^nL_n^{\alpha}(z).
\end{array}\]
Integral representations:
\[\begin{array}{rl}
L_n^\alpha(z)&=\frac{1}{2\pi\i}\int\limits_{[0^+]}
\e^{-tz}(1+t)^{\alpha+n}t^{-n-1}\d t\\[5mm]
&=\frac{1}{2\pi\i}\int\limits_{[0^+]}
(1-t)^{-\alpha-1}\exp(\frac{tz}{t-1})t^{-n-1}\d t.
\end{array}\]
Expression in terms of the Bessel polynomials (to be defined in the next subsection):
\[\begin{array}{rl}
L_n^\alpha(z)&=z^nB_n^{-2n-\alpha-1}(-z^{-1}).
\end{array}\]
Recurrence relations:
\begin{eqnarray*}
\p_zL_n^\alpha(z)&=&-L_{n-1}^{\alpha+1}(z),\\
\left(z\p_z+\alpha-z\right)L_n^\alpha(z)&=&(n+1)L_{n+1}^{\alpha-1}(z),\\[3ex]
\left(z\p_z+\alpha\right)L_n^\alpha(z)&=&(\alpha+n)L_n^{\alpha-1}(z),\\
\left(\p_z-1\right)L_n^\alpha(z)&=&-L_n^{\alpha+1}(z),
\\[3ex]
\left(z\p_z-n\right)L_n^\alpha(z)&=&-(n+\alpha)L_{n-1}^\alpha(z),\\
\left(z\p_z+n+\alpha+1-z\right)L_n^\alpha(z)&=&(n+1)L_{n+1}^\alpha(z).
\end{eqnarray*}

The first, resp. second integral representation is easily seen to be
equivalent to the first, resp. second generating function.

The differential equation, the Rodriguez-type formula, the first
generating function, the first integral representation and the first
pair of recurrence relations are 
 special cases of the  corresponding formulas of Subsect.  
\ref{Hypergeometric type polynomials}.

We have several alternative expressions for Laguerre polynomials:
\begin{eqnarray*}
L_n^\alpha(z)&=&\lim\limits_{\nu\to n}
(-1)^n(\nu-n)
{\bf F}_{1+\alpha-2\nu,\alpha}^{\I}(z)
=\frac{(1+\alpha)_n}{n!}F(-n;1+\alpha;z)\\
&=&z^n\lim\limits_{\nu\to n}(\nu-n)\tilde{\bf F}_{1+\alpha-2\nu,\alpha}^{\I}(z)
=\frac{1}{n!}(-z)^nF(-n,-\alpha-n;-;-z^{-1})
\\
&=&\sum\limits_{j=0}^n\frac{(1+\alpha+j)_{n-j}}{j!(n-j)!}(-z)^j.
\end{eqnarray*}

Let us derive the above identity using the integral representation
(\ref{kumme}).
Using that $a$ is an integer we can replace the open curve $[1,0^+,1]$
with a closed loop $[\infty^-]$:
\begin{eqnarray*}&&
\lim\limits_{\nu\to n}
(-1)^n(\nu-n)
{\bf F}_{1+\alpha-2\nu,\alpha}^{\I}(z)\\
&=&
\lim\limits_{\nu\to n}\frac{\sin \nu\pi}{\pi}{\bf F}_{1+\alpha-2\nu,\alpha}^{\I}(z)\\
&=&\frac{1}{2\pi\i}
\int_{[\infty^-]}\e^{\frac{z}{s}} (-s)^{-1-\alpha}(1-s)^{\alpha+n}\d s.
\end{eqnarray*}
Then we set $s=-\frac{1}{t}$, resp. $s=1-\frac1t$ to obtain the integral representations.

 The 
value at 0 and behavior at $\infty$:
\[
L_n^\alpha(0)=\frac{(\alpha+1)_n}{n!},\ \  \ \ 
\lim\limits_{z\to\infty}\frac{L_n^\alpha(z)}{z^n}=\frac{(-1)^n}{n!}.
\]

An additional identity valid in the degenerate case:
\begin{eqnarray*}
L_n^{\alpha}(z)&=&(n+1)_\alpha(-z)^{-\alpha}
L_{n+\alpha}^{-\alpha}(z),\ \ \alpha\in\zz.
\end{eqnarray*}

\subsection{Bessel polynomials}

The ${}_2F_0$ functions for $-a=n=0,1,2,\dots$ are
polynomials. Appropriately normalized they are called {\em Bessel
  polynomials}. They are seldom used in the literature, because they
do not form an orthonormal basis in any weighted space and they are
easily expressed in terms of Laguerre polynomials.

Following
Subsect. 
\ref{Hypergeometric type polynomials}, they can be defined by the following version of the Rodriguez-type formula:
\[\begin{array}{l}
B_n^\theta(z):=\frac{1}{n!}
z^{-\theta}\e^{z^{-1}}\p_z^n\e^{-z^{-1}}z^{\theta+2n}.
\end{array}\]
Differential equation:
\begin{eqnarray*}
\F(-n,n+\theta+1;-;\p_z,z)B_n^\theta(z)&&\\
=\Big(z^2\p_z^2+(-1+(2+\theta)z)\p_z-\frac12n(1+\theta-\alpha)\Big)B_n^\theta(z)&=&0.
\end{eqnarray*}
Generating functions:
\[\begin{array}{l}
\e^{-t}(1-tz)^{-\theta-1}
=\sum\limits_{n=0}^\infty t^nB_n^{\theta-n}(z),\\[3mm]
(1+tz)^{\theta}\exp(\frac{-t}{1+tz})=\sum\limits_{n=0}^\infty
t^nB_n^{\theta-2n}(z).\end{array}\]
Integral representations:
\[\begin{array}{rl}
B_n^\theta(z)&=\frac{1}{2\pi\i}
\int\limits_{[0^+]}\e^{t}(1-tz)^{-\theta-n-1}t^{-n-1}\d t\\[3mm]
&=\frac{1}{2\pi\i}\int\limits_{[0^+]}
(1+tz)^{\theta+2n}\exp(\frac{-t}{1+tz})t^{-n-1}\d t.
\end{array}\]
Expression in terms of the Laguerre polynomials:
\[\begin{array}{l}
B_n^\theta(z)=(-z)^nL_n^{-\theta-2n-1}(-z^{-1}).
\end{array}\]
Recurrence relations:
\begin{eqnarray*}
\left(z\p_z+n+\theta+1\right)B_n^\theta(z)&=&(n+\theta+1)B_n^{\theta+1}(z),
\\
\left(z^2\p_z-1-nz\right)B_n^\theta(z)&=&-B_n^{\theta-1}(z),\\[3ex]
\left( z\p_z-n\right)B_n^\theta(z)&=&-B_{n-1}^{\theta+1}(z),\\
\left(z^2\p_z-1+(n+\theta+1)z\right)B_n^\theta(z)&=&-(n+1)B_{n+1}^{\theta-1}(z),
\\[3ex] 
\p_zB_n^\theta(z)&=&-(n+\theta+1)B_{n-1}^{\theta+2}(z),\\
\left(z^2\p_z-1-\theta z\right)B_n^\theta(z)&=&-(n+1)B_{n+1}^{\theta-2}(z).
\end{eqnarray*}

Most of the above identities can be directly obtained from the corresponding identities about Laguerre polynomials.

The  differential equation, the Rodriguez-type formula, the second
generating function, the second integral representation and the last
pair of recurrence relations are  special cases of the  corresponding formulas of Subsect.  
\ref{Hypergeometric type polynomials}.

We have several alternative expressions for Bessel polynomials:
\begin{eqnarray*}
B_n^\theta(z)&=&\lim\limits_{\nu\to n}(-1)^n(\nu-n)\tilde{\bf F}_{\theta,-1-\theta-2n}^{\I}(z)=\frac{1}{n!}F(-n,n+\theta+1;-;z)\\
&=&z^n\lim\limits_{\nu\to n}(\nu-n){\bf F}_{\theta,-1-\theta-2\nu}^{\I}(-z^{-1})\\&
=&\frac{(1+\theta+n)_n}{n!}(-z)^nF(-n;-\theta-2n;-z^{-1}).
\end{eqnarray*}
The value  at zero and behavior at $\infty$:
\[
B_n^\theta(0)=\frac{1}{n!},\ \ \ \ 
\lim\limits_{z\to\infty}\frac{B_n^\theta(z)}{z^n}=\frac{(-1)^n(n+\theta+1)_n}{n!}.\]

Both for Laguerre and Bessel polynomials there exist additional recurrence relations and a generating function. Below we give a pair of such recurrence relations for Bessel polynomials.
\begin{eqnarray*}
\Big((2+2n+\theta)z^2\p_z+(2+2n+\theta)(n+\theta+1)z &&\\-(n+\theta+1)\Big)B_n^\theta(z)
&=&-(n+1)(n+\theta+1)B_{n+1}^\theta(z),\\[2ex]
\Big(-(2n+\theta)z^2\p_z+(2n+\theta)nz+n\Big)B_n^\theta(z)&=&
B_{n-1}^\theta(z).\end{eqnarray*}
They correspond to an additional generating function
\[\begin{array}{rl}2^\theta r^{-1}(1+r)^{-\theta}\exp(\frac{2t}{1+r})
&=\sum\limits_{n=0}^\infty t^n B_n^\theta(z),\\[3mm]
&\hbox{where}\ \ \ \ r:=\sqrt{1+4zt}.
\end{array}\]

\subsection{Special cases}

Apart from the polynomial case and the degenerate case,
 the confluent equation has  some other   cases with special properties.

\subsubsection{Bessel equation}

If $\theta=0$, the confluent equation is equivalent to the (modified) Bessel equation, which we already remarked in Rem. \ref{bess}. By a square root substitution, it is also equivalent to the ${}_0F_1$ equation; see (\ref{gas}).

\subsubsection{Hermite equation}

If $\alpha=\pm\frac12$, the confluent equation is equivalent to the Hermite equation by the quadratic substitutions (\ref{ha6a}) and  (\ref{ha6}).

\section{The ${}_0F_1$  equation}
\label{sa5}


\subsection{Introduction}
Let $c\in\cc$. In this section we will consider
the {\em ${}_0F_1$ equation} given by the operator
\[\F(c;z,\p_z):=z\p_z^2+c\p_z-1.\]
It is a limiting case of the ${}_1F_1$ and ${}_2F_1$ operator:
\[\lim_{a,b\to\infty}\frac1{ab}\F(a,b;c;z/ab,\p_{(z/ab)})=\lim_{a\to\infty}\frac1a\F(a;c;z/a,\p_{(z/a)})=\F(c;z,\p_z).\]

Instead of $c$ it is often more natural to use its {\em Lie-algebraic parameter }
\beq \alpha:=c-1,\ \ \ c=\alpha+1.\label{newnot3}\eeq
Thus we obtain the operator
\begin{eqnarray*}
\F_\alpha (z,\p_z)&:=&z\p_z^2+(\alpha +1)\p_z-1.
\end{eqnarray*}

The Lie-algebraic parameter  has well-known interpretation in terms of the ``Cartan element'' of the Lie algebra $aso(2)$, \cite{V,Wa,DM}.


\subsection{Equivalence with a subclass of the confluent equation}

The ${}_0F_1$ equation can be reduced to a special class of 
the confluent equation by the so-called {\em Kummer's 2nd transformation}:
\begin{eqnarray}
\F(c;z,\p_z)
=\frac{4}{w}\e^{-w/2}\F\Big(c-\12;2c-1;w,\p_w\Big)\e^{w/2},
\label{gas}\end{eqnarray}
where $w=\pm 4\sqrt{z}$, $z=\frac{1}{16}w^2$.
Using the Lie-algebraic parameters this can be rewritten as
\beq
\F_\alpha (z,\p_z)
=\frac{4}{w}\e^{-w/2}\F_{0,2\alpha }(w,\p_w)\e^{w/2}.\label{gas1}\eeq

\subsection{Integral representations}

There are two kinds of integral representations of solutions to the
${}_0F_1$ equation. Thm \ref{schl} describes representations of the first kind, which
 will be called {\em
  Bessel-Schl\"afli representations}. They will be treated as the main ones.
\bet\label{schl}
Suppose that $[0,1]\ni t\mapsto\gamma(t)$ satisfies
\[\e^t\e^{\frac{z}{t}}t^{-c}\Big|_{\gamma(0)}^{\gamma(1)}=0.\]
Then
\beq\F(c;z,\p_z)
\int_\gamma\e^t\e^{\frac{z}{t}}t^{-c}\d t=0.\label{dad5}\eeq
\eet

\proof
We check that for any contour $\gamma$ (\ref{dad5}) equals
\[-\int_\gamma\Big(
\p_t\e^t\e^{\frac{z}{t}}t^{-c}\big)\d t.\]
\qed


Integral representations that can be derived from the
representations for the confluent equation by  2nd Kummer's identity will be called  {\em
  Poisson-type representations}. They will be treated as secondary ones.
They are described in the following theorem.
\bet\ben\item Let the contour $\gamma$ satisfy
\[(t^2-z)^{-c+3/2}\e^{2t}\Big|_{\gamma(0)}^{\gamma(1)}=0.\]
Then
\[\F(c;z,\p_z)\int_\gamma(t^2-z)^{-c+1/2}\e^{2t}\d t=0.\]
\item
 Let the contour $\gamma$ satisfy
\[(t^2-1)^{c-1/2}\e^{2t\sqrt z}\Big|_{\gamma(0)}^{\gamma(1)}=0.\]
Then
\[\F(c;z,\p_z)\int_\gamma(t^2-1)^{c-3/2}\e^{2t\sqrt z}\d t=0.\]\een
\label{dad7}\eet

\proof  By (\ref{gas}) and (\ref{dad1}), for appropriate contours $\gamma$ and $\gamma'$,
\begin{eqnarray*}
&&\e^{-2\sqrt z}\int_\gamma\e^s s^{-c+\frac12}(s-4\sqrt z)^{-c+\frac12}\d s\\
&=&2^{-2c+2}\int_{\gamma'}\e^{2t}(t^2-z)^{-c+\frac12}\d t
\end{eqnarray*}
is annihilated by  $\F(c)$, where we set
$t=\frac{s}{2}-\sqrt z$.  This proves 1.

By (\ref{gas}) and (\ref{dad}), for appropriate contours 
$\gamma$ and $\gamma'$,
\begin{eqnarray*}
&&\e^{-2\sqrt z}\int_\gamma\e^{\frac{4\sqrt z}{s}} s^{-2c+1}(1-s)^{c-\frac32}\d s\\
&=&-2^{-2c+2}\int_{\gamma'}\e^{2t\sqrt z}(1-t^2)^{c-\frac32}\d t
\end{eqnarray*}
is annihilated by  $\F(c)$, where we set
$t=\frac2s-1$.  This proves 2.
\qed

\subsection{Symmetries}
\label{symcom2}

The only nontrivial symmetry is
\[\begin{array}{lcr}z^{-\alpha }\ \F_{-\alpha }\ z^\alpha  &=&\F_\alpha .\end{array}\]
It can be interpreted as  a ``Weyl symmetry'' of $aso(2)$.

\subsection{Factorizations and  and commutation relations}
\label{symcom2a}

There are two ways to factorize the ${}_0F_1$ operator:
\begin{eqnarray*}
\F_\alpha&=&
\big(z\partial_z+\alpha+1\big)\partial_z-1\\
&=&\p_z
\big(z\partial_z+\alpha\big)-1.
\end{eqnarray*}

The factorizations can be used to derive
the following commutation relations:
\[\begin{array}{rl}
\p_z&\F_\alpha \\[1ex]
=\ \ \F_{\alpha +1}&\p_z,\\[3ex]
(z\p_z+\alpha )&\F_\alpha \\[1ex]
=\ \ \F_{\alpha -1}&(z\p_z+\alpha ).\end{array}\]
 Each commutation relation can be associated with a ``root'' of   the Lie
algebra $aso(2)$.

\subsection{Canonical forms}

The natural weight of the ${}_0F_{1}$  operator is $z^\alpha$, so that
\[\F_{\alpha}=
z^{-\alpha}\partial_zz^{\alpha+1}\partial_z
-1.\]
The balanced form of the ${}_0F_{1}$  operator is
\begin{eqnarray*}
z^{\frac{\alpha}{2}}\F_{\alpha}
z^{-\frac{\alpha}{2}}&=&
\partial_zz\partial_z-1-\frac{\alpha^2}{4z}.
\end{eqnarray*}

The symmetry $\alpha\to-\alpha$ is obvious in the balanced form.

\ber\label{bessel0}
In the literature, the 
 ${}_0F_1$  equation is seldom used. Much more frequent is
the {\em modified Bessel equation}, which  is equivalent to the ${}_0F_1$ equation:
\begin{eqnarray*}
z^{\frac{\alpha}{2}}\F_{\alpha}(z,\partial_z)
z^{-\frac{\alpha}{2}}&=&
\partial_w^2+\frac{1}{w}\partial_w-1-\frac{\alpha^2}{w^2},\end{eqnarray*}
where $z=\frac{w^2}{4}$, $w=\pm 2\sqrt z$.

Even more frequent is the {\em Bessel equation}:
\begin{eqnarray*}
-z^{\frac{\alpha}{2}}\F_{\alpha}(z,\partial_z)
z^{-\frac{\alpha}{2}}&=&
\partial_u^2+\frac{1}{u}\partial_u+1-\frac{\alpha^2}{u^2},
\end{eqnarray*}
where $z=-\frac{u^2}{4}$, $u=\pm 2\i\sqrt z$.
Clearly, we can pass from the modified Bessel to the Bessel equation by
  $w=\pm\i u$.
\eer

\subsection{The ${}_0F_1$ function}

The ${}_0F_1$ equation
 has a regular singular point at  $0$.
Its indices at $0$ are equal to $0$, $1-c$.

If $c\neq0,-1,-2,\dots$, then the only solution of the ${}_0F_1$ equation
$\sim1$ at 0 
is called
the {\em ${}_0F_1$ hypergeometric function}.
It is 
\[F(c;z):=\sum_{j=0}^\infty
\frac{1}{
(c)_j}\frac{z^j}{j!}.\]
It is defined for $c\neq0,-1,-2,\dots$.
Sometimes it is more convenient to consider
the function
\[ {\bf F}  (c;z):=\frac{F(c;z)}{\Gamma(c)}=
\sum_{j=0}^\infty
\frac{1}{
\Gamma(c+j)}\frac{z^j}{j!}\]
defined for all $c$.

We can express the ${}_0F_1$ function in terms of the confluent function
\begin{eqnarray*}F(c;z)&=&
\e^{-2\sqrt{z}}F\Big(\frac{2c-1}{2};
2c-1;4\sqrt{z}\Big)\\&=&
\e^{2\sqrt{z}}F\Big(\frac{2c-1}{2};
2c-1;-4\sqrt{z}\Big).
\end{eqnarray*}
It is also  a limit of the confluent function.
\[F(c;z)=\lim_{a\to\infty} F(a;c;z/a).\]

For all parameters we have an integral representation called the {\em
  Schl\"afli formula}:
\begin{eqnarray*}
\frac{1}{2\pi \i}\int\limits_{]-\infty,0^+,-\infty[}
\e^t\e^{\frac{z}{t}}t^{-c}\d t
&=& {\bf F}  (c,z),\ \ \ \ \Re z>0.\end{eqnarray*}
For $\Re c>\frac12$ we have a representation called the {\em Poisson formula}:
\begin{eqnarray*}
\int_{-1}^1(1-t^2)^{c-\frac32}\e^{
2t\sqrt z}&=&\Gamma(c-\frac12)\sqrt\pi {\bf F}  (c,z).\end{eqnarray*}



We will usually prefer to use the Lie-algebraic parameters:
\begin{eqnarray*}
F_\alpha (z)&:=&F(\alpha+1;z),\\
 {\bf F}  _\alpha (z)&:=& {\bf F} (\alpha+1;z)
.
\end{eqnarray*}

\ber In the literature the ${}_0F_1$ function is seldom used. Instead,
one uses the {\em modified Bessel function} and, even more frequently,
the {\em Bessel function}:
\begin{eqnarray*}
I_\alpha(w)&=&\Big(\frac{w}{2}\Big)^\alpha {\bf F}  _\alpha\Big(\frac{w^2}{4}\Big),\\[3mm]
J_\alpha(w)&=&\Big(\frac{w}{2}\Big)^\alpha {\bf F}  _\alpha\Big(-\frac{w^2}{4}\Big).
\end{eqnarray*}
They solve the modified Bessel, resp. the Bessel equation.
\eer

\subsection{Standard solutions}

$z=0$  is a regular singular point. We have two standard solutions corresponding to its two  indices. 
Besides, we have an additional solution with a special behavior at $\infty$.

We know from Thm \ref{schl} that for appropriate contours $\gamma$ the integrals
\begin{eqnarray*}
\int\limits_{\gamma}
\e^t\e^{\frac{z}{t}}t^{-\alpha -1}\d t
\end{eqnarray*}
solve the ${}_0F_{1}$ equation.
The integrand goes to zero as $t\to-\infty$ and $t\to0-0$ (the latter for $\Re z>0$). Therefore, contours ending at these points yield solutions. We will see that in this way we can obtain all 3 standard solutions.

Besides, we can use Thm \ref{dad7} to obtain other integral representations, which are essentially special cases of  representations for the ${}_1F_1$ and ${}_2F_0$ functions.

\subsubsection{Solution $\sim1$ at 0}

If $\alpha \neq-1,-2,\dots$, then the only solution of the ${}_0F_1$ equation
$\sim1$ at 0 
is
\begin{eqnarray*}F_\alpha (z)&=&
\e^{-2\sqrt{z}}F_{0,2\alpha }\big(4\sqrt{z}\big)\\&=&
\e^{2\sqrt{z}}F_{0,2\alpha }\big(-4\sqrt{z}\big).
\end{eqnarray*}



For all parameters we have an integral representation
\begin{eqnarray*}
\frac{1}{2\pi \i}\int\limits_{]-\infty,0^+,-\infty[}
\e^t\e^{\frac{z}{t}}t^{-\alpha -1}\d t
&=& {\bf F}  _\alpha (z),\ \ \ \ \Re z>0;\end{eqnarray*}
and for $\Re \alpha >-\frac12$ we have another integral representation
\begin{eqnarray*}
\int_{-1}^1(1-t^2)^{\alpha -\frac12}\e^{
2t\sqrt z}\d t&=&\Gamma(\alpha +\frac12)\sqrt\pi  {\bf F}  _\alpha (z),\ \ 
\ \ z\not\in]-\infty,0]
.\end{eqnarray*}

\subsubsection{Solution  $\sim z^{-\alpha }$ at 0}

If $\alpha \neq 1,2,\dots$, then the only solution to the ${}_0F_1$ equation
 $\sim z^{-\alpha}$ at 0 is 
\begin{eqnarray*}z^{-\alpha }F_{-\alpha} (z)&=&
z^{-\alpha }\e^{-2\sqrt{z}}F_{0,-2\alpha }\big(4\sqrt{z}\big)\\&=&
z^{-\alpha }\e^{2\sqrt{z}}F_{0,-2\alpha }\big(-4\sqrt{z}\big).
\end{eqnarray*}


For all parameters we have 
\begin{eqnarray*}
\frac{1}{2\pi \i}\int\limits_{[(0-0)^+]}\e^t\e^{\frac{z}{t}}t^{-\alpha -1}\d t
&=&z^{-\alpha } {\bf F}  _{-\alpha }(z),\ \ \Re z>0;\end{eqnarray*}
and for $\frac12>\alpha$ we have
\begin{eqnarray*}
\int_{-\sqrt z}^{\sqrt z}(z-t^2)^{-\alpha -\frac12}\e^{
2t}\d t&=&\Gamma\Big({-}\alpha +\frac12\Big)\sqrt\pi z^{-\alpha }
    {\bf F}  _{-\alpha }(z),\ \ \ z\not\in]-\infty,0].\end{eqnarray*}

\subsubsection{Solution 
$\sim{\rm exp}(- 2z^\12) z^{-\frac{ \alpha }2-\frac14}$
for $z\to+\infty$} 
\label{s5.8}

The following function is
 also a solution of the ${}_0F_1$ equation:
\begin{eqnarray*}
\tilde F_\alpha (z)&:=&\e^{-2\sqrt z} z^{-\frac{\alpha}{2} -\frac14}
\tilde F_{0,2\alpha }\Big(-\frac{1}{4\sqrt z}\Big).
\end{eqnarray*}

We have the identity
\[\tilde F_\alpha (z)=z^{-\alpha }\tilde F_{-\alpha }(z).\]

Integral representations for all parameters:
\begin{eqnarray*}
\int_{-\infty}^{0}\e^t\e^{\frac{z}{t}}(-t)^{-\alpha -1}\d t
&=&\pi^{\12}\tilde F_\alpha (z),\ \ \ \ \Re z>0;\end{eqnarray*}
for $ \Re \alpha >-\frac12
$:
\begin{eqnarray*}\int_{-\infty}^{-1}(t^2-1)^{\alpha -\frac12}\e^{
2t\sqrt z}\d t&=&
\frac12\Gamma\Big(\alpha +\frac12\Big) \tilde F_\alpha (z),\ \ \ z\not\in]-\infty,0]; \end{eqnarray*}
for $ \Re \alpha <\frac12$:
\begin{eqnarray*}\
\int_{-\infty}^{-\sqrt z}(t^2-z)^{-\alpha -\frac12}\e^{
2t}\d t&=&\frac12\Gamma\Big(-\alpha +\frac12\Big) \tilde F_\alpha (z)
,\ \ \ z\not\in]-\infty,0].
\end{eqnarray*}

As $|z|\to\infty$ and 
$|\arg z|<\pi/2-\epsilon$, we have 
 \beq
\tilde F_\alpha (z)\sim{\rm exp}(- 2z^\12) z^{-\frac{\alpha }2-\frac14}.
\label{saddle1}\eeq
$F_\alpha$  is a unique solution with this property.

To prove (\ref{saddle1}) we can use the saddle point method. We write the left
hand side as
\[\int_0^\infty\e^{\phi(s)}s^{-\alpha-1}\d s\]
with $\phi(s)=-s-\frac{z}{s}$. We compute:
\[\phi'(s)=-1+\frac{z}{s^2},\ \ \phi''(s)=-2\frac{z}{s^3}.\]
We find the stationary point at $s_0=\sqrt z$ with
$\phi''(s_0)=-2z^{-\frac12}$ and
$\phi(s_0)=-2\sqrt z$. Hence the left hand side of (\ref{saddle1}) can be
approximated by
\[\int_{-\infty}^\infty
\e^{\phi(s_0)+\frac{1}{2}(s-s_0)^2\phi''(s_0)}s_0^{-\alpha-1}\d s=
\pi^{\frac12} z^{-\frac{\alpha}{2}-\frac14}\e^{-2\sqrt z}.\]

\ber In the literature, instead of the $\tilde F$ function one uses
the {\em MacDonald function}, solving the modified Bessel equation:
\begin{eqnarray*}
K_\alpha(w)&=&\sqrt\pi\Big(\frac{
  w}{2}\Big)^\alpha \tilde
F_\alpha\Big(\frac{w^2}{4}\Big),\end{eqnarray*}
 and the Hankel functions of the 1st and 2nd kind, solving the Bessel equation:
\begin{eqnarray*}
H_\alpha^{(1)}(w)&=&\frac{\i}{\sqrt\pi}\Big(\frac{\e^{-\i\pi/2}
  w}{2}\Big)^\alpha \tilde F_\alpha\Big(\e^{-\i\pi}\frac{w^2}{4}\Big),\\ 
H_\alpha^{(2)}(w)&=&-\frac{\i}{\sqrt\pi}\Big(\frac{\e^{\i\pi/2}
  w}{2}\Big)^\alpha\tF_\alpha\Big(\e^{\i\pi}\frac{w^2}{4}\Big).\end{eqnarray*}
\eer

\subsection{Connection formulas}

We can use the solutions with a simple behavior at zero as the basis:
\begin{eqnarray*}
\tilde F_\alpha (z)
&=&\frac{\sqrt\pi}{\sin\pi (-\alpha )} {\bf F}  _\alpha (z)
+\frac{\sqrt \pi}{\sin\pi \alpha }
z^{-\alpha } {\bf F}  _{-\alpha }(z).
\end{eqnarray*}
Alternatively, we can use the $\tilde F$ function and its analytic
continuation around $0$ in the clockwise or anti-clockwise direction
as the basis:
\begin{eqnarray*}
{\bf F}_\alpha(z)&=&\frac{1}{ 2\pi^{\frac{3}{2}}}\left(\e^{-\i\pi
 ( \alpha-\12)}
\tF_\alpha(z)-
\e^{\i\pi (\alpha-\12)}\tF_\alpha(\e^{-\i 2\pi}z)\right)\\
&=&\frac{1}{ 2\pi^{\frac{3}{2}}}\left(\e^{-\i\pi (\alpha-\12)}\tF_\alpha(\e^{\i 2\pi}z)-
\e^{\i\pi (\alpha-\12)}\tF_\alpha(z)\right),\\[3ex]
z^{-\alpha}{\bf F}_{-\alpha}(z)
&=&\frac{1}{ 2\pi^{\frac{3}{2}}}\left(\e^{\i\pi (\alpha+\12)}\tF_\alpha(z)-
\e^{-\i\pi( \alpha+\12)}\tF_\alpha(\e^{-\i 2\pi}z)\right)
\\
&=&\frac{1}{2\pi^{\frac{3}{2}}}\left(\e^{\i\pi( \alpha+\12)}
\tF_\alpha(\e^{\i 2\pi}z)-
\e^{-\i\pi (\alpha+\12)}\tF_\alpha(z)\right).
\end{eqnarray*}

\subsection{Recurrence relations}
The following recurrence relations easily follow from the commutation
relations of Subsect. \ref{symcom2a}:
 \begin{eqnarray*}
 \p_z {\bf F}  _\alpha (z)&=& {\bf F}  _{\alpha +1}(z),
 \\[3mm]
 \left(z\p_z+\alpha\right) {\bf F}  _\alpha (z)&=& {\bf F}  _{\alpha -1}(z).
\end{eqnarray*}

\subsection{Degenerate case}

$\alpha=m\in\zz$ is the degenerate case of the ${}_0F_1$ equation at $0$.
We have then
\[{\bf F}
(1+m;z)=\sum_{n=\max(0,-m)}\frac{1}{n!(m+n)!}z^n.\]
This easily implies the identity
\[{\bf F}(1+m;z)=z^{-m}
{\bf F}(1-m;z).
\]Thus the two standard solutions determined by the behavior at zero are proportional to one another. 

We have an integral representation, called the {\em Bessel formula},  and a generating function:
\begin{eqnarray*}
\frac{1}{2\pi i}\int\limits_{[0^+]}
\e^{t+z\slash t}t^{-m-1}\d t&= &{\bf F}  _m(z)=z^{-m} {\bf F} _{-m}(z),\\
\e^{t}\e^{z\slash t}&
=&\sum_{m\in\zz}t^m {\bf F}  _m(z).
\end{eqnarray*}

\subsection{Special cases}

If $\alpha=\pm\12$, then the ${}_0F_1$ equation can be reduced to an equation easily solvable in terms of elementary functions:
\begin{eqnarray*}
\F_{-\frac12}(z,\partial_z)&=&\partial_u^2-1,\\
\F_{\frac12}(z,\partial_z)&=&u^{-1}(\partial_u^2-1)u
,\end{eqnarray*}
where $u=2\sqrt z$.
They have solutions
\begin{eqnarray*}
F_{-\12}(z)=\cosh2\sqrt z,&&\tilde F_{-\12}(z)=\exp(-2\sqrt z),\\
F_{\12}(z)=\frac{\sinh2\sqrt z}{2\sqrt z},&&\tilde F_{\12}(z)=\frac{\exp(-2\sqrt z)}{\sqrt z}.
\end{eqnarray*}

\section{The Gegenbauer equation}
\label{s7}
\subsection{Introduction}
The   hypergeometric equation can be moved by an affine transformation
so that its finite singular points are placed at $-1$ and $1$. If in
addition the equation is reflection invariant,
 then it will be called the {\em Gegenbauer equation}.

Because of the reflection invariance, the third classical parameter can be obtained from the first two: $c=\frac{a+b+1}{2}$. Therefore, we will use only  ${a},{b}\in\cc$
 as the (classical) parameters of the Gegenbauer equation. It will
 be given by the operator
\beq
\S({a},{b};z,\p_z):=(1-z^2)\p_z^2-({a}+{b}+1)z\p_z-{a}{b}.
\label{geg}\eeq

To describe the symmetries of the Gegenbauer operator it is convenient to use
its Lie-algebraic  parameters
\[\begin{array}{rl}
\alpha :=\frac{{a}+{b}-1}{2},& \lambda: =\frac{{b}-{a}}{2},\\[3ex]
{a}=\frac12+\alpha -\lambda ,&{b}=\frac12+\alpha +\lambda .
\end{array}\]
Thus (\ref{geg}) becomes
\begin{eqnarray*}
\S_{\alpha ,\lambda }(z,\p_z)
&:=&(1-z^2)\p_z^2-2(1+\alpha )z\p_z
+\lambda ^2-\Big(\alpha +\12\Big)^2.\end{eqnarray*}
The Lie-algebraic parameters  have an interesting interpretation in terms of the natural basis of  the Cartan algebra of the Lie algebra $so(5)$ \cite{DM}.


\subsection{Equivalence with the hypergeometric equation}

The Gegenbauer equation is equivalent to certain subclasses of the hypergeometric
equation by a number of different substitutions. 

First of all, we can  reduce the Gegenbauer equation to the hypergeometric equation by
two affine transformations. They move the singular points from $-1$, $1$ to
$0$, $1$ or $1$, $0$:
\beq\begin{array}{l}
\S({a},{b};z,\p_z)=\F({a},{b};\frac{{a}+{b}+1}{2}
;u,\p_u),\end{array}\label{ha2}\eeq
where
\[\begin{array}{rl}&u=\frac{1-z}{2},\ \ \ z=1-2u,\\[3mm]
\hbox{or}&u=\frac{1+z}{2},\ \ \ z=-1+2u.
\end{array}\]
In the Lie-algebraic parameters
\[\S_{\alpha,\lambda}(z,\p_z)=\F_{\alpha,\alpha,2\lambda}(u,\p_u).\]

Another pair of substitutions  is a
consequence of the reflection invariance of the Gegenbauer equation (see
Subsect. \ref{a.5}): 
\beq\begin{array}{rl}\S({a},{b};z,\p_z)&=
4\F(\frac{{a}}{2},\frac{{b}}{2};\12;w,\p_w),\\[2ex]
z^{-1}\S({a},{b};z,\p_z)z&=
4\F(\frac{{a}+1}{2},\frac{{b}+1}{2};\frac32;w,\p_w),
\end{array}\label{ha1}\eeq
where \[w=z^2,\ \ \  \ 
z=\sqrt w.\]
In the Lie-algebraic parameters
\begin{eqnarray}
\label{1/2}\S_{\alpha,\lambda}(z,\p_z)&=&\F_{-\12,\alpha,\lambda}(w,\p_w),\\
\label{-1/2}z^{-1}\S_{\alpha,\lambda}(z,\p_z)z&=&\F_{\12,\alpha,\lambda}(w,\p_w).
\end{eqnarray}

\subsection{Symmetries}
\label{symcom}

All the operators below equal $\S_{\alpha ,\lambda }(w,\p_w)$ for an
appropriate $w$:
\[\begin{array}{rrcl}
w=\pm z:&&&\\
&&\S_{\alpha ,\pm\lambda },&\\[2ex]
w=\pm z:&&&\\
&(z^2-1)^{-\alpha }&\S_{-\alpha ,\pm\lambda }&
(z^2-1)^{\alpha },\\[2ex]
w=\frac{\pm z}{(z^2-1)^{\12}}:&&&\\
& (z^2-1)^{\12(\alpha +\lambda +\frac52)}
&\S_{\lambda ,\pm\alpha }
& (z^2-1)^{\12(-\alpha -\lambda -\12)},\\[2ex]
w=\frac{\pm z}{(z^2-1)^{\12}}:&&&\\
& (z^2-1)^{\12(\alpha -\lambda +\frac52)}&
\S_{-\lambda ,\pm\alpha }&
 (z^2-1)^{\12(-\alpha +\lambda -\12)}.
\end{array}\]

The symmetries of the Gegenbauer operator
have an interpretation in terms of  the Weyl group of the Lie
algebra $so(5)$.

Note that the first two symmetries from the above table are
inherited from the hypergeometric equation through the substitution 
(\ref{ha2}). 

The  symmetries involving $w=\frac{\pm z}{(z^2-1)^{\12}}$ go under the name of the {\em Whipple transformation}. To obtain them we first use
the substitution (\ref{ha1}) $z\to z^2$, then
$z^2\to\frac{z^2}{1-z^2}$, which is one of the symmetries from the Kummer's table, finally 
the substitution (\ref{ha1}) in the opposite direction
$\frac{z^2}{1-z^2}\to\sqrt{\frac{z^2}{1-z^2}}$.
We will continue our discussion of  the Whipple transformation in  Subsect. 
\ref{The Riemann surface of the Gegenbauer equation}.

\subsection{Factorizations and   commutation relations}
\label{symcoma}
There are several ways of factorizing the Gegenbauer operator:
\begin{eqnarray*}\S_{\alpha,\lambda}&=&\Big((1-z^2)\partial_z
-2(1+\alpha)z\Big)\partial_z\\
&&+\Big(\alpha+\lambda+\frac12\Big)\Big(-\alpha+\lambda-\frac12\Big)\\
&=&\partial_z\Big((1-z^2)\partial_z
-2\alpha z\Big)\\
&&+\Big(\alpha+\lambda-\frac12\Big)\Big(-\alpha+\lambda+\frac12\Big),\\
(1-z^2)\S_{\alpha,\lambda}&=&\Big((1-z^2)\partial_z
+\big(\alpha-\lambda+\frac32\big)z\Big)\Big((1-z^2)\partial_z
+\big(\alpha+\lambda+\frac12\big)z\Big)\\
&&-\Big(\alpha+\lambda+\frac12\Big)\Big(\alpha-\lambda+\frac32\Big)\\
&=&\Big((1-z^2)\partial_z
+\big(\alpha+\lambda+\frac32\big)z\Big)\Big((1-z^2)\partial_z
+\big(\alpha-\lambda+\frac12\big)z\Big)\\
&&-\Big(\alpha-\lambda+\frac12\Big)\Big(\alpha+\lambda+\frac32\Big);
\end{eqnarray*}
\begin{eqnarray*}z^2\S_{\alpha,\lambda}&=&\Big(z(1-z^2)\partial_z
-\alpha-\lambda-\frac32+\big(-\alpha+\lambda-\frac12\big)z^2\Big)
\Big(z\partial_z+\alpha+\lambda+\frac12\Big)\\
&&+\Big(\alpha+\lambda+\frac12\Big)\Big(\alpha+\lambda+\frac32\Big)\\
&=&
\Big(z\partial_z+\alpha+\lambda-\frac32\Big)
\Big(z(1-z^2)\partial_z
-\alpha-\lambda+\frac12+\big(-\alpha+\lambda-\frac12\big)z^2\Big)\\
&&+\Big(\alpha+\lambda-\frac12\Big)\Big(\alpha+\lambda-\frac32\Big)\\
&=&\Big(z(1-z^2)\partial_z
-\alpha+\lambda-\frac32+\big(-\alpha-\lambda-\frac12\big)z^2\Big)
\Big(z\partial_z+\alpha-\lambda+\frac12\Big)\\
&&+\Big(\alpha-\lambda+\frac12\Big)\Big(\alpha-\lambda+\frac32\Big)\\
&=&
\Big(z\partial_z+\alpha-\lambda-\frac32\Big)
\Big(z(1-z^2)\partial_z
-\alpha+\lambda+\frac12+\big(-\alpha-\lambda-\frac12\big)z^2\Big)\\
&&+\Big(\alpha-\lambda-\frac12\Big)\Big(\alpha-\lambda-\frac32\Big).
\end{eqnarray*}

The following commutation relations can be derived from the  factorizations:
\[\begin{array}{rrl}
&\p_z&
\S_{\alpha ,\lambda } \\[0.4ex]
&=\ \ \ \S_{\alpha +1,\lambda }&\p_z,\\[2ex]
&((1-z^2)\p_z-2\alpha z)&
\S_{\alpha ,\lambda }\\[0.4ex]
&=\ \ \ \S_{\alpha -1,\lambda }&((1-z^2)\p_z-2\alpha z),\\[2ex]
&((1-z^2)\p_z-
(\alpha +\lambda +\frac12 )z)&(1-z^2)\S_{\alpha ,\lambda }\\[0.4ex]
&=\ \ \ (1-z^2)\S_{\alpha ,\lambda +1}&
((1-z^2)\p_z-
(\alpha +\lambda +\12)z),\\[2ex]
&((1-z^2)\p_z-(\alpha -\lambda +\frac12)z)&
(1-z^2)\S_{\alpha ,\lambda }\\[0.4ex]
&=\ \ \ (1-z^2)\S_{\alpha ,\lambda -1}&((1-z^2)\p_z-(\alpha -\lambda +\12)z);
\\[3ex]
&(z\p_z+
\alpha -\lambda +\frac12)&z^2\S_{\alpha ,\lambda }\\[0.4ex]
&=\ \ \ z^2\S_{\alpha +1,\lambda -1}&(z\p_z+
\alpha -\lambda +\12),\\[2ex]
&(z(1{-}z^2)\p_z{-}\alpha{+}\lambda{+}\frac12{-}(\alpha{+}\lambda{+}\frac12)z^2)
&z^2\S_{\alpha ,\lambda }\\[0.4ex]
&=\ \ \ z^2\S_{\alpha -1,\lambda +1}
&(z(1{-}z^2)\p_z
{-}\alpha{+}\lambda{+}\frac12{-}(\alpha{+}\lambda{+}\frac12)z^2),
\\[2ex]
&(z\p_z+
\alpha +\lambda +\frac12)&z^2\S_{\alpha ,\lambda }\\[0.4ex]
&=\ \ \ z^2\S_{\alpha +1,\lambda +1}&(z\p_z+
\alpha +\lambda +\12),\\[2ex]
&(z(1{-}z^2)\p_z
{-}\alpha{-}\lambda{+}\frac12{-}(\alpha{-}\lambda{+}\frac12)z^2
)&
z^2\S_{\alpha ,\lambda }\\[0.4ex]
&=\ \ \ z^2\S_{\alpha -1,\lambda -1}&(z(1{-}z^2)\p_z
{-}\alpha{-}\lambda{+}\frac12{-}(\alpha{-}\lambda{+}\frac12)z^2).
\end{array}\]
Each of these commutation relations is associated with a root
 of the Lie algebra $so(5)$. 

Note that only the first pair of commutation relations is directly inherited from the basic commutation relations of the hypergeometric equation of Subsect. 
\ref{commu}. The next pair comes from what we called additional
commutation relations (see  Subsect. \ref{addi}), which in the reflection invariant case simplify, so that they can be counted as basic commutation relations (see a discussion in Subsect. \ref{Properties of hypergeometric type operators}).
Note that the Whipple transformation transforms the first pair of the
commutation relations into the second, and the other way around. 

The last four commutation relations  form a separate class -- they can be obtained by applying consecutively an appropriate pair from the first four commutation relations.

\subsection{The Riemann surface of the Gegenbauer equation}
\label{The Riemann surface of the Gegenbauer equation}
Let us analyze more closely the Whipple symmetry.

First let us precise the meaning of the holomorphic function involved in this symmetry. If $z\in\Omega_+:=\cc\backslash[-1,1]$, then 
$1-z^{-2}\in\cc\backslash]-\infty,0]$. Therefore,
\beq
\frac{z}{(z^2-1)^{\12}}:=\frac{1}{(1-z^{-2})^{\12}}\label{anala}\eeq
defines a unique anlytic function on
$z\in\Omega_+$ (where on the right we have the principal branch of the square root). Note that, for $z\to\infty$, (\ref{anala}) converges to $1$.

Consider a second copy of $\Omega_+$, denoted $\Omega_-$. Glue them together along $]-1,1[$, so that crossing $]-1,1[$ we go from $\Omega_\pm$ to $\Omega_\mp$. The resulting complex manifold will be called $\Omega$. 
The elements of $\Omega_\pm$ corresponding to $z\in
\cc\backslash]-1,1[$ will be denoted
 will be denoted $z_\pm$.
$\Omega$
 is biholomorphic to the sphere with 4
    punctures, which correspond to the points $-1,1,\infty_+,\infty_-$.

It is easy to see that $\Omega$ is the Riemann surface of the maximal holomorphic function extending (\ref{anala}).
On $\Omega_-$ it equals
$-\frac{z}{(z^{2}-1)^{\12}}$. 

It is useful to reinterpret this holomorphic function  as a biholomorphic function from $\Omega$ into itself:
\begin{eqnarray*}
\tau(z_+)&:=&
\left\{\begin{array}{ll}
\left(\frac{z}{\sqrt{z^2-1}}\right)_+, &{} {\rm Re} z>0,\\
\left(\frac{z}{\sqrt{z^2-1}}\right)_-,& \Re z<0,
\end{array}\right.\\
\tau(z_-)&:=&
\left\{\begin{array}{ll}
\left(-\frac{z}{\sqrt{z^2-1}}\right)_-, &{} {\rm Re} z>0,\\
\left(-\frac{z}{\sqrt{z^2-1}}\right)_+,& \Re z<0.
\end{array}\right.\end{eqnarray*}

We also introduce
\begin{eqnarray*}
\epsilon(z_\pm)&:=&z_\mp,\\
(-1)z_\pm&:=&(-z)_\pm
.\end{eqnarray*}
Note that $\tau^2=\id$, $\epsilon^2=\id$, $(-1)^2=\id$, $\tau\epsilon=(-1)\epsilon\tau$.
$\tau$ and $\epsilon$ generate
a group isomorphic to the group of
the symmetries of the square. The vertices of this square can be identified with
$(1,\infty_+,-1,\infty_-)$.
They are permuted by these transformations as follows:
\[\begin{array}{r}
 \epsilon(1,\infty_+,-1,\infty_-)=(1,\infty_-,-1,\infty_+),\\[3mm]
 (-1)(1,\infty_+,-1,\infty_-)=(-1,\infty_+,1,\infty_-),\\[3mm]
\tau(1,\infty_+,-1,\infty_-)=(\infty_+,1,\infty_-,-1).\end{array}\]

It is useful to view the Gegenbauer equation as defined on $\Omega$.

\subsection{Integral representations}

\bet\ben\item
Let $[0,1]\ni t\mapsto\gamma(t)$ satisfy
\[(t^2-1)^{\frac{{b}-{a}+1}{2}}(t-z)^{-{b}-1}\Big|_{\gamma(0)}^{\gamma(1)}=0.\]
Then
\beq
\S({a},{b};z,\p_z)
\int_\gamma (t^2-1)^{\frac{{b}-{a}-1}{2}}(t-z)^{-{b}}\d t=0
.\label{dad9}\eeq
\item
Let  $[0,1]\ni t\mapsto\gamma(t)$ satisfy
\begin{eqnarray*}&&
(t^2+2tz+1)^{\frac{-{b}-{a}}{2}+1}t^{b-2}
\Big|_{\gamma(0)}^{\gamma(1)}=0.\end{eqnarray*}
Then
\beq
\S({a},{b};z,\p_z)
\int_\gamma(t^2+2tz+1)^{\frac{-{b}-{a}}{2}}t^{{b}-1}\d t=0
.\label{dad8}\eeq
\een
\label{gqw}\eet

\proof
We compute that (\ref{dad9}) and (\ref{dad8}) equal
\begin{eqnarray*}
&&{a}
\int_\gamma \Big(\p_t(t^2-1)^{\frac{{b}-{a}+1}{2}}(t-z)^{-{b}-1}\Big)\d t,\\
&&\int\limits_\gamma\Big(
\p_t(t^2+2tz+1)^{\frac{-{b}-{a}}{2}+1}t^{b-2}\Big)\d t\end{eqnarray*}
respectively.

Note that
(\ref{dad9}) is essentially a special case of Theorem \ref{intr}.

(\ref{dad8}) can be derived from (\ref{dad9}). In fact, using
 the Whipple symmetry
 we see that, for an appropriate contour $\tilde \gamma$,
\beq\begin{array}{l}
(z^2-1)^{-\frac{{a}}{2}}\int\limits_{\tilde\gamma}(s^2-1)^{\frac{-{b}-{a}}{2}}
(s-\frac{z}{\sqrt{z^2-1}})^{{b}-1}\d s\end{array}\label{sol}\eeq
solves the Gegenbauer equation.
Then we change the variables
\[\begin{array}{l}
t=s\sqrt{z^2-1}-z,\ \ \ \ s=\frac{t+z}{\sqrt{z^2-1}},
\end{array}\]
and we obtain that (\ref{sol}) equals
\[\int_\gamma(t^2+2tz+1)^{\frac{-{b}-{a}}{2}}t^{{b}-1}\d t,\]
with an appropriate contour $\gamma$. \qed


Note that in the above theorem we can interchange $a$ and $b$. Thus we obtain four kinds of integral representations.

\subsection{Canonical forms}

The natural weight of the Gegenbauer operator equals $(z^2-1)^{\alpha}$,
so that
\[\S_{\alpha,\lambda}=
-(z^2-1)^{-\alpha}\partial_z(z^2-1)^{\alpha+1}\partial_z
+\lambda^2-\Big(\alpha+\frac{1}{2}\Big)^2.\]
The balanced form  of the Gegenbauer operator is
\begin{eqnarray*}
&&
(z^2-1)^{\frac{\alpha}{2}}\S_{\alpha,\lambda}
(z^2-1)^{-\frac{\alpha}{2}}
\\
&=&
\partial_z(1-z^2)\partial_z-\frac{\alpha^2}{1-z^2}+\lambda^2-\frac{1}{4}.
\end{eqnarray*}
Note that the symmetries $\alpha\to-\alpha$ and $\lambda\to-\lambda$ are obvious in the balanced form.

\ber\label{Legen}
 In the literature  the Gegenbauer equation is used mostly in the context of Gegenbauer polynomials, that is for $-a=0,1,2,\dots$.  In the general case,
instead of the Gegenauer equation one usually considers the so-called
{\em associated Legendre equation}. It coincides with  the balanced form of the Gegenbauer equation, except that one of its parameters is shifted by $\frac12$. In the standard form it is
\[
(1-z^2)\p_z^2-2z\p_z-\frac{m^2}{1-z^2}+l(l+1),\]
so that $m$, $l$  correspond to $\alpha$, $\lambda-\frac12$ according to our convention.
\eer

\subsection{Even  solution}

Inserting a power series into equation we see that
the Gegenbauer equation possesses an even solution equal to
\begin{eqnarray*}
S_{\alpha,\lambda}^+(z)&:=&\sum_{j=0}^\infty\frac{(\frac{{a}}{2})_j(\frac{{b}}{2})_j}{(2j)!}(2z)^{2j}\\
=\ F\Big(\frac{{a}}{2},\frac{{b}}{2};\12;z^2\Big)
&=&F_{-\12,\alpha,\lambda}(z^2)
.\end{eqnarray*}

It is the unique solution of the Gegenbauer equation satisfying
\beq S_{\alpha,\lambda}^+(0)=1,\ \ \ \ \ \frac{\d}{\d z}S_{\alpha,\lambda}^+(0)=0.\label{in1}\eeq
One way to  derive the expression in terms of the hypergeometric
function
  is to use the transformation (\ref{ha1}).

We have the identities
\begin{eqnarray*}
S_{\alpha ,\lambda }^+(z)&
=&(1-z^2)^{\frac{-1-2\alpha \pm2\lambda }{4}}
S_{\mp\lambda ,\alpha }^+
\Big(\frac{\i z}{\sqrt{1-z^2}}\Big)
\\
&=&(1-z^2)^{-\alpha }
S_{-\alpha ,\lambda }^+(z),\end{eqnarray*}
beside the obvious ones
\[S_{\alpha,\lambda}^+(z)=
S_{\alpha,-\lambda}^+(z)=S_{\alpha,\lambda}^+(-z)=S_{\alpha,-\lambda}^+(-z),\]

\subsection{Odd solution}

Similarly, the Gegenbauer equation possesses an odd solution equal to
\begin{eqnarray*}
S_{\alpha,\lambda}^-(z)&:=&\sum
\limits_{j=0}^\infty\frac{(\frac{{a}+1}{2})_j(\frac{{b}+1}{2})_j}{(2j+1)!} 
(2z)^{2j+1} \\
=\ 2zF\Big(\frac{{a}+1}{2},\frac{{b}+1}{2};\frac{3}{2};z^2\Big)
&=& 2zF_{\12,\alpha,\lambda}(z^2).
\end{eqnarray*}
It is the unique solution of the Gegenbauer equation satisfying
\beq S_{\alpha,\lambda}^-(0)=0,\ \ \ \ \ \frac{\d}{\d z}S_{\alpha,\lambda}^-(0)=2.\label{in2}\eeq

We have the identities
\begin{eqnarray*}
S_{\alpha ,\lambda }^-(z)&
=&-\i(1-z^2)^{\frac{-1-2\alpha \pm2\lambda }{4}}
S_{\mp\lambda ,\alpha }^-
\Big(\frac{\i z}{\sqrt{1-z^2}}\Big)\\
&=&(1-z^2)^{-\alpha }
S_{-\alpha ,\lambda }^-(z),\end{eqnarray*}
beside the obvious ones:
\[S_{\alpha,\lambda}^+(z)=
S_{\alpha,-\lambda}^+(z)=-S_{\alpha,\lambda}^+(-z)=-S_{\alpha,-\lambda}^+(-z),\]


\subsection{Standard solutions}

As usual, by standard solutions we mean solutions with a simple behavior around singular points. The singular points of the Gegenbauer equation are at $\{1,-1,\infty\}$. The discussion of the point $-1$ can be easily reduced to that of $1$. Therefore, it is enough to discuss $2\times2$ solutions corresponding to two indices at $1$ and $\infty$.

By Thm \ref{gqw}, for appropriate $\gamma_1$, $\gamma_2$ the integrals
\begin{eqnarray}
\int\limits_{\gamma_1}(t^2-1)^{-\frac12+\lambda }(t-z)^{-\frac12-\alpha -\lambda }\d
  t, 
&&
\label{ka1a-}\\
\int\limits_{\gamma_2}(t^2+2tz+1)^{-\alpha-\frac12}(-t)^{-\frac12+\alpha +\lambda }\d t
&&
\label{ka2a-}\end{eqnarray}
are solutions. 

The natural endpoints of $\gamma_1$ are $-1,1,z,\infty$. We will see that all standard solutions can be obtained from such integrals.

The natural endpoints of $\gamma_2$ are $z+\sqrt{z^2-1},z-\sqrt{z^2-1},
0,\infty$. Similarly, all standard solutions can be obtained from the integrals over contours with these endpoints.

It is interesting to note that in some aspects the theory of the Gegenbauer equation is more complicated than that of the hypergeometric equation. One of its manifestations is a relatively big number of natural normalizations of solutions. Indeed, let us consider e.g. integral representations of the type (\ref{ka1a-}).
The natural endpoints fall into two categories: $\{1,-1\}$ and $\{0,\infty\}$.
Therefore, we have 3 kinds of contours joining two of these endpoint:
$[-1,1]$, $[0,\infty[$ and the contours joining  two distinct categories. This corresponds two three distinct natural normalizations, which we describe in the wht follows.

\subsubsection{Solution $\sim1$  at $1$}

\label{subsub1}

If $\alpha \neq-1,-2,\dots$, then the unique solution of the Gegenbauer
equation equal to 1 at $1$  is the following function:
\begin{eqnarray*}
S_{\alpha ,\lambda }(z):&=&F_{\alpha ,\alpha ,2\lambda }
\Big(\frac{1-z}{2}\Big)=F\Big({a},{b};\frac{{a}+{b}+1}{2};\frac{1-z}{2}\Big)\\&=&
F_{\alpha ,-\frac12,\lambda }(1-z^2)=
F\Big(\frac{{a}}{2},\frac{{b}}{2};\frac{{a}+{b}+1}{2};1-z^2\Big).\end{eqnarray*}
We will also introduce several alternatively normalized functions:
\begin{eqnarray*}
{\bf S}_{\alpha,\lambda}(z)&:=&\frac{1}{\Gamma(\alpha+1)}S_{\alpha,\lambda}(z)\\
&=&\frac{1}{\Gamma(\frac{a+b+1}{2})}
F\Big({a},{b};\frac{{a}+{b}+1}{2};\frac{1-z}{2}\Big)=
{\bf F}_{\alpha,\alpha,2\lambda}\Big(\frac{1-z}{2}\Big),\\[4ex]
{\bf S}_{\alpha,\lambda}^\I(z)&:=&2^{-\frac12-\alpha+\lambda}\frac{\Gamma(\frac{1+2\alpha-2\lambda}{2})\Gamma(\frac{1+2\lambda}{2})}{\Gamma(\alpha+1)}S_{\alpha,\lambda}
(z)
\\
&=&
2^{-a}\frac{\Gamma(a)\Gamma(\frac{-a+b+1}{2})}{\Gamma(\frac{a+b+1}{2})}
F\Big({a},{b};\frac{{a}+{b}+1}{2};\frac{1-z}{2}\Big)=2^{-\frac12-\alpha+\lambda}
{\bf F}_{\alpha,\alpha,2\lambda}^\I\Big(\frac{1-z}{2}\Big),\\[4ex]
{\bf S}_{\alpha,\lambda}^{\II}(z)&:=&\frac{\Gamma(\frac{1+2\alpha-2\lambda}{2})\Gamma(\frac{1+2\alpha+2\lambda}{2})}{\Gamma(2\alpha+1)}S_{\alpha,\lambda}(z)\\
&=&
\frac{\Gamma(a)\Gamma(b)}{\Gamma(a+b)}
F\Big({a},{b};\frac{{a}+{b}+1}{2};\frac{1-z}{2}\Big),\\[4ex]
{\bf S}_{\alpha,\lambda}^{\0}(z)&:=&2^{2\alpha}\frac{\Gamma(\frac{1+2\alpha}{2})^2}{\Gamma(2\alpha+1)}S_{\alpha,\lambda}(z)=
\sqrt{\pi}
\frac{\Gamma(\frac{1+2\alpha}{2})}{\Gamma(\alpha+1)}S_{\alpha,\lambda}(z)
\\
&=&
2^{a+b-1}\frac{\Gamma(\frac{a+b}{2})^2}{\Gamma(a+b)}
F\Big({a},{b};\frac{{a}+{b}+1}{2};\frac{1-z}{2}\Big).
\end{eqnarray*}

Assuming that  $z\not\in]-\infty,-1]$, we have the following
integral representations: for $\Re\alpha+\frac12>\Re\lambda>-\12$
\begin{eqnarray}
\int\limits_{-\infty}^{-1}(t^2-1)^{-\frac12+\lambda }(z-t)^{-\frac12-\alpha -\lambda }\d
  t\nonumber 
&=&
 {\bf S}^\I_{\alpha ,\lambda }(z)
,\label{ka1a}\end{eqnarray}
and for $\Re\alpha+\12>|\Re\lambda|$
\begin{eqnarray}
\int\limits_0^{\infty}(t^2+2tz+1)^{-\alpha-\frac12}t^{-\frac12+\alpha +\lambda }\d t
\nonumber&=&
{\bf S}_{\alpha ,\lambda }^{\II}(z).
\label{ka2a}\end{eqnarray}

\subsubsection{Solution $\sim 2^{-\alpha}(1-z)^{-\alpha }$ at $1$}

\label{subsub2}

If $\alpha \neq1,2,\dots$, then the unique solution of the Gegenbauer
equation behaving as
$(1-z)^{-\alpha }$
 at $1$  is the following function:
\begin{eqnarray*}
(1-z^2)^{-\alpha }S_{-\alpha ,-\lambda }(z)&=&2^{-\alpha}(1-z)^{-\alpha }F_{-\alpha ,\alpha ,-2\lambda }
\Big(\frac{1-z}{2}\Big)\\&=&
(1-z^2)^{-\alpha }F_{-\alpha ,-\frac12,-\lambda }(1-z^2).\end{eqnarray*}

Assuming that  $z\not\in]-\infty,-1]\cup[1,\infty[$, we have the following
integral representations: for $-\Re\alpha+\frac12>\Re\lambda>-\12$
\begin{eqnarray}
\int\limits_{-1}^z(1-t^2)^{-\frac12+\lambda }(z-t)^{-\frac12-\alpha -\lambda }\d
  t\nonumber 
&=&
(1-z^2)^{-\alpha } {\bf S}^\I_{-\alpha ,\lambda }(z)
,\label{ka1b}
\end{eqnarray}
and for $\12>\Re\alpha$
\begin{eqnarray}
\int\limits_{-\i\sqrt{1-z^2}-z}^{\i\sqrt{1-z^2}-z}
(t^2+2tz+1)^{-\alpha -\frac12}(-t)^{-\frac12+\alpha +\lambda }\d t
\nonumber
&=&
(1-z^2)^{-\alpha } 
{\bf S}_{\alpha ,\lambda }^{\0}(z).
\label{ka2b}\end{eqnarray}

\subsubsection{Solution $\sim  z^{-{a}}$ at $\infty$}

\label{subsub3}
If $2\lambda \neq-1,-2,\dots$, then the unique solution of the Gegenbauer
equation behaving as $z^{-a}=z^{-\frac12-\alpha +\lambda }$ at $\infty$
is the following function:
\begin{eqnarray*}
(z^2-1)^{\frac{-1-2\alpha +2\lambda }{4}}S_{-\lambda ,-\alpha }\Big(\frac{z}{\sqrt{z^2-1}}\Big)
&=&(1+z)^{-\frac12-\alpha +\lambda }F_{-2\lambda ,\alpha ,-\alpha }
\Big(\frac{2}{1+z}\Big)\\&=&
z^{-\frac12-\alpha +\lambda }F_{-\lambda ,\alpha ,\frac12}(z^{-2}).\end{eqnarray*}

Assuming that  $z\not\in]-\infty,1]$, we have the following
integral representations: for $\frac12>\Re\lambda$
\begin{eqnarray}
\int\limits_{-1}^1(t^2-1)^{-\frac12-\lambda }(z-t)^{-\frac12-\alpha +\lambda }\d
  t\nonumber 
&=&
(z^2-1)^{\frac{-1-2\alpha +2\lambda }{4}}{\bf S}_{-\lambda ,\alpha }^{\0}\Big(\frac{z}{\sqrt{z^2-1}}\Big)
,\label{ka1c}
\end{eqnarray}
and for $-\lambda+\12>-\Re\alpha>-\12$
\begin{eqnarray}
\int\limits_{\sqrt{z^2-1}-z}^{0}
(t^2+2tz+1)^{-\alpha-\frac12}(-t)^{-\frac12+\alpha -\lambda }\d t
\nonumber
&=&(z^2-1)^{\frac{-1-2\alpha +2\lambda }{4}}{\bf S}_{-\lambda ,\alpha }^\I\Big(\frac{z}{\sqrt{z^2-1}}\Big).
\label{ka2c}\end{eqnarray}

\subsubsection{Solution $\sim  z^{-{b}}$ at $\infty$}

\label{subsub4}
If $2\lambda \neq1,2,\dots$, then the unique solution of the Gegenbauer
equation behaving as $z^{-b}=z^{-\frac12-\alpha -\lambda }$ at $\infty$
is the following function:
\begin{eqnarray*}
(z^2-1)^{\frac{-1-2\alpha -2\lambda }{4}}S_{\lambda ,\alpha }\Big(\frac{z}{\sqrt{z^2-1}}\Big)
&=&(1+z)^{-\frac12-\alpha -\lambda }F_{2\lambda ,\alpha ,\alpha }
\Big(\frac{2}{1+z}\Big)\\&=&
z^{-\frac12-\alpha -\lambda }F_{\lambda ,\alpha ,\frac12}(z^{-2}).\end{eqnarray*}

Assuming that  $z\not\in]-\infty,1]$, we have the following
integral representations: for $\Re\lambda+\frac12>|\Re\alpha|$
\begin{eqnarray}
\int\limits_z^\infty(t^2-1)^{-\frac12-\lambda }(t-z)^{-\frac12-\alpha +\lambda }\d
  t\nonumber 
&=&(z^2-1)^{\frac{-1-2\alpha -2\lambda }{4}}{\bf S}_{\lambda ,\alpha }^{\II}\Big(\frac{z}{\sqrt{z^2-1}}\Big)
,\label{ka1d}
\end{eqnarray}
and for $\lambda+\12>-\Re\alpha>-\12$
\begin{eqnarray}
\int\limits_\infty^{-\sqrt{z^2-1}-z}
(t^2+2tz+1)^{-\alpha -\frac12}t^{-\frac12+\alpha -\lambda }\d t
\nonumber
&=&(z^2-1)^{-\frac14-\frac{\alpha}{2} -\frac{\lambda}{2}}{\bf S}_{\lambda ,\alpha }^\I\Big(\frac{z}{\sqrt{z^2-1}}\Big)
.
\label{ka2d}\end{eqnarray}

\ber As mentioned in Remark \ref{Legen}, in the literature instead of the Gegenbauer equation the associated Legendre equation  usually appears.
One class of its standard  solutions  are  the
{\em associated Legendre function of the 1st kind}
\begin{eqnarray*}
{\bf P}_l^m(z)&=&\left(\frac{z+1}{z-1}\right)^{\frac{m}{2}}
{\bf F} \Big(-l,l+1;1-m;\frac{1-z}{2}\Big)\\
&=&\frac{2^m}{(z^2-1)^{\frac{m}{2}}}
{\bf F} \Big(1-m+l,-m-l;1-m;\frac{1-z}{2}\Big)\\
&=&\frac{2^m}{(z^2-1)^{\frac{m}{2}}}
{\bf S}_{-m,l+\frac12}(z),\end{eqnarray*}
which up to a constant are $(z^2-1)^{\frac{m}{2}}$ times
the solutions of Subsubsect.  \ref{subsub2}.
Another class of solutions are the {\em associated Legendre function of the 2nd kind}
\begin{eqnarray*}
{\bf Q}_l^m(z)&=&
\frac{(z^2-1)^{\frac{m}{2}}}{2^{l+1}z^{l+m+1}} {\bf F} \Big(\frac{l+m+2}{2},\frac{l+m+1}{2};l+\frac{3}{2};z^{-2}\Big)\\
&=&
2^{-l-1}(z^2-1)^{-\frac{l+1}{2}}{\bf S}_{l+\frac12,m}\Big(\frac{z}{\sqrt{z^2-1}}\Big),\end{eqnarray*}
which up to a constant are  $(z^2-1)^{\frac{m}{2}}$ times
the solutions  of Subsubsect. 
 \ref{subsub4}.
(In the literature one can find a couple of other varieties of associated Legendre functions of the 1st and 2nd kind, differing by their normalization, see e.g. \cite{NIST}).
\end{remark}

\subsection{Connection formulas}

We can express the standard solutions in terms of the even and odd solutions
\begin{eqnarray*}\nonumber
 {\bf S}_{\alpha ,\lambda }(z)
&=&\frac{\sqrt\pi}
{\Gamma(\frac{3}{4}+\frac{\alpha}{2} -\frac{\lambda }{2})\Gamma(\frac{3}{4}+\frac{\alpha}{2} +\frac{\lambda }{2})}
S_{\alpha ,\lambda }^+(z)\\[1ex]&&
+\frac{\sqrt{\pi}}
{\Gamma(\frac{1}{4}+\frac{\alpha}{2} -\frac{\lambda}{2})\Gamma(\frac{1}{2}+\frac{\alpha}{2} +\frac{\lambda }{2})}
S_{\alpha ,\lambda }^-(z),\\[2ex]
(1-z^2)^{-\alpha }
 {\bf S}_{-\alpha ,-\lambda }(z)
&=&\frac{\sqrt\pi}
{\Gamma(\frac{3}{4}-\frac{\alpha}{2} +\frac{\lambda}{2})\Gamma(\frac{3}{4}-\frac{\alpha}{2} -\frac{\lambda }{2})}
S_{\alpha ,\lambda }^+(z)\\[1ex]
&&+\frac{\sqrt{\pi}}
{\Gamma(\frac{1}{4}-\frac{\alpha}{2} +\frac{\lambda }{2})\Gamma(\frac{1}{4}-\frac{\alpha}{2} -\frac{\lambda }{2})}
S_{\alpha ,\lambda }^-(z),\\[2ex]
(1-z^2)^{-\frac{1}{4}-\frac{\alpha}{2} +\frac{\lambda }{2}}{\bf S}_{-\lambda ,-\alpha } (\frac{z}{\sqrt{z^2-1}})
&=&\frac{\sqrt\pi}
{\Gamma(\frac{3}{4}-\frac{\alpha}{2} -\frac{\lambda }{2})\Gamma(\frac{3}{4}+\frac{\alpha}{2} -\frac{\lambda }{2})}
S_{\alpha ,\lambda }^+(z)\\[1ex]
&&+\frac{\i\sqrt{\pi}}
{\Gamma(\frac{1}{4}-\frac{\alpha}{2} -\frac{\lambda }{2})\Gamma(\frac{1}{4}+\frac{\alpha}{2} -\frac{\lambda }{2})}
S_{\alpha ,\lambda }^-(z),\\[2ex]
(1-z^2)^{-\frac{1}{4}-\frac{\alpha}{2} -\frac{\lambda }{2}}{\bf S}_{\lambda ,\alpha } \Big(\frac{z}{\sqrt{z^2-1}}\Big)
&=&\frac{\sqrt\pi}
{\Gamma(\frac{3}{4}-\frac{\alpha}{2} +\frac{\lambda }{2})\Gamma(\frac{3}{4}+\frac{\alpha}{2} +\frac{\lambda }{2})}
S_{\alpha ,\lambda }^+(z)\\[1ex]
&&+\frac{\i\sqrt{\pi}}
{\Gamma(\frac{1}{4}-\frac{\alpha}{2} +\frac{\lambda}{2})\Gamma(\frac{1}{4}+\frac{\alpha}{2} +\frac{\lambda }{2})}
S_{\alpha ,\lambda }^-(z).
\end{eqnarray*}

\subsection{Recurrence relations}

The following recurrence relations can be easily  derived from the
commutation properties of Subsect. \ref{symcoma}
\begin{eqnarray*}
\p_z {\bf S}_{\alpha ,\lambda }(z)&=&-\12\Big(\frac12+\alpha -\lambda \Big)\Big(\frac12+\alpha +\lambda \Big)
 {\bf S}_{\alpha +1,\lambda }(z),\\
\left((1-z^2)\p_z -2\alpha z\right)
{\bf S}_{\alpha ,\lambda }(z)
&=&-2{\bf S}_{\alpha -1,\lambda }(z) ,\\[2ex]
\left((1-z^2)\p_z -\Big(\frac12+\alpha +\lambda \Big) z\right){\bf S}_{\alpha ,\lambda }(z)
&=&-\Big(\frac12+\alpha +\lambda \Big) {\bf S}_{\alpha ,\lambda +1}(z),
\\
\left((1-z^2)\p_z -\Big(\frac12+\alpha -\lambda \Big) z\right){\bf S}_{\alpha ,\lambda }(z)
&=&-\Big(\frac12+\alpha -\lambda \Big){\bf S}_{\alpha ,\lambda -1}(z),
\\[2ex]
\left(z\p_z+\frac12+\alpha -\lambda \right){\bf S}_{\alpha ,\lambda }(z)
&=&\frac12\Big(\frac12+\alpha -\lambda \Big)\Big(\frac32+\alpha-\lambda\Big) {\bf S}_{\alpha +1,\lambda -1}(z),\\
\left(z(1{-}z^2)\p_z{+}\Big(\frac12{-}\alpha {+}\lambda \Big)(1{-}z^2){-}2\alpha z^2
\right){\bf S}_{\alpha ,\lambda }(z)
&=&-2{\bf S}_{\alpha -1,\lambda +1}(z),\\[2ex]
\left(z\p_z+\frac12+\alpha +\lambda \right){\bf S}_{\alpha ,\lambda }(z)
&=&\Big(\frac12+\alpha +\lambda \Big)(\alpha+1) {\bf S}_{\alpha +1,\lambda +1}(z),\\
\left(z(1{-}z^2)\p_z{+}\Big(\frac12{-}\alpha {-}\lambda \Big)(1{-}z^2){-}2\alpha  z^2
\right){\bf S}_{\alpha ,\lambda }(z)
&=&-2{\bf S}_{\alpha -1,\lambda -1}(z).
\end{eqnarray*}


\subsection{Gegenbauer polynomials}

If $-a=n=0,1,2,\dots$, then  Gegenbauer functions are polynomials.

 We will use two distinct normalizations of these polynomials.
The $C_n^\I$ polynomials have a natural Rodriguez-type definition:
\[
C_n^{\I,\alpha}(z):=\frac{1}{2^nn!}(z^2-1)^{-\alpha}\p_z^n
(z^2-1)^{n+\alpha}.
\]
The $C_n^\II$ polynomials are defined as 
\[C_n^{\II,\alpha}(z):=\frac{(2\alpha+1)_n}{(\alpha+1)_n}
C_n^{\I,\alpha}(z).\]
\ber
The first kind polynomials is just the special case of the conventional Jacobi polynomials (see Rem. \ref{jaco}) with $\alpha=\beta$:
\[C_n^{\I,\alpha}(z)=P_n^{\alpha,\alpha}(z).\]
The second kind of polynomials is called in the literature the {\em Gegenbauer polynomials}. In the standard notation its parameter is shifted by $\frac12$:
\[C_n^{\II,\alpha}(z)=C_n^{\alpha+\frac12}(z).\]
\eer

When describing the properties of Gegenbauer polynomials we can choose
between  $C_n^{\I}$ and  $C_n^{\II}$. We either give properties of
both kinds of polynomials or choose those that give simpler
formmulas.

Both kinds of polynomials solve the Gegenbauer equation:
\begin{eqnarray*}
\Big(
(1-z^2)\partial_z^2-2(1+\alpha)z\partial_z+n(n+2\alpha+1)\Big)C_n^{\I/\II}(z)&&\\
=\ \ {\mathcal S}(-n,n+2\alpha+1;z,\partial_z)C_n^{\I/\II}(z)&=&0.
\end{eqnarray*}
Generating functions:
\begin{eqnarray*}
(1-2tz+t^2(z^2-1))^{-\alpha}&=&\sum\limits_{n=0}^\infty
(2t)^nC_n^{\I,-\alpha-n}(z),\\
(1-2zt+t^2)^{-\alpha-\frac12}&
=&\sum\limits_{n=0}^\infty C_n^{\II,\alpha}(z)t^n.
\end{eqnarray*}
Integral representations:
\begin{eqnarray*}C_n^{\I,\alpha}(z)&=&\frac{1}{2\pi\i}\int\limits_{[0^+]}
\Big(1-2tz+t^2(z^2-1)\Big)^{\alpha+n}t^{-n-1}\d t,\\
C_n^{\II,\alpha}(z)&
=&\frac{1}{2\pi\i}\int\limits_{[0^+]}(1-2zt+t^2)^{-\alpha-\frac12}t^{-n-1}\d t
.\end{eqnarray*}
We give symmetries for both kinds of polynomials:
\begin{eqnarray*}C_n^{\I,\alpha}(z)
&=&(- 1)^nC_n^{\I,\alpha}(- z)\\
&=&\frac
{(2\alpha+1+n)_n}
{(\mp 2)^{n}(\alpha+\frac12)_n}(z^2-1)^{\frac{n}{2}}C_n^{\I,-\12-\alpha-n}
\Big(\frac{\pm z}{\sqrt{z^2-1}}\Big)
.\end{eqnarray*}
\begin{eqnarray*}C_n^{\II,\alpha}(z)
&=&(- 1)^nC_n^{\II,\alpha}(- z)\\
&=&\frac{(\mp 2)^n(\alpha+\frac12)_n}
{(2\alpha+1+n)_n}(z^2-1)^{\frac{n}{2}}C_n^{\II,-\12-\alpha-n}
\Big(\frac{\pm z}{\sqrt{z^2-1}}\Big)
.\end{eqnarray*}
We give recurrence relations only for $C_n^{\II,\alpha}$, those for $C_n^{\I,\alpha}$ differ by coefficients on the right, but have a comparable level of complexity: 
\begin{eqnarray*}
\p_z C_n^{\II,\alpha}(z)&=&(2\alpha+1) C_{n-1}^{\II,\alpha+1}(z),
\\
\left((1-z^2)\p_z-2\alpha z\right)C_n^{\II,\alpha}(z)
&=&\frac{-(n+1)(n+2\alpha)}{2\alpha}C_{n+1}^{\II,\alpha-1}(z),
\\[5ex]
\left((1-z^2)\p_z-(n+2\alpha+1)z\right)C_n^{\II,\alpha}(z)
&=&-(n+1)C_{n+1}^{\II,\alpha}(z),
\\
\left((1-z^2)\p_z+nz\right)C_n^{\II,\alpha}(z)
&=&(n+2\alpha)C_{n-1}^{\II,\alpha}(z),\\[7ex]
(z\p_z-n)C_n^{\II,\alpha}(z)&=&(2\alpha+1) C_{n-2}^{\II,\alpha+1}(z),\\
\left(z(1-z^2)\p_z+1+n-(n+2\alpha+1)z^2\right)C_n^{\II,\alpha}(z)
&=&-\frac{(n+1)(n+2)}{2\alpha-1}C_{n+2}^{\II,\alpha-1}(z),\\[5ex]
\left(z\p_z+n+2\alpha+1\right)C_n^{\II,\alpha}(z)&=&(2\alpha+1) C_n^{\II,\alpha+1}(z),\\
\left(z(1-z^2)\p_z-n-2\alpha+nz^2\right)C_n^{\II,\alpha}(z)
&=&-\frac{(2\alpha+n-1)(2\alpha+n)}{2\alpha-1}C_{n}^{\II,\alpha-1}(z).
\end{eqnarray*}

The differential equation, the Rodriguez-type formula, the first generating function and the first integral representation 
are special cases of 
 the  corresponding formulas of Subsect.  
\ref{Hypergeometric type polynomials}. Thus the polynomials
 $C^\I$ belong to the scheme of Subsect.
\ref{Hypergeometric type polynomials}.
$C^{\II}$ do not  have a natural Rodriguez-type formula, and  do not belong to the scheme of Subsect. \ref{Hypergeometric type polynomials}.

The  $C^{\I}$ polynomials have  simple expressions in terms of the Jacobi polynomials:
\begin{eqnarray*}
C_n^{\I,\alpha}(z)&
=&(\pm 1)^nR_n^{\alpha,\alpha}\Big(\frac{1\mp z}{2}\Big)\\
&=&
\Big(\frac{\pm 1-z}{2}\Big)^n
R_n^{\alpha,-2\alpha-2n-1}\Big(\frac{2}{1\mp z}\Big)\\[3mm]
&=&
\Big(\frac{z\mp 1}{2}\Big)^nR_n^{-2\alpha-2n-1,\alpha}\Big(\frac{\pm 1+z}{\mp 1+z}\Big).
\end{eqnarray*}
We have several alternative expressions for $C^\I$ and $C^\II$ polynomials:
\begin{eqnarray*}
C_n^{\I,\alpha}(z)&:=&\lim\limits_{\nu\to n}(-1)^n(\nu{-}n){\bf S}_{\alpha,\nu+\alpha+\frac12}^{\I}(z)=
\lim\limits_{\nu\to n}(\nu{-}n){\bf F}_{\alpha,\alpha,2\nu+2\alpha+1}^{\I}\Big(\frac{1\mp z}{2}\Big)\\&=&(\pm1)^n\frac{(\alpha+1)_n}{n!}F\Big(-n,n+2\alpha+1;\alpha+1;
\frac{1\mp z}{2}\Big),\\
C_n^{\II,\alpha}(z)&:=&\lim\limits_{\nu\to n}(-1)^n(\nu{-}n){\bf S}_{\alpha,\nu+\alpha+\frac12}^{\II}(z)\\
&=&(\pm1)^n\frac{(2\alpha+1)_n}{n!}F\Big(-n,n+2\alpha+1;\alpha+1;
\frac{1\mp z}{2}\Big)\\
&=&\sum\limits_{k=0}^{[\frac{n}{2}]}\frac{(-1)^k(\alpha+\frac12)_{n-k}}{k!(n-2k)!}(2z)^{n-2k}
.
\end{eqnarray*}

Values at $\pm1$, behavior at infinity we give for both kinds of polynomials:
\begin{eqnarray*}
C_n^{\I,\alpha}(\pm1)&=\ (\pm 1)^n\frac{(\alpha+1)_n}{n!},\ \ \ \
\lim\limits_{z\to\infty}\frac{C_n^{\I,\alpha}(z)}{z^n}&=\ \frac{2^{-n}(2\alpha+n+1)_n}{n!},\\
C_n^{\II,\alpha}(\pm1)&=\ (\pm 1)^n\frac{(2\alpha+1)_n}{n!},\ \ \ \
\lim\limits_{z\to\infty}\frac{C_n^{\II,\alpha}(z)}{z^n}&=\ \frac{2^n(\alpha+\frac12)_n}{n!}.
\end{eqnarray*}

The degenerate case has a simple expression in terms of $C^\I$
polynomials:
\[C_n^{\I,\alpha}=\Big(\frac{2}{z^2-1}\Big)^{\alpha}
C_{n+2\alpha}^{\I,-\alpha}(z),\ \ \alpha\in\zz.\]

The initial conditions at $0$ and the identities for the even and odd case are given only for $C_n^{\II,\alpha}$, since those for $C_n^{\I,\alpha}$ are more complicated:
\begin{eqnarray*}
C_{2m}^{\II,\alpha}(0)&=\ \frac{(-1)^m(\alpha+\frac12)_m}{m!},\ \ 
\p_zC_{2m}^{\II,\alpha}(0)&=\ 0;\\
C_{2m+1}^{\II,\alpha}(0)&=\ 0,\ \ \ \ \ \ \ \ \ \ \ \p_zC_{2m+1}^{\II,\alpha}(0)
&=\ \frac{(-1)^m2(\alpha+\frac12)_m}{m!}.\end{eqnarray*}

\begin{eqnarray*}
C_{2m}^{\II,\alpha}(z)&=&(-1)^m\frac{(\alpha+\frac12)_m}{(\alpha+1)_m} R_m^{\alpha,-\frac12}(z^2)\\
&=&(-1)^m\frac{(\alpha+\frac12)_m}{m!}S_{\alpha,2m+\frac12+\alpha}^+(z)\\
&=&(-1)^m\frac{(\alpha+\frac12)_m}{m!}F\Big(-m,m+\frac12+\alpha;\frac12;z^2\Big)
,\\
C_{2m+1}^{\II,\alpha}(z)&=&(-1)^m\frac{(\alpha+\frac12)_m}{(\alpha+1)_m}
2zR_m^{\alpha,\frac12}(z^2)\\ 
&=&(-1)^m\frac{(\alpha+\frac12)_m}{m!}S_{\alpha,2m+\frac32+\alpha}^-(z)\\
&=&(-1)^m\frac{(\alpha+\frac12)_m}{m!}2zF\Big(-m,m+\frac32+\alpha;\frac32;z^2\Big)  
.\end{eqnarray*}

We have the following special cases:
\ben
\item
If $\alpha\in\zz$, $-n\leq\alpha\leq\frac{-n-1}{2}$, then $C_n^{\I,\alpha}=0$.
\item
If $\alpha\in\zz+\frac12$, $\frac{-n-1}{2}\leq\alpha\leq-\frac{1}{2}$, then $C_n^{\II,\alpha}=0$.
\item
If $\alpha\in\zz$, $\frac{-n+1}{2}\leq\alpha\leq-1$, then $C_n^{\I/\II,\alpha}=(1-z^2)^{-\alpha} W$, where $W$ is a polynomial not divisible by $1-z^2$.
\een

\subsection{Special cases}

When describing special cases of the Gegenbauer quation we will primarily use the Lie-algebraic parameters.

\subsubsection{The Legendre equation}

Suppose that one of the parameters is an integer. Using, if necessary, recurrence relations we can assume it is zero. After applying an appropriate symmetry, we can assume that $\alpha=0$. We obtain then the {\em Legendre operator}:
\begin{eqnarray}
\S_{0,\lambda }(z,\p_z)
&=&(1-z^2)\p_z^2-2z\p_z
+\lambda ^2-\frac14.\label{legen}\end{eqnarray}
For the particular case $\lambda=0$ its solutions can be expressed by the so called {\em complete elliptic functions}.

The Legendre operator for polynomials of degree $n$ has the form
\begin{eqnarray*}
&&(1-z^2)\p_z^2-2z\p_z
+n(n+1).\end{eqnarray*}
The {\em Legendre polynomials} are special cases of both $C^\I$ and $C^\II$:
\begin{eqnarray*}
P_n(z)&=&C_n^{\I,0}(z)=C_n^{\II,0}(z)\\&=&\frac{1}{2^nn!}\p_z^n
(z^2-1)^{n}.
\end{eqnarray*}
Their generating function is a special case of the generating function for $C^\II$:
\begin{eqnarray*}
(1-2zt+t^2)^{-\frac12}&
=&\sum\limits_{n=0}^\infty P_n(z)t^n.
\end{eqnarray*}

\subsubsection{Chebyshev equation of the 1st kind}
Suppose that one of the parameters belongs to $\zz+\frac12$. Using, if necessary, recurrence relation, we can assume it equals $-\frac12$. After applying  an appropriate symmetry we can assume that $\alpha=-\frac12$. We obtain then the {\em Chebyshev operator of the 1st kind}:
\begin{eqnarray}
\S_{0,\lambda }(z,\p_z)
&=&(1-z^2)\p_z^2-2z\p_z
+\lambda ^2.\label{cheb1}\end{eqnarray}
After substitution $z=\cos\phi$ it becomes
\[\partial_\phi^2+\lambda^2.\]
Thus the coresponding 
equation can be solved in terms of elementary functions.

To obtain an operator that annihilates a polynomial of degree $n$ we simply set $\lambda=n$:
\begin{eqnarray*}
&&(1-z^2)\p_z^2-2z\p_z
+n^2.\end{eqnarray*}
The {\em Chebyshev polynomials of the 1st kind} are 
\begin{eqnarray*}
T_n(z)&=& 
\frac{n!}{(1/2)_n}C_n^{\I,-\frac12}(z)=\frac{\d}{\d\alpha}C_n^{\II,-\12}(z)\\
&=&
\frac12\Big((z+\i\sqrt{1-z^2})^n+
(z-\i\sqrt{1-z^2})^n\Big).
\end{eqnarray*}
Note that $C_n^{\II,-\frac12}=0$, therefore the usual generating function for $C^\II$ cannot be applied for the Chebyshev polynomials of the 1st kind. Instead, we have generating functions
\begin{eqnarray*}
-\log(1-2zt+t^2)&
=&\sum\limits_{n=0}^\infty T_n(z)\frac{t^n}{n},\\
\frac{1-zt}{1-2zt+t^2}&
=&\sum\limits_{n=0}^\infty T_n(z)t^n.
\end{eqnarray*}
\subsubsection{Chebyshev equation of the 2nd kind}
If one of the parameters belongs to $\zz+\frac12$, instead of $\alpha=-\frac12$ we can reduce ourselves to the case
$\alpha=\frac12$. We obtain then the {\em Chebyshev operator of the 2nd kind}:
\begin{eqnarray}\label{cheb2}
\S_{0,\lambda }(z,\p_z)
&=&(1-z^2)\p_z^2-2z\p_z
+\lambda ^2-1.\end{eqnarray}
After substitution $z=\cos\phi$ it becomes
\[\sin\phi(\partial_\phi^2+\lambda^2)(\sin\phi)^{-1}.\]
Clearly, the corresponding  equation can also be solved in  elementary functions.

To obtain an operator that annihilates a polynomial of degree $n$ we  set $\lambda=n+1$:
\begin{eqnarray*}
&&(1-z^2)\p_z^2-2z\p_z
+n(n+2).\end{eqnarray*}
The
{\em Chebyshev polynomials of the 2nd kind} are
\begin{eqnarray*}
U_n(z)&=& 
\frac{n!}{(3/2)_n}C_n^{\I,\frac12}(z)=C_n^{\II,\12}(z)\\
&=&
\frac{(z+\i\sqrt{1-z^2})^{n+1}-
(z-\i\sqrt{1-z^2})^{n+1}}{2\i\sqrt{1-z^2}}.
\end{eqnarray*}
Their generating function is a special case of the generating function for $C^\II$:
\begin{eqnarray*}
(1-2zt+t^2)^{-1}&
=&\sum\limits_{n=0}^\infty U_n(z)t^n.
\end{eqnarray*}

\section{The Hermite equation}
\label{s8}
\subsection{Introduction}
Let ${a}\in\cc$.
In this section we study the {\em Hermite equation}, which is given by
the operator
\[\S({a},z,\p_z):=\p_z^2-2z\p_z-2{a}.\]
The choice of the parameter $a$ is dictated by the analogy with the
parameters of the Gegenbauer. It will be called a {\em classical parameter},
even though it is not the usual one in the literature.

The Hermite operator can be obtained as the limit of the Gegenbauer operator:
\beq
\lim_{b\to\infty}\frac{2}{b}\S\Big(a,b;z\sqrt{2/b},\partial_{\big(z\sqrt{2/b}\big)}\Big)
=\S(a;z,\partial z).\label{psda}\eeq

To describe the symmetries it is convenient to use  its  {\em Lie-algebraic parameter}:
\[\lambda={a}-\frac12,\ \ \ a=\lambda+\frac12.\]
In the new parameter the Hermite operator equals
\begin{eqnarray*}
\S_\lambda (z,\p_z)
&=&\p_z^2-2z\p_z-2\lambda -1.
\end{eqnarray*}
The Lie-algebraic parameter  has an interesting interpretation in terms of  a ``Cartan element'' of the Lie algebra $sch(1)$ \cite{DM}.


\subsection{Equivalence with a subclass of the confluent equation}

The Hermite equation is reflection invariant. By using the quadratic
 transformation
 we can reduce it to a special case of the confluent equation:
\begin{eqnarray}
\S({a};z,\p_z)&=&4\F(\frac{{a}}{2};\12;w,\p_w),\label{ha6}\\[2ex]
z^{-1}\S({a};z,\p_z)z&=&4\F(\frac{{a}+1}{2};\frac32;w,\p_w),
\label{ha6a}\end{eqnarray}
where
\[w=z^2,\ \ \ \ z=\sqrt w.\]
In the Lie-algebraic parameters
\begin{eqnarray*}
\S_\lambda(z,\p_z)&=&4\F_{\lambda,-\frac12}(w,\p_w),\\
z^{-1}\S_\lambda(z,\p_z)z&=&4\F_{\lambda,\frac12}(w,\p_w).
\end{eqnarray*}

\subsection{Symmetries}
\label{symcom5}
 
The following operators equal $\S_\lambda (w,\p_w)$ for an
appropriate $w$: 
 \[\begin{array}{rrcl}
w=\pm z:&&&\\
&&\S_\lambda (z,\p_z),&\\[1ex]
w=\pm\i z:&&&\\
&-\exp(-z^2)&\S_{-\lambda }(z,\p_z)&\exp(z^2).
\end{array}\]
The group of symmetries of the Hermite equation is isomorphic to $\zz_4$ and can be interpreted as the ``Weyl group'' of $sch(1)$.

\subsection{Factorizations and commutation properties}
\label{symcom5a}
There are several ways to factorize the Hermite operator:
\begin{eqnarray*}
\S_\lambda&=&\big(\p_z-2z\big)\p_z-2\lambda-1\\
&=&\p_z\big(\p_z-2z\big)-2\lambda+1,\\[2ex]
z^2\S_\lambda&=&
\Big(z\p_z+\lambda-\frac32\Big)\Big(z\p_z-\lambda+\frac12-2z^2\Big)\\
&&+\Big(\lambda-\frac32\Big)\Big(\lambda-\frac12\Big)\\
&=&
\Big(z\p_z-\lambda-\frac32-2z^2\Big)
\Big(z\p_z+\lambda+\frac12\Big)\\
&&+\Big(\lambda+\frac32\Big)\Big(\lambda+\frac12\Big).
\end{eqnarray*}
The factorizations can be used to derive the following commutation relations:
\[\begin{array}{rl}
\p_z&\S_\lambda \\
=\ \ \S_{\lambda +1}&\p_z,\\[2ex]
(\p_z-2z)
&\S_\lambda \\
=\ \ \ \S_{\lambda -1}&(\p_z-2z),\\[3ex]
(z\p_z+\lambda +\frac12)&z^2\S_\lambda \\
=\ \ z^2\S_{\lambda +2}& (z\p_z+\lambda +\12),\\[2ex]
(z\p_z-\lambda +\frac12-2z^2)&z^2\S_\lambda \\
=\ \ z^2\S_{\lambda -2}&(z\p_z-\lambda +\12-2z^2).
\end{array}\]
Each of these commutations relations is associated with a ``root'' of the Lie algebra $sch(1)$.

\subsection{Convergence of the Gegenbauer equation to the Hermite equation}

It is interesting to describe the transition from the symmetries of the 
Gegenbauer equation to the symmetries of the Hermite equation. We
consider the limit (\ref{psda}).  We also consider  the surface
$\Omega$ described in Subsect.
\ref{The Riemann surface of the Gegenbauer equation}.

Let us look  only at the part of $\Omega$ given by the union of
$\Omega_+\cap\{\Im z>0\}_+$ and 
 $\Omega_-\cap\{\Im z>0\}_-$ glued along $]-1,1[$.
The scaling involved in the limit (\ref{psda}) transforms this part
of $\Omega$ into $\cc$.

$\tau(\Omega_+\cap\{\Im z>0\})$ 
 is equal to the union of $\Omega_-\cap\{\Im z>0,\Re z>0\}$ and
 $\Omega_-\cap\{\Im z<0,\Re z>0\}$ glued along $]0,1[$. Thus 
the limit of $\tau$ on
$\Omega_+\cap\{\Im z>0\}$ equals the multiplication by $-\i$.

$-\tau(\Omega_-\cap\{\Im z<0\})$ 
 is equal to the union of
 $\Omega_+\cap\{\Im z>0,\Re z<0\}$  and  $\Omega_-\cap\{\Im z<0,\Re z<0\}$ 
glued along $]-1,0[$. Thus the limit of $-\tau$ on
$\Omega_-\cap\{\Im z>0\}$ also equals the multiplication by $-\i$.

Thus the multiplication by $-\i$ is not the limit of a single element of the
group of the symmetries of ther Gegenbauer equation, but a combination of the
limits of two symmetries. 

\subsection{Integral representations}

Below we describe  two kinds of integral representations 
of the Hermite equation.

\bet\ben\item
Let $[0,1]\ni t\mapsto\gamma(t)$ satisfy
\[\e^{t^2}(t-z)^{-{a}-1}\Big|_{\gamma(0)}^{\gamma(1)}=0.\]
Then
\beq\S({a};z,\p_z)\int_\gamma\e^{t^2}(t-z)^{-{a}}\d t.\label{dad10}\eeq
\item
Let $[0,1]\ni t\mapsto\gamma(t)$ satisfy
\[\e^{-t^2-2zt}t^{{a}}\Big|_{\gamma(0)}^{\gamma(1)}=0.\]
Then\beq
\S({a};z,\p_z)
\int_\gamma
\e^{-t^2-2zt}t^{a-1}\d t=0.\label{dad11}\eeq
\een\label{dad12}\eet

\proof
We check that for any contour $\gamma$, (\ref{dad10}) and (\ref{dad11}) equal
\begin{eqnarray*}
&&-{a} \int_\gamma\left(\p_t
\e^{t^2}(t-z)^{-{a}-1}\right)\d t,\\
&&-2\int_\gamma
\Big(\p_t\e^{-t^2-2zt}t^{{a}}\Big)\d t
\end{eqnarray*}
respectively.

We can also deduce the second  representation from the first
by the symmetry involving the multiplication by $\e^{z^2}$ and the change of variables
$z\mapsto\i z$.
\qed

\subsection{Canonical forms}

The natural weight of the Hermite operator equals $\e^{-z^2}$, so that
\[\S_{\lambda}=
\e^{z^2}\partial_z\e^{-z^2}\partial_z-2\lambda-1
.\]
The balanced (as well as Schr\"odinger-type) form 
of the Hermite operator is
\begin{eqnarray*}
\e^{-\frac{z^2}{2}}\S_{\lambda}\e^{\frac{z^2}{2}}
&=&
\partial_z^2-z^2-2\lambda.
\end{eqnarray*}
Note that the symmetry $(z,\lambda)\mapsto(\i z,-\lambda)$ is obvious
in the balanced form. 

\ber
The balanced form of the Hermite equation is known in the literature
as the {\em Weber} or {\em parabolic cylinder equation}.
It is usually written in one of two forms
\begin{eqnarray*}
\partial_z^2-\frac14z^2-k,&&\partial_z^2+\frac14z^2-k.
\end{eqnarray*}
\eer

\subsection{Even solution}

Inserting a power series in the equation we see that the Hermite equation has
an even solution
\begin{eqnarray*}
S_\lambda^+(z)&:=&\sum\limits_{j=0}^\infty\frac{(\frac{{a}}{2})_j}{(2j)!}(2z)^{2j}
\\=
\ 
F\Big(\frac{{a}}{2};\12;z^2\Big)&=&F_{-\frac12,\lambda}(z^2).
\end{eqnarray*}

It is the unique solution
satisfying
\beq S_\lambda^+(0)=1,\ \ \ \ \frac{\d}{\d z}S_\lambda^+(0)=0.\label{inh1}\eeq
It has the properties
\[S_\lambda ^+(z)=S_\lambda ^+(-z)
=\e^{z^2}S_{-\lambda }^+(\i z)
.\]

\subsection{Odd solution}

The Hermite equation has
an odd solution
\begin{eqnarray*}
S_\lambda^-(z)&:=&
\sum\limits_{j=0}^\infty\frac{(\frac{{a}+1}{2})_j}{(2j+1)!}(2z)^{2j+1},\\
=\ 2zF\Big(\frac{{a}+1}{2};\frac32;z^2\Big)&=&2zF_{\frac12,\lambda}(z^2)
.\end{eqnarray*}

It is the unique solution of the Hermite
equation satisfying
\beq S_\lambda^-(0)=0,\ \ \ \ \frac{\d}{\d z}S_\lambda^-(0)=2.\label{inh2}\eeq
It has the properties
\[S^-_\lambda (z)=-S_\lambda ^-(-z)
=-\i \e^{z^2}S_{-\lambda }^-(\i z)
.\]

\subsection{Standard solutions}

The Hermite equation has only one singular point, $\infty$. We will
see that one can define two kinds of solutions with a simple asymptotics at $\infty$.

By Thm \ref{dad12}, for appropriate $\gamma_1$ and $\gamma_2$ the following integrals are solutions:
\begin{eqnarray*}
\int\limits_{\gamma_1}\e^{-t^2-2tz}t^{\lambda -\frac12}\d t,
&&\\
\int\limits_{\gamma_2}\e^{t^2}
(z-t)^{-\lambda -\frac12}
\d
  t. 
&&
\end{eqnarray*}

In the first case the integrand has a singular point at $0$ and goes to zero as $t\to \pm \infty$. We can thus use $\gamma_1$ with such endpoints. We will see that they give all standard solutions.

In the second case the integrand has a singular point at $z$ and goes
to zero as $t\to \pm \i\infty$. Using $\gamma_2$ with such
endpoints we will also obtain all standard solutions.

\subsubsection{Solution $\sim z^{-{a}}$ for $z\to+\infty$} 

The following function is the solution of the Hermite equation that behaves as
$z^{-{a}}=z^{-\lambda-\frac12}$ for $|z|\to\infty$, $|\arg z|<\pi/2-\epsilon$:
\begin{eqnarray*}
S_\lambda (z)&:=&z^{-\lambda -\frac12}\tilde F_{-\frac12,\lambda }(-z^{-2})
\ = \
z^{-{a}}F\Big(\frac{{a}}{2},\frac{{a}+1}{2};-;-z^{-2}\Big).
\end{eqnarray*}
We will also introduce
 alternatively normalized solutions:
\begin{eqnarray*}
{\bf S}_\lambda^\I (z)&:=&2^{-\lambda-\frac12}\Gamma\Big(\lambda+\frac12\Big)S_\lambda(z)\\
& = &
2^{-a}z^{-{a}}\frac{1}{\Gamma(a)}F\Big(\frac{{a}}{2},\frac{{a}+1}{2};-;-z^{-2}\Big),\\
{\bf S}_\lambda^\0 (z)&:=&\sqrt{\pi}S_\lambda(z).
\end{eqnarray*}
(The normalization of ${\bf S}_\lambda^\0$ is somewhat trivial -- we introduce it to preserve the analogy with the Gegenbauer equation, which had a less trivially normalized solution
 ${\bf S}_{\alpha,\lambda}^\0$.)

Assuming that $z\not\in]-\infty,0]$,
 we have an integral representation valid for $ -\frac12<\lambda$:
\begin{eqnarray*}
\int\limits_0^\infty\e^{-t^2-2tz}t^{\lambda -\frac12}\d t
&=&{\bf S}_\lambda ^\I(z),
\end{eqnarray*} and for all parameters:
\begin{eqnarray*}
-\i\int\limits_{]-\i\infty,z^-,\i\infty[}\e^{t^2}
(z-t)^{-\lambda -\frac12}
\d
  t 
&=&{\bf S}^\0_\lambda (z).
\end{eqnarray*}

\subsubsection{Solution $\sim(-\i z)^{{a}-1}\e^{z^2}$ for 
$z\to+\i\infty$}
The following function is the solution of the Hermite equation that behaves as
$(-\i z)^{{a}-1}\e^{z^2}=(-\i z)^{\lambda -\frac12}\e^{z^2}$ for $|z|\to\infty$,  $|\arg z-\pi/2|<\pi/2-\epsilon$:
\begin{eqnarray*}
\e^{z^2}S_{-\lambda }(-\i z)&=&
(-\i z)^{\lambda -\frac12}\e^{z^2}\tilde F_{-\frac12,-\lambda}(z^{-2}).
\end{eqnarray*}

Assuming that $z\not\in[0,\infty[$, we have an
integral representation valid for all parameters:
\begin{eqnarray*}
\int\limits_{]-\infty,0^+,\infty[}\e^{-t^2-2tz}(-\i t)^{\lambda -\frac12}\d t&=&
\e^{z^2}{\bf S}_{-\lambda }^\0(-\i z),
\end{eqnarray*} and for $ \lambda<\frac12$:
\begin{eqnarray*}
-\i\int\limits_{[z,\i\infty[}\e^{t^2}(-\i (t-z))^{-\lambda -\frac12}\d t
&=& \e^{z^2}{\bf S}_{-\lambda }^\I(-\i z).
\end{eqnarray*}

\subsection{Connection formulas}

We can decompose the standard solutions into
the even and odd solutions:
\begin{eqnarray*}
S_\lambda (z)&=&\frac{\sqrt\pi}{\Gamma(\frac{2\lambda +3}{4})}S_\lambda ^+(z)
-\frac{\sqrt\pi}{\Gamma(\frac{2\lambda +1}{4})}S_\lambda ^-(z);
\\
\e^{z^2}S_{-\lambda }(-\i z)
&=&\frac{\sqrt{\pi}}{\Gamma(\frac{3-2\lambda }{4})}S_\lambda ^+(z)
+\i\frac{\sqrt{\pi}}{\Gamma(\frac{1-2\lambda }{4})}S_\lambda ^-(z).
\end{eqnarray*}

\subsection{Recurrence relations}

The following recurrence relations follow easily from the commutation properties of Subsect. \ref{symcom5a}:
\begin{eqnarray*}
\p_z S_\lambda(z)&=&-\Big(\frac12+\lambda \Big) S_{\lambda +1}(z),\\
(\p_z -2z)S_\lambda (z)&=&-2S_{\lambda -1}(z),\\[3ex]
(z\p_z+\frac12-\lambda -2z^2) S_\lambda (z)&=&-2
S_{\lambda -2}(z),\\
(z\p_z+\frac12+\lambda )S_\lambda (z)&=&-\frac12\Big(\frac12+\lambda \Big)
\Big(\frac32+\lambda\Big) S_{\lambda +2}(z).
\end{eqnarray*}

\subsection{Hermite polynomials}

If $-a=n=0,1,2,\dots$, then  Hermite functions are polynomials.

Following
Subsect. 
\ref{Hypergeometric type polynomials}, they can be defined by the following version of the Rodriguez-type formula:
\[
H_n(z):=\frac{(-1)^n}{n!}\e^{z^2}\p_z^n\e^{-z^2}.
\]

\begin{remark}
The Hermite polynomials found usually in the literature equal
\[n!H_n(z).\]
The advantage of  our convention is that the
 Rodriguez-type formula has the same form for all classes of
 hypergeometric type polynomials.
\end{remark}
The differential equation:
\begin{eqnarray*}
&&\big(\p_z^2-2z\p_z+2n\big)H_n(z)\\
&=&(-n;z,\p_z)H_n(z)=0.
\end{eqnarray*}
The generating function:
\[\begin{array}{l}\exp(2tz-t^2)=
\sum\limits_{n=0}^\infty
t^nH_n(z).
\end{array}\]
The integral representation:
\[\begin{array}{l}H_n(z)
=\frac{1}{2\pi\i}\int\limits_{[0^+]}\exp(2tz-t^2)t^{-n-1}\d t.
\end{array}\]
Recurrence relations:
\begin{eqnarray*}
\p_zH_n(z)&=&2H_{n-1}(z),\\
\left(\p_z-2z\right)H_n(z)&=&-(n+1)H_{n+1}(z),\\[5ex]
\left(z\p_z-n\right)H_n(z)&=&2H_{n-2}(z),\\
\left(z\p_z+n+1-2z^2\right)H_n(z)&=&-(n+1)(n+2)H_{n+2}(z).
\end{eqnarray*}

The differential equation, the Rodriguez-type formula, the  generating
function, the integral representation and the first pair of recurrence relations
are special cases of 
 the  corresponding formulas of Subsect.  
\ref{Hypergeometric type polynomials}.

We have several alternative expressions for Hermite polynomials:
\begin{eqnarray*}
H_n(z)&=&-\lim_{\nu\to n}(-1)^n(\nu{-}n){\bf S}_{-n-\frac12}^{\I}(z)\ =\ \frac{2^n}{n!}S_{-n-\frac12}(z)\\
&=&\frac{2^n}{n!}z^nF\Big(-\frac{n}{2},\frac{-n+1}{2};-;z^{-2}\Big)\\
&=&\sum\limits_{k=0}^{[\frac{n}{2}]}
\frac{(-1)^k(2z)^{n-2k}}{k!(n-2k)!}.
\end{eqnarray*}

Behavior at $\infty$.
\[\begin{array}{l}
\lim\limits_{n\to\infty}\frac{H_n(z)}{z^n}=\frac{2^n}{n!}.\end{array}\]

Initial conditions at $0$.
\[\begin{array}{l}
H_{2m}(0)=\frac{(-1)^m}{m!},\ \ \ \ H'_{2m}(0)=0,\\[3mm]
H_{2m+1}(0)=0,\ \ \ \ H'_{2m+1}(0)=\frac{(-1)^m2}{m!}.\end{array}\]

Identities for even and odd  polynomials.
\begin{eqnarray*}
H_{2m}(z)
&=&\frac{(-1)^m2^{2m}m!}{(2m)!}L_m^{-1\slash2}(z^2)
\
=\ \frac{(-1)^m(2z)^{2m}m!}{(2m)!}B_m^{-2m-\12}(-z^{-2}),\\
&=&\frac{(-1)^m}{m!}S_{-2m-\frac12}^+(z)=\frac{(-1)^m}{m!}F\Big(-m;\12;z^2\Big),\\
H_{2m+1}(z)
&=&\frac{(-1)^m2^{2m+1}m!}{(2m+1)!}zL_m^{1\slash2}(z^2)
\ =\ \frac{(-1)^m(2z)^{2m+1}m!}{(2m+1)!}B_m^{-2m-\frac{3}{2}}(-z^{-2})\\
& =&\ \frac{(-1)^m}{m!}S_{-2m-\frac32}^-(z)= \frac{(-1)^m}{m!}2zF\Big(-m;\frac{3}{2};z^2\Big).
\end{eqnarray*}

\appendix
\section{Contours for integral representations}

In this appendix we collect  contours used in various integral representations of hypergeometric type functions.

 For each basic type of integral representations considered in our text
we give at least one contour for every standard representation. 
We give the priority to type (a) contours. If they are unavailable, we show a type (b) contour. In some cases we present both a type (a) and type (b).

We also show contours that yield the degenerate solutions and the polynomial solutions. They are given by closed loops.

Here is the explanation of basic elements of our figures:

\begin{center}
\includegraphics{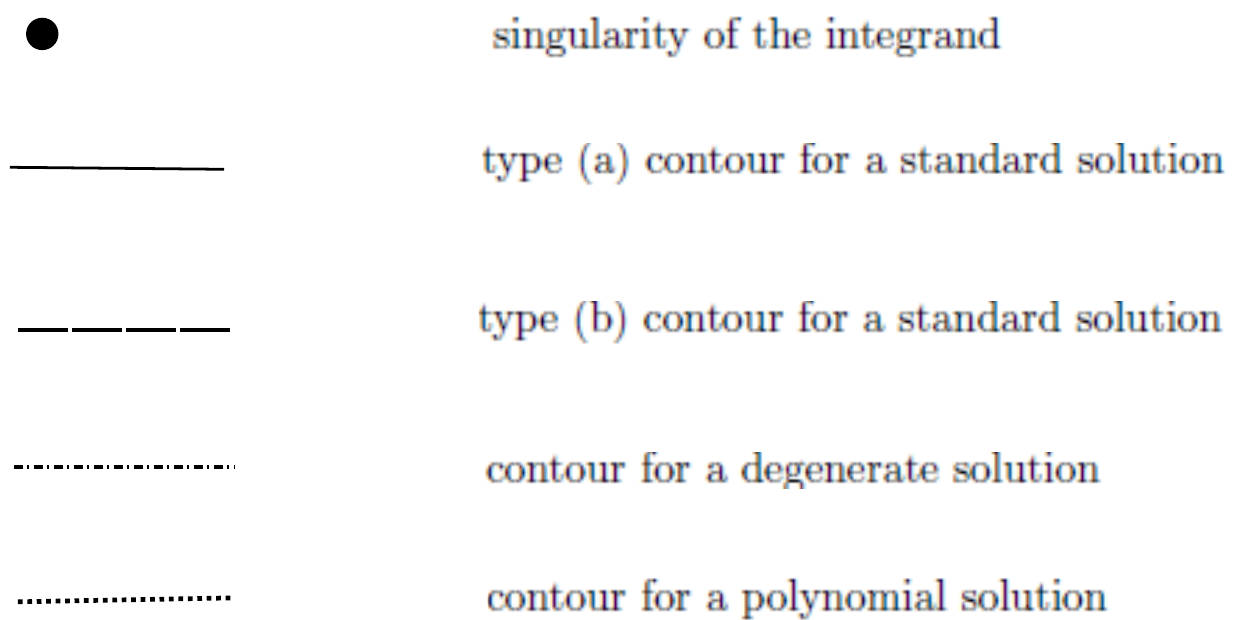}
\end{center}
\begin{center}
\includegraphics[width=16cm,totalheight=18cm]{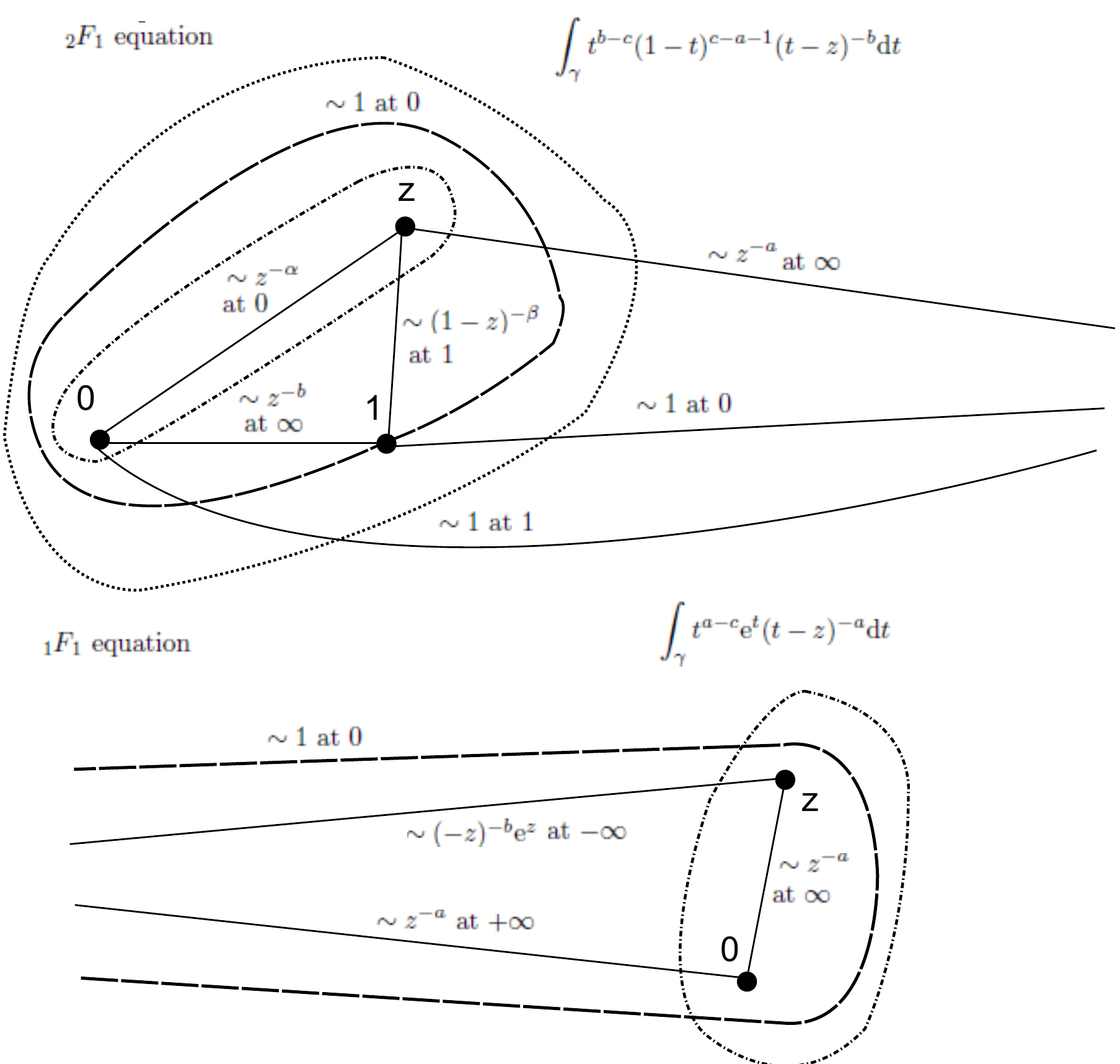}
\end{center}
\begin{center}
\includegraphics[width=16cm,totalheight=18cm]{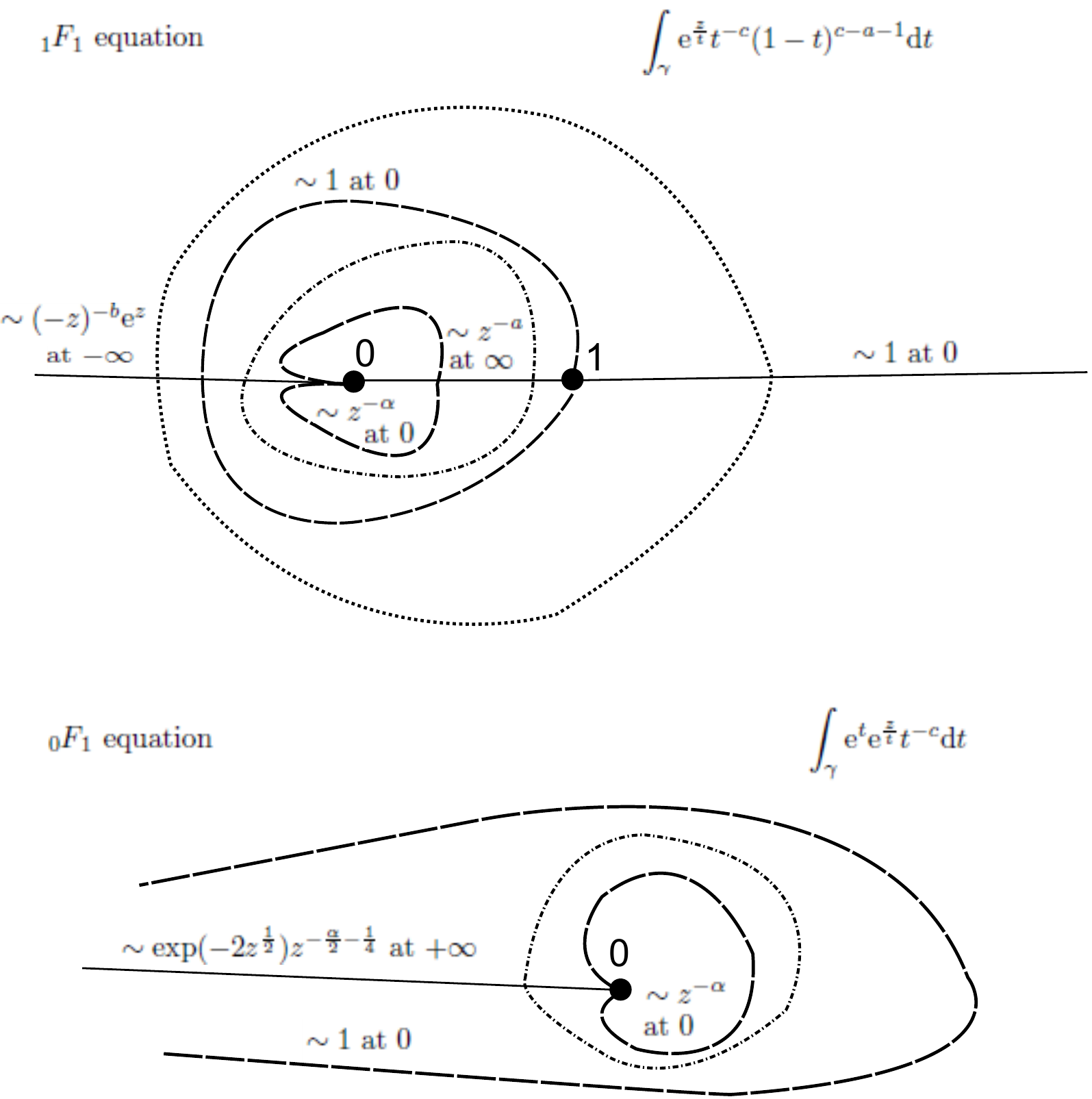}
\end{center}
\begin{center}
\includegraphics[width=16cm,totalheight=18cm]{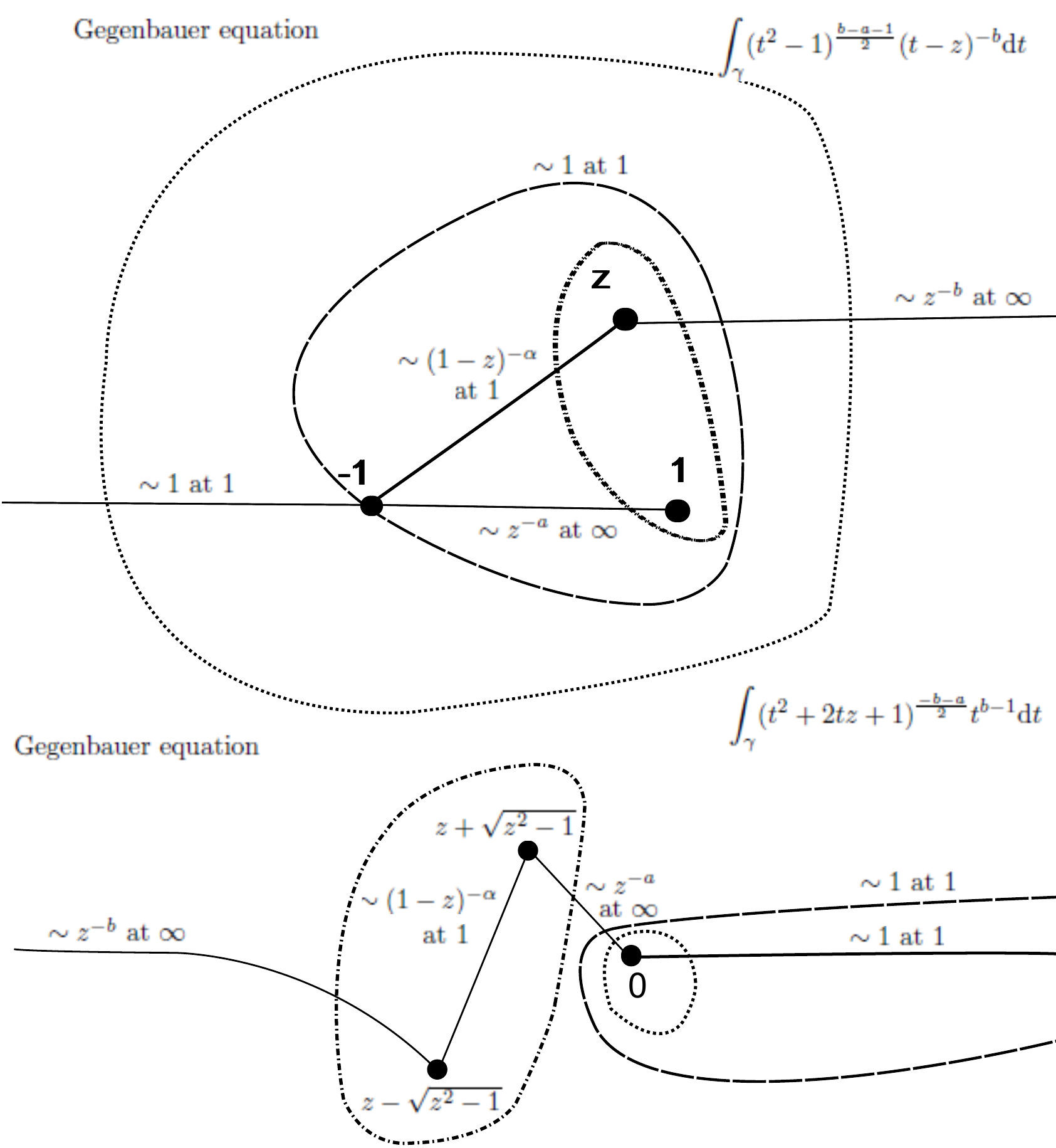}
\end{center}
\begin{center}
\includegraphics[width=16cm,totalheight=18cm]{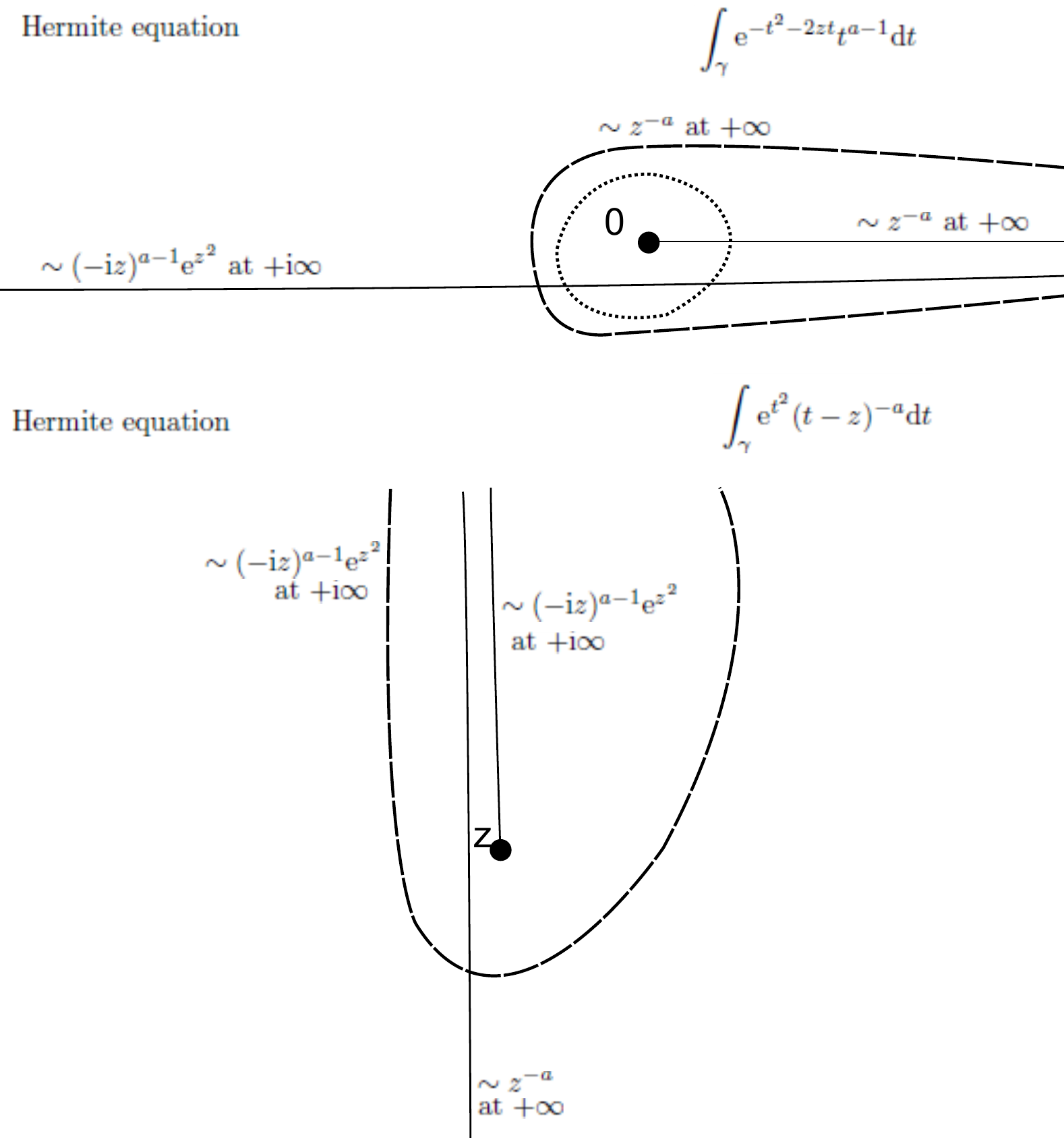}
\end{center}

\end{document}